\documentclass[11pt]{amsart}
\usepackage{pst-all}
\usepackage{amsmath}
\usepackage[latin1]{inputenc}
\usepackage{amsfonts}
\usepackage{amssymb}
\usepackage{multicol}
\usepackage{verbatim}
\usepackage{amsthm}
\usepackage{amscd}
\usepackage{pb-diagram}
\usepackage[mathscr]{eucal}
\usepackage[latin1]{inputenc}

\usepackage[T1]{fontenc}

\usepackage{amsmath}

\usepackage{amsfonts}

\usepackage{amssymb}

\usepackage{amsthm}

\usepackage{graphics}

\newcommand{\mk}{\medskip}

\newcommand{\ZZ}{\mathbb{Z}}

\newcommand{\CC}{\mathbb{C}}

\newcommand{\NN}{\mathbb{N}}

\newcommand{\yl}{15pt}
\newcommand{\yh}{8pt}

\newcommand{\ffbox}[1]{
\setbox9=\hbox{$\scriptstyle\overline{1}$}
\framebox[\yl][c]{\rule{0mm}{\ht9}${\scriptstyle #1}$}
}

\newcommand{\QQ}{\mathbb{Q}}

\def\dnl{=\hspace{-.2cm}>\hspace{-.2cm}=}
\def\dnr{=\hspace{-.2cm}<\hspace{-.2cm}=}
\def\sn{-\hspace{-.17cm}-}

\newcommand{\Glie}{\mathfrak{g}}

\newcommand{\Hlie}{\mathfrak{h}}

\newcommand{\demo}{\noindent {\it \small Proof:}\quad}

\newcommand{\U}{\mathcal{U}}

\newcommand{\Lo}{\mathcal{L}}

\newtheorem{thm}{Theorem}[section]

\newtheorem{defi}[thm]{Definition}

\newtheorem{cor}[thm]{Corollary}

\newtheorem{prop}[thm]{Proposition}

\newtheorem{lem}[thm]{Lemma}

\newtheorem{rem}[thm]{Remark}

\usepackage[all]{xy}
\title{Kirillov-Reshetikhin conjecture : the general case}

\author{David Hernandez}
\address{CNRS - UMR 8100 : Laboratoire de Math\'ematiques de Versailles, 45 avenue des Etats-Unis , Bat. Fermat, 78035 VERSAILLES, 
FRANCE}
\email{hernandez @ math . cnrs . fr}
\urladdr{http://www.math.uvsq.fr/\textasciitilde hernandez}

\begin{document}

\begin{abstract} We prove the Kirillov-Reshetikhin (KR) conjecture in the general case : for all twisted quantum affine algebras we prove that the characters of KR modules solve the twisted $Q$-system \cite{hkott} and we get explicit formulas for the character of their tensor products (the untwisted case was treated in \cite{Nab, Nad, her06}).  The proof is uniform and provides several new developments for the representation theory of twisted quantum affine algebras, including twisted Frenkel-Reshetikhin $q$-characters (expected in \cite{Fre, Fre2}). We also prove the twisted $T$-system \cite{ks}.  As an application we get explicit formulas for the twisted $q$-characters of fundamental representations for all types, including the conjectural formulas \cite{r} for types $D_4^{(3)}$, $E_6^{(2)}$. We prove the conjectural formulas \cite{ks} for KR modules in types $A_n^{(2)}$ ($n\geq 2$) and $D_4^{(3)}$. Eventually our results imply the conjectural branching rules \cite{hkott} to the quantum subalgebra of finite type.

\vskip 4.5mm


\end{abstract}

\maketitle

\tableofcontents

\section{Introduction} 

Characters and Frenkel-Reshetikhin $q$-characters of finite dimensional representations of untwisted quantum affine algebras have attracted much attention in recent years (see for example \cite{c0, cm, dk, fl, Fre2, Fre, her06, Naams, Nab, nn1} and references therein). Twisted quantum affine algebras and their representations have also been intensively studied (see for example \cite{a, bn, Cha5, da, hn, jin, jm, kas, ns, kmoy, oss, s} and references therein), and many conjectures have been formulated for the character of their finite dimensional representations \cite{hkott, ks, r}, but many of them are still open.

Let $\U_q(\hat{\Glie})$ be a (untwisted or twisted) quantum affine algebra of rank $n$. The Kirillov-Reshetikhin 
modules form a certain infinite class of simple finite dimensional representations of $\U_q(\hat{\Glie})$. The main question answered in this paper is the following : what is the character of the Kirillov-Reshetikhin modules and of their tensor products for the quantum group of finite type $\U_q(\Glie)\subset\U_q(\hat{\Glie})$ ? This problem goes back to 1931 as Bethe \cite{be} solved it for certain modules of type $A_1$ in another language. The methods to solve physical models involved here are now known as ``Bethe  Ansatz''. In the 80s, in a series of fundamental and striking papers, Kirillov and Reshetikhin \cite{kr, ki1, 
ki2, ki3} solved the problem for type $A$ and proposed formulas for all untwisted types by analyzing the Bethe Ansatz. These papers became the starting point of an intense research. The description of these characters is called the Kirillov-Reshetikhin conjecture. 

For untwisted cases, there are many results for these conjectures and related problems (see \cite{knt} for an historic and a guide through the huge literature on this subject; to name a few of very interesting recent ones see \cite{kl, hkoty, k, c0}).
The Kirillov-Reshetikhin conjecture was proved in \cite{Nab, Nad} for untwisted simply-laces cases, and in \cite{her06} for general untwisted case with a different uniform algebraic approach. 

In this paper we prove the twisted cases and so we get a complete uniform proof of the conjecture for all quantum affine algebras. As the representation theory in twisted cases is far less understood than in the untwisted cases, more work is needed than in \cite{her06}, and in this paper we also prove new general results for finite dimensional representations of twisted quantum affine algebras. In particular we develop the theory of twisted Frenkel-Reshetikhin $q$-characters expected in \cite{Fre, Fre2}, and we prove additional results and conjectures :

\begin{itemize}
\item the characterization of the image of twisted $q$-characters,

\item the special property of Kirillov-Reshetikhin for all twisted types,

\item explicit formulas for twisted $q$-characters of fundamental representations for all types, including the Conjectural formulas \cite{r} in types $D_4^{(3)}$, $E_6^{(2)}$,

\item the Conjecture \cite{hkott} of twisted $Q$-systems for all types,

\item the Conjecture \cite{ks} of twisted $T$-system for all types,

\item the Conjecture \cite{ks} of explicit formulas for twisted $q$-characters of Kirillov-Reshetikhin modules for types $A_n^{(2)}$ ($n\geq 2$), $D_4^{(3)}$,

\item the Conjectural branching rules to the quantum subalgebra of finite type conjectured in \cite{hkott} for types $A_{2n}^{(2)}$, $A_{2n-1}^{(2)}$, $D_n^{(2)}$, $D_4^{(3)}$,

\item the Conjecture \cite[Section 6]{ks} that in type $D_4^{(3)}$ the conjectural explicit formulas satisfy the twisted $T$-system.

\end{itemize}

 Note that the set of solutions that we get for twisted $Q$-systems and twisted $T$-systems are Laurent polynomials with positive coefficients (as characters), and so is an example of the Laurent phenomena for these systems (see \cite{fz, fzd} for a description in another context).

Let us describe the conjectures and recall some previous results on this problem in more details :

The Kirillov-Reshetikhin modules are simple $l$-highest weight modules (notion analog to the notion of highest weight 
module adapted to the Drinfeld realization of quantum affine algebras \cite{Cha2, Cha5}). They are 
characterized by their Drinfeld polynomials $(P_j(u))_{1\leq j\leq n}$ (analogs of the highest weight). In types different than $A_{2n}^{(2)}$, they  are of the form  :
\begin{equation*}
\begin{split}
P_j(u) = 
\begin{cases}
(1-au)(1-aq_i^2u)\cdots (1-q_i^{2k-2}u)&\text{ if $j = i$,}
\\0&\text{ if $j\neq i$,}
\end{cases}
\end{split}
\end{equation*}
(where $q_i=q^{d_{\overline{i}}}$, see section \ref{un}). In types $A_{2n}^{(2)}$, they are of form :
\begin{equation*}
\begin{split}
P_j(u) = 
\begin{cases}
(1-au)(1-aq^2u)\cdots (1-q^{2k-2}u)&\text{ if $j = i$,}
\\0&\text{ if $j\neq i$.}
\end{cases}
\end{split}
\end{equation*}

For $k=1$ the Kirillov-Reshetikhin modules are called fundamental representations. The
Kirillov-Reshetikhin conjecture predicts the character of these modules 
and of their tensor products for 
the subalgebra $\U_q(\Glie)\subset\U_q(\hat{\Glie})$ ($\U_q(\Glie)$ is a quantum Kac-Moody algebra of finite type) : 
in \cite{kr} conjectural formulas were given for these characters (they were obtained by observation of the Bethe 
Ansatz related to solvable lattice models). A conjectural induction rule called $Q$-system was also given in \cite{kr} (the $Q$-system for untwisted exceptional types was given in \cite{ki3}, and twisted $Q$-systems in \cite{hkott} for twisted types). We will denote by $\mathcal{F}(\nu)$ the formulas for the 
characters (we use here the version of \cite{ki1, ki2, hkoty, hkott, knt}. The version of \cite{kr} is 
slightly different because the definition of binomial coefficients is a little changed, see remark 1.3 of \cite{knt}). The Kirillov-Reshetikhin conjecture can be stated in the following form (see \cite{knt}), proved in the present paper for all types : the characters of Kirillov-Reshetikhin modules solve the (twisted) $Q$-system and any of their tensors products are given by the formulas $\mathcal{F}(\nu)$.

In the proof of the Kirillov-Reshetikhin conjecture for simply-laced cases \cite{Nab, Nad} and for the general untwisted cases \cite{her06}, one crucial tool in the theory of Frenkel-Reshetikhin $q$-characters. For untwisted quantum affine algebras this theory was introduced in \cite{Fre} (see also \cite{kn}) as analogs of characters for quantum affine algebras adapted to the Drinfeld presentation, and then developed in \cite{Fre2, Nab}. 

Let us remind some points of the theory for untwisted quantum affine algebras. In \cite{kns, ks} functional relations called $T$-system were defined (they are motivated by the observation of transfer matrices in solvable lattice models). Transposed in the language of Frenkel-Reshetikhin $q$-characters, and motivated by the relations between 
$q$-characters and Bethe-Anzatz in \cite{Fre}, it was naturally conjectured that the $q$-characters of 
Kirillov-Reshetikhin modules solve the $T$-system \cite{kosy}. As in terms of usual characters the $T$-system becomes the $Q$-system, this conjecture implies the Kirillov-Reshetikhin conjecture.

\noindent In general no explicit formulas for $q$-characters of finite dimensional simple modules are known. However Frenkel and Mukhin \cite{Fre2} defined an algorithm to compute the $q$-character of a class of simple modules (satisfying the ``special'' property, term defined in \cite{Nab}, see below), and Nakajima \cite{Nab} defined an algorithm to compute the $q$-character of arbitrary simple modules in simply-laced cases. Although in general no explicit formula has been obtained from them (they are very complicated), in some cases they give useful informations.

\noindent In particular Nakajima \cite{Nab, Nad} made a remarkable advance by proving the Kirillov-Reshetikhin conjecture and $T$-systems in simply-laced cases : using the main result of \cite{Nab} (the algorithm) he proved that in simply-laced cases the Kirillov-Reshetikhin modules are special and noticed that this property is useful in the study of these modules (the algorithm of \cite{Nab} is drastically simplified in this situation). The main result of \cite{Nab} is based on the study of quiver varieties \cite{Naams} and is not known in non simply-laced cases (see the conjecture of \cite{her02}), and so Nakajima's proof can not be used for all types.

\noindent In \cite{her04} we gave a general uniform proof for all untwisted types of these conjectures (Kirillov-Reshetikhin conjecture, special property, $T$-systems). Our proof is purely algebraic (the results of \cite{Nab} are not used; see the introduction of \cite{her04} for the main ideas of the proof) and so extends uniformly to untwisted non simply-laced cases. 

Let us go back to general quantum affine algebras. In the present paper we propose a general uniform proof for all twisted types of these conjectures (Kirillov-Reshetikhin conjecture, special property, twisted $T$-systems of \cite{ks}) which extends the proof of \cite{her04}. 

\noindent For the proof the main technical tool is the theory of twisted Frenkel-Reshetikhin $q$-characters that we develop in the present paper.  The existence of an analog theory for the twisted case was expected in \cite{Fre, Fre2}. Many technical non straightforward modifications and new developments are given compared to \cite{her04}, in particular because a twisted quantum affine algebras may have "elementary" subalgebras of type $A_1^{(1)}$ {\it and} of type $A_2^{(2)}$ (see the introduction of \cite{da} for general comments on this point) : this is one of the crucial technical point for the results of this paper. For example a characterization of the image of the twisted $q$-character morphism is proved.

After having developed the theory of twisted $q$-characters, the proof of the main results of the present paper has two steps : the special property and then the use of a twisted analog of the Frenkel-Mukhin algorithm to prove the twisted $T$-system. In the proof we have to study in details the "elementary" type $A_2^{(2)}$ as the structure of the corresponding Kirillov-Reshetikhin modules is drastically different than in the untwisted "elementary" $A_1^{(1)}$ type. 

As an application we prove several conjectural explicit formulas listed above. Although we see many important differences between representation theory of untwisted quantum affine algebras and twisted quantum affine algebras, we also prove precise relations between the ($q$-)characters from both theory : we construct an isomorphism between the Grothendieck rings of finite dimensional representations preserving Kirillov-Reshetikhin modules. As for the untwisted case \cite{Fre}, the Grothendieck ring of finite dimensional representations is proved to be commutative.

For a geometric side, in analogy with the untwisted simply-laced cases, our result gives an explicit formula of what would be the Euler number of a ``quiver variety'' in twisted cases. Nakajima told to the author that he knows a geometric approach to the representation theory of twisted quantum affine algebras.

The proof of the twisted $T$-systems was announced in \cite{her06}. For related results for current algebras see \cite{cmk}, and see \cite{herl} for the case of general untwisted quantum affinizations.

After the first version of this paper had been written, the author was asked by Nakajima how he completed the proof of the Drinfeld-Serre relations, sketched in \cite{Dri2, jin}. At the time he wrote this paper the author could not complete it, but he noticed that these relations are not needed for the results of this paper : in the present paper we start with Drinfeld-Jimbo generators and we give a reference or we prove all relations that are used in the proofs. In particular we adapt the proof of the classification of simple finite dimensional representations of twisted quantum affine algebras in \cite{Cha5} without using these relations. 

Let us describe the organization of this paper. It can be divided in two parts of different nature, the first one where the general theory of twisted $q$-character is developed and the general uniform proof of twisted $T$-systems is given, and the second one where as an application we explain explicit formulas and additional conjectures that follow from our results. Let us explain the structure of the paper in more details :

In Section \ref{un} we give backgrounds on quantum Kac-Moody algebras, Drinfeld realization, finite dimensional representations of (twisted) quantum affine algebras. 
In Section \ref{tw} we develop the theory of twisted Frenkel-Reshetikhin $q$-characters. We characterize the image and study the restriction to finite quantum subalgebras.
In Section \ref{deux} we give twisted $Q$-systems and twisted $T$-systems. The type $A_2^{(2)}$ is studied in Section \ref{adde}.
The special property is proved in Section \ref{special}.
The end of the proof of the main theorem (the twisted $T$-system) is given in the 
Section \ref{proofth}. In Section \ref{expform} we establish a connection with untwisted types and we give explicit formulas for characters of arbitrary tensor product of Kirillov-Reshetikhin modules (the Kirillov-Reshetikhin conjecture). Then we enter the second part of the paper. The branching rules to finite quantum subalgebra \cite{hkott} that follow from our results are given in Section \ref{branch}.
In Section \ref{expkr} we give the explicit formulas for twisted $q$-characters of Kirillov-Reshetikhin modules in several types that follow from our results (many of them had been conjectured in other papers).
In Section \ref{expfund} we give explicit formulas for twisted $q$-characters of fundamental representations for all types.

{\bf Acknowledgments :} The author would like to thank H. Nakajima for useful comments and references, and M. Okado for interesting questions on this work.

\section{Quantum Kac-Moody algebras and their representations}\label{un}

\subsection{Cartan matrix and quantized Cartan matrix} 

We consider a
generalized Cartan matrix $C=(C_{i,j})_{1\leq i,j\leq n}$, i.e., \label{carmat} $C_{i,j}\in\ZZ$, $C_{i,i}=2$, $C_{i,j}\leq 0$ for $i\neq j$ and $C_{i,j}=0$ if and only if $C_{j,i}=0$. We set $I=\{1,\dots,n\}$ and $l=\text{rank}(C)$. In the following we suppose that $C$ is symmetrizable, that is to say that there is a matrix $D=\text{diag}(r_1,\dots,r_n)$ ($r_i\in\NN^*$)\label{ri} such that $B=DC$\label{symcar} is symmetric. 
We consider a realization $(\Hlie, \Pi, \Pi^{\vee})$ of
$C$ (see \cite{kac}): $\Hlie$ is a $2n-l$ dimensional $\QQ$-vector
space, $\Pi=\{\alpha_1,\dots,\alpha_n\}\subset \Hlie^*$ (set of the
simple roots) and
$\Pi^{\vee}=\{\alpha_1^{\vee},\dots,\alpha_n^{\vee}\}\subset \Hlie$
(set of simple coroots) are set so that
$\alpha_j(\alpha_i^{\vee})=C_{i,j}$ for $1\leq i,j\leq n$. Let
$\Lambda_1,\dots,\Lambda_n\in\Hlie^*$ (resp.\ the
$\Lambda_1^{\vee},\dots,\Lambda_n^{\vee}\in\Hlie$) be the fundamental
weights (resp.\ coweights) :
$\Lambda_i(\alpha_j^{\vee})=\alpha_i(\Lambda_j^{\vee})=\delta_{i,j}$.

Let 
$$P =\{\lambda \in\Hlie^* \mid \text{$\lambda(\alpha_i^{\vee})\in\ZZ$ for all $i\in I$}\},$$  
$$P^+=\{\lambda \in P \mid \text{$\lambda(\alpha_i^{\vee})\geq 0$ for all $i\in I$}\},$$ 
be respectively the weight lattice and the semigroup of dominant weights.
Let $Q=\bigoplus_{i\in I} \ZZ \alpha_i\subset P$ (the root lattice) and $Q^+=\sum_{i\in I}\NN \alpha_i\subset Q$. For $\lambda,\mu\in \Hlie^*$, write $\lambda \geq \mu$ if $\lambda-\mu\in Q^+$. $\Delta$ is the set of roots and $\Delta^+$ is the set of positive roots.  Let $\nu:\Hlie^*\rightarrow \Hlie$ linear such that for all $i\in I$ we have $\nu(\alpha_i)=r_i\alpha_i^{\vee}$. For $\lambda,\mu\in\Hlie^*$, $\lambda(\nu(\mu))=\mu(\nu(\lambda))$.

The matrix $C$ is said to be of finite type if all its principal minors are positive. The matrix $C$ is said to be of affine type if all its proper principal minors are positive and $\text{det}(C) = 0$. The indecomposable affine Cartan matrices can be divided into two classes, the untwisted affines and the twisted affines (see \cite{kac} for the definition; the list of twisted affine Cartan matrices will be reminded in Section \ref{detw}). They give rise to algebras and representation theories which are analog but of different nature for many aspects, as for example we will see in the present paper.

\noindent In the following we suppose that $q\in\CC^*$ is not a root of unity. We choose $h\in\CC$ such that $q = \text{exp}(h)$ so that $q^r$ makes sense for $r\in\QQ$. For $l\in\ZZ, r\geq 0, m\geq m'\geq 0$ we define in $\ZZ[q^{\pm}]$ :
$$[l]_q=\frac{q^l-q^{-l}}{q-q^{-1}}\text{ , }[r]_q!=[r]_q[r-1]_q...[1]_q\text{ ,
}\begin{bmatrix}m\\m'\end{bmatrix}_q=\frac{[m]_q!}{[m-m']_q![m']_q!}.$$
Note that we can also define $[l]_q$ for $l\in\QQ$ (but then $[l]_q$ is no more in $\ZZ[q^{\pm 1}]$).

\subsection{Quantum Kac-Moody algebras}

$\Glie$ is a Kac-Moody Lie algebra of Cartan matrix $C$.

\begin{defi}\label{dkm} The quantum Kac-Moody algebra $\U_q(\Glie)$ is the $\CC$-algebra with generators $k_i^{\pm 1}$, $x_i^{\pm}$ ($i\in I$) and 
relations: 
$$k_ik_j=k_jk_i\text{ , } k_ix_j^{\pm}=q^{\pm r_i C_{i,j}}x_j^{\pm}k_i,$$
\begin{equation}\label{aidcont}[x_i^+,x_j^-]=\delta_{i,j}\frac{k_i-k_i^{-1}}{q^{r_i}-q^{-r_i}},\end{equation}
\begin{equation}\label{aidcontd}\underset{r=0... 1-C_{i,j}}{\sum}(-1)^r\begin{bmatrix}1-C_{i,j}\\r\end{bmatrix}_{q^{r_i}}(x_i^{\pm})^{1-C_{i,j}-r}x_j^{\pm}(x_i^{\pm})^r=0 \text{ (for $i\neq j$)}.\end{equation}
\end{defi}

\noindent This algebra was introduced independently by Drinfeld \cite{Dri1} and Jimbo \cite{jim}. The last relations are called quantum Serre relations and can be rewritten in the form
$$\underset{r=0... 1-C_{i,j}}{\sum}(-1)^r(x_i^{\pm})^{(1-C_{i,j}-r)}x_j^{\pm}(x_i^{\pm})^{(r)}=0 \text{ (for $i\neq j$)},$$
where we denote $(x_i^{\pm})^{(r)} = (x_i^{\pm})^r/[r]_{q^{r_i}}!$ for $r\geq 0$ (we will also use the notation $(x_i^\pm)^{(r)} = 0$ for $r < 0$).

It is remarkable that 
one can define a Hopf algebra structure on $\U_q(\Glie)$ by : 
$$\Delta(k_i)=k_i\otimes k_i,$$ 
$$\Delta(x_i^+)=x_i^+\otimes 1 + k_i\otimes x_i^+\text{ , }\Delta(x_i^-)=x_i^-\otimes 
k_i^{-1} + 1\otimes x_i^-,$$ 
$$S(k_i)=k_i^{-1}\text{ , }S(x_i^+)=-x_i^+k_i^{-1}\text{ , }S(x_i^-)=-k_ix_i^-,$$ 
$$\epsilon(k_i)=1\text{ , }\epsilon(x_i^+)=\epsilon(x_i^-)=0.$$
The quantum Kac-Moody algebras corresponding to Cartan matrices of finite (resp. of untwisted affine, twisted affine) types are called finite (resp. untwisted affine, twisted affine) quantum algebras.

\noindent Let $\U_q(\Hlie)$ be the commutative subalgebra of $\U_q(\Glie)$ generated by the $k_i^{\pm 1}$ ($i\in I$).

\noindent For $V$ a $\U_q(\Hlie)$-module and $\omega\in P$ we denote by $V_{\omega}$ the weight space of weight 
$\omega$: 
$$V_{\omega}=\{v\in V|\forall i\in I, k_i.v=q^{r_i\omega(\alpha_i^{\vee})}v\}.$$ 
In particular we have $x_i^{\pm}.V_{\omega}\subset V_{\omega \pm \alpha_i}$.

\noindent We say that $V$ is $\U_q(\Hlie)$-diagonalizable if $V=\underset{\omega\in P}{\bigoplus}V_{\omega}$ (in particular $V$ is of "type $1$").

\noindent For $V$ a finite dimensional $\U_q(\Hlie)$-diagonalizable module we define the usual character : 
$$\chi(V)=\underset{\omega\in P}{\sum}\text{dim}(V_{\omega})e^{\omega}\in\ZZ[e^{\omega}]_{\omega\in P}.$$

\subsection{Untwisted quantum loop algebras}\label{unt} Let be $\Glie$ a simple finite-dimensional Lie algebra of rank $n$. To $\Glie$ is associated an untwisted affine Kac-Moody algebra $\hat{\Glie}$ of rank $n+1$ and its quotient the loop algebra $\Lo\Glie$. The nodes of the Dynkin diagram of $\hat{\Glie}$ are indexed by $\{0,1,\cdots,n\}$ where $\{1,\cdots,n\}$ correspond to the Dynkin diagram of $\Glie$. There is a unique $c\in \sum_{0\leq i\leq n}\NN \alpha_i^{\vee}$ and $\delta\in Q^+$ such that 
$$\{h\in\sum_{0\leq i\leq n}\ZZ \alpha_i^{\vee} | \forall 0\leq i\leq n, \alpha_i(h) = 0\} = \ZZ c,$$
$$\{h\in Q|\forall 0\leq i\leq n,h(\alpha_i^\vee) = 0\} = \ZZ \delta.$$
We consider the $a_i^\vee,a_i\geq 0$ defined uniquely by $c = \sum_{0\leq i\leq n}a_i^\vee \alpha_i^{\vee}$ and $\delta = \sum_{0\leq i\leq n}a_i\alpha_i$. Then there is $p\in\NN$ such that for any $0\leq i\leq n$, $p a_i^\vee = r_i a_i$.

Let us give a quick review on the untwisted quantum loop algebras $\U_q(\Lo\Glie)$. The element $\prod_{0\leq i\leq n} k_i^{p a_i^{\vee}/r_i} = \prod_{0\leq i\leq n}k_i^{a_i}$ is central in $\U_q(\hat{\Glie})$. $\U_q(\Lo\Glie)$ is defined as a quotient of the untwisted quantum affine algebra $\U_q(\hat{\Glie})$ by the relation 
$$\prod_{0\leq i\leq n} k_i^{a_i} = 1.$$ 
It is clear that $\U_q(\Lo\Glie)$ inherits a Hopf algebra structure from $\U_q(\hat{\Glie})$. We use this algebra as all finite dimensional representations of the quantum affine algebra $\U_q(\hat{\Glie})$ can be factorized through $\U_q(\Lo\Glie)$ (see \cite{Cha2}). Let us give the Drinfeld presentation of the quantum loop algebra $\U_q(\Lo\Glie)$ \cite{Dri2, bec, jin}:

\begin{thm} $\U_q(\Lo\Glie)$ is the algebra with
generators $x_{i,r}^{\pm}$ ($i\in I, r\in\ZZ$), $k_i^{\pm 1}$ ($i\in I$), $h_{i,m}$ ($i\in I, m\in\ZZ-\{0\}$) and the
following relations ($i,j\in I, r,r'\in\ZZ, m\in\ZZ-\{0\}$): 
$$[k_i,k_j]=[k_{h},h_{j,m}]=[h_{i,m},h_{j,m'}]=0,$$
$$k_ix_{j,r}^{\pm}=q^{\pm r_iC_{i,j}}x_{j,r}^{\pm}k_i,$$
$$[h_{i,m},x_{j,r}^{\pm}]=\pm \frac{1}{m}[mB_{i,j}]_q x_{j,m+r}^{\pm},$$
$$[x_{i,r}^+,x_{j,r'}^-]=\delta_{i,j}\frac{\phi^+_{i,r+r'}-\phi^-_{i,r+r'}}{q^{r_i}-q^{-r_i}},$$
$$x_{i,r+1}^{\pm}x_{j,r'}^{\pm}-q^{\pm B_{i,j}}x_{j,r'}^{\pm}x_{i,r+1}^{\pm}=q^{\pm B_{i,j}}x_{i,r}^{\pm}x_{j,r'+1}^{\pm}-x_{j,r'+1}^{\pm}x_{i,r}^{\pm},$$
$$\underset{\pi\in\Sigma_s}{\sum}\underset{k=0..s}{\sum}(-1)^k\begin{bmatrix}s\\k\end{bmatrix}_{q^{r_i}}x_{i,r_{\pi(1)}}^{\pm}...x_{i,r_{\pi(k)}}^{\pm}x_{j,r'}^{\pm}x_{i,r_{\pi(k+1)}}^{\pm}...x_{i,r_{\pi(s)}}^{\pm}=0,$$
where the last relation holds for all $i\neq j$, $s=1-C_{i,j}$, all sequences of integers $r_1,...,r_s$. $\Sigma_s$ is
the symmetric group on $s$ letters. For $i\in I$ and $m\in\ZZ$, $\phi_{i,m}^{\pm}\in \U_q(\Lo\Glie)$ is determined
by the formal power series in $\U_q(\Lo\Glie)[[z]]$ (resp. in $\U_q(\Lo\Glie)[[z^{-1}]]$): 
$$\underset{m\geq 0}{\sum}\phi_{i,\pm m}^{\pm}z^{\pm m}=k_i^{\pm 1}\text{exp}(\pm(q^{r_i} - q^{-r_i})\underset{m'\geq
1}{\sum}h_{i,\pm m'}z^{\pm m'})$$ 
and $\phi_{i,m}^+=0$ for $m<0$, $\phi_{i,m}^-=0$ for $m>0$. \end{thm}

\noindent One has a triangular decomposition $$\U_q(\Lo\Glie)\simeq\U_q^+(\Lo\Glie)\otimes\U_q(\Lo\Hlie)\otimes\U_q^-(\Lo\Glie)$$ 
where $\U_q^\pm(\Lo\Glie)$ (resp. $\U_q(\Lo\Hlie)$) is the subalgebra generated by the $x_{j,r}^{\pm}$ (resp. the $k_i^{\pm 1}$ and the $h_{i,m}$).

A $\U_q(\Lo\Glie)$-module $V$ is said to be of $l$-highest weight if there is $v\in V$ eigenvector of all $\phi_{i,m}^{\pm}$ such that $V=\U_q(\Lo\Glie)^-.v$ and $\forall i\in I, m\in\ZZ, x_{i,m}^+.v=0$. For such a vector, denote by $\gamma = (\gamma_{i,\pm m}^{\pm})_{i\in I, m\geq 0}$ the corresponding eigenvalues of the $\phi_{i,\pm m}^{\pm}$. Note that we have necessarily $\gamma_{i,0}^+\gamma_{i,0}^- = 1$. Moreover if this condition is satisfied, there is a unique non trivial simple module $L(\gamma)$ of $l$-highest weight $\gamma$.

\begin{thm}\label{finu}\cite{Cha2}
The dimension of $L(\gamma)$ is finite if and only if there exist polynomials $(P_i)_{i\in I}$, $P_i(u)\in\CC[u]$, such that $P_i(0) = 1$ and $\gamma$ satisfies in $\CC[[u]]$ (resp. in $\CC[[u^{-1}]]$) :
$$\underset{m\geq 0}{\sum}\gamma_{i,\pm m}^{\pm}u^{\pm 
m}=q^{r_i\text{deg}(P_i)}\frac{P_i(uq^{-r_i})}{P_i(uq^{r_i})}.$$
Moreover all (type $1$) simple finite dimensional representations of $\U_q(\Lo\Glie)$ are of this form.
\end{thm}

For example for $i\in I, a\in\CC^*,k\geq 1$ the simple module corresponding to 
\begin{equation*}
\begin{split}
P_j(u) = \begin{cases} (1-ua)(1-uaq^{2r_i})\cdots (1-uaq^{2(k-1)r_i})&\text{ for $j = i$,}
\\1&\text{ for $j\neq i$,} 
\end{cases}
\end{split}
\end{equation*}
is called a Kirillov-Reshetikhin module.

Let us define the Frenkel-Reshetikhin $q$-characters morphism $\chi_q$ \cite{Fre} (see also \cite{kn}). Let $V$ be a finite dimensional representation of $\U_q(\Lo\Glie^\sigma)$ and $\gamma$ such that 
$$V_{\gamma} = \{v\in V|\exists p\geq 0,\forall i\in I,\forall m\geq 0,(\phi_{i,\pm m}^\pm -\gamma_{i,\pm m}^\pm)^p = 0\}\neq\{0\}.$$ 
Then there exist \cite{Fre} polynomials $(P_i)_{i\in I}, (Q_i)_{i\in I}$, $P_i(u),Q_i(u)\in\CC[u]$, such that $P_i(0) = Q_i(0) = 1$ and $\gamma$ satisfies in $\CC[[u]]$ (resp. in $\CC[[u^{-1}]]$) :
$$\gamma_i^{\pm}(u) = q^{r_i(\text{deg}(P_i) - \text{deg}(Q_i))}\frac{P_i(uq^{-r_i})Q_i(uq^{r_i})}{P_i(uq^{r_i})Q_i(uq^{-r_i})}.$$
Let $\text{Rep}(\U_q(\Lo\Glie))$ be the Grothendieck group of finite dimensional representations of $\U_q(\Lo\Glie)$. 

\begin{defi}\cite{Fre} The $q$-character morphism is the group morphism
$$\chi_q : \text{Rep}(\U_q(\Lo\Glie))\rightarrow \ZZ[Y_{i,a}^{\pm}]_{i\in I, a\in\CC^*}\text{ , }\chi_q(V)=\underset{\gamma}{\sum}\text{dim}(V_{\gamma})m_{\gamma},$$
where 
$$m_{\gamma}=\underset{i\in I, a\in\CC^*}{\prod}Y_{i,a}^{q_{i,a}-r_{i,a}}\text{ , }P_i(u)=\underset{a\in\CC^*}{\prod}(1-ua)^{q_{i,a}}\text{ , }Q_i(u)=\underset{a\in\CC^*}{\prod}(1-ua)^{r_{i,a}}.$$
\end{defi}

\begin{thm}\cite{Fre} 
$\chi_q$ is an injective ring morphism. In particular the Grothendieck ring $\text{Rep}(\U_q(\Lo\Glie))$ is commutative.
\end{thm}
For $i\in I$, $a\in\CC^*$ let 
\begin{equation*}
\begin{split}
A_{i,a} &= Y_{i,aq^{r_i}}Y_{i,aq^{-r_i}}\times \prod_{j|C_{j,i} = -1}Y_{j,a}^{-1}\times \prod_{j|C_{j,i}=-2}Y_{j,aq}^{-1}Y_{j,aq^{-1}}^{-1}
\\&\times\prod_{j|C_{j,i} = -3}Y_{j,aq^{-2}}^{-1}Y_{j,a}^{-1}Y_{j,aq^2}^{-1}.
\end{split}
\end{equation*}

\begin{thm}\label{imchiqu}\cite{Fre, Fre2} We have : 
$$\text{Im}(\chi_q) = \bigcap_{i\in I}(\ZZ[Y_{i,a}(1 + A_{i,aq^{r_i}}^{-1})]_{a\in\CC^*}\times \ZZ[Y_{j,a}]_{j\neq i,a\in\CC^*}).$$
\end{thm}

This result was obtained for non-simply quantum affine algebras in \cite{Naams} with a different method.

\subsection{Twisted quantum loop algebras}\label{detw} $\Glie$ is a simple finite dimensional Lie algebra as above. Consider an automorphism $\sigma$ of the Dynkin diagram, that is to say a bijection $\sigma : I\rightarrow I$ of the set $I$ of nodes of the Dynkin diagram of $\Glie$ such that $C_{\sigma(i),\sigma(j)} = C_{i,j}$ for all $i,j\in I$. Let $M$ be the order of $\sigma$. We study the twisted cases, that is to say we suppose that $M\geq 2$. From the classification of Dynkin diagram of finite type, we have $M\in\{2,3\}$ and $\Glie$ is simply-laced (in fact $\Glie$ is of type $A_n$ ($n\geq 2$), $D_n$ ($n\geq 4$), or $E_6$).

$I_{\sigma}$ is the set of orbits of $\sigma$. For $i\in I$ we denote by $\overline{i}\in I_{\sigma}$ the orbit of $i$. 

The corresponding twisted affine Lie algebra $\hat{\Glie}^{\sigma}$ is a Kac-Moody algebra with a Cartan matrix $(C_{i,j}^{\sigma})_{i,j\in \hat{I}_{\sigma}}$ of affine type, where $\hat{I}_{\sigma} = I_{\sigma}\sqcup \{\epsilon\}$ and $\epsilon$ is an additional node (see \cite{kac}) : for type not equal to $A_{2n}^{(2)}$ we have $\epsilon = 0$ and $I_\sigma = \{1,\cdots n\}$, and for type $A_{2n}^{(2)}$ we have $\epsilon = n$ and $I_\sigma = \{0,\cdots, n-1\}$. We denote by $\Glie^\sigma$ the semi-simple Lie algebra of Cartan matrix $(C_{i,j}^{\sigma})_{i,j\in I_{\sigma}}$ (it is of finite type, but not simply-laced except for type $A_2^{(2)}$).

The twisted quantum loop algebra $\U_q(\Lo\Glie^\sigma)$ is defined as a quotient of the twisted quantum affine algebra $\U_q(\hat{\Glie}^\sigma)$ by the relation $\prod_{i\in \hat{I}_\sigma} k_i^{a_i} = 1$ as for the untwisted quantum loop algebra. It is clear that $\U_q(\Lo\Glie^\sigma)$ inherits a Hopf algebra structure from $\U_q(\hat{\Glie}^\sigma)$. In the following we denote by $X_i^{\pm}$, $K_i$ ($0\leq i\leq n$) the Drinfeld-Jimbo generators of $\U_q(\Lo\Glie^\sigma)$.

The root vectors of $\U_q(\Lo\Glie^\sigma)$, and in particular the Drinfeld generators $x_{i,m}^{\pm}\in\U_q(\Lo\Glie^\sigma)$ for $i\in I$, $m\in\ZZ$ (resp. $h_{i,m}\in\U_q(\Lo\Glie^\sigma)$ for $m\in\ZZ$) corresponding to the roots $\pm \alpha_i + m \delta$ (resp. $m\delta$) are defined in \cite{da} by using Drinfeld-Jimbo generators and Lusztig automorphisms. 

Let us introduce some notations. We set $(d_0,\cdots,d_n)$ equal to:
\begin{equation*}
\begin{split}
(\frac{1}{2},2)&\text{ for type $A_{2}^{(2)}$,}
\\(\frac{1}{2},1,\cdots,1,2)&\text{ for type $A_{2n}^{(2)}$ ($n\geq 2$),}
\\(1,\cdots,1,2)&\text{ for type $A_{2n-1}^{(2)}$ ($n\geq 2$),} 
\\(1,2,\cdots,2,1)&\text{ for type $D_{n+1}^{(2)}$ ($n\geq 2$),}
\\(1,1,1,2,2)&\text{ for type $E_6^{(2)}$,}
\\(1,1,3)&\text{ for type $D_4^{(3)}$.}
\end{split}
\end{equation*}
Here we use the numbering of \cite{kac} given in the following diagrams.

\vspace{1cm}

{\centerline {\bf Twisted Affine Dynkin diagrams}}

\begin{multicols}{2}
\begin{itemize}

\item $A_2^{(2)}$ : {\large \vspace{-.55cm} 
$$\hspace{-1.4cm}\stackrel{0}{\circ}\hspace{-.08cm}=\hspace{-.17cm}<\hspace{-.21cm}= 
\hspace{-.08cm}\stackrel{1}{\circ}$$
\vspace{-1cm}
$$\hspace{-1.42cm}-\hspace{-.28cm}-\hspace{-.23cm}-$$
\vspace{-1.21cm}
$$\hspace{-1.42cm}-\hspace{-.28cm}-\hspace{-.23cm}-$$
}

\item $A_{2n}^{(2)}$ ($n\geq 2$) : {\large \vspace{-.15cm} 
$$\stackrel{0}{\circ}\hspace{-.08cm}\dnr\hspace{-.07cm}\stackrel{1}{\circ} \hspace{-.18cm}\sn\hspace{-.18cm}\stackrel{2}{\circ} 
\dots 
\stackrel{\text{n-1}}{\circ}\hspace{-.18cm}\dnr\hspace{-.07cm}\stackrel{\text{n}}{\circ}$$}

\item $A_{2n-1}^{(2)}$ ($n\geq 3$) : {\large \vspace{-.15cm} 
$$\stackrel{0}{\circ}\hspace{-.18cm}\sn\hspace{-.18cm}\stackrel{2}{\circ} \hspace{-.18cm}\sn\hspace{-.18cm}\stackrel{3}{\circ} 
\dots 
\stackrel{\text{n-1}}{\circ}\hspace{-.18cm}\dnr\hspace{-.07cm}\stackrel{\text{n}}{\circ}$$
\vspace{-.86cm}$$\hspace{-1.7cm}|$$
\vspace{-.78cm}$$\hspace{-1.38cm}\circ\text{ \scriptsize 1}$$}

\item $D_{n+1}^{(2)}$ ($n\geq 2$) : {\large \vspace{-.15cm} 
$$\stackrel{0}{\circ}\hspace{-.08cm}\dnr\hspace{-.07cm}\stackrel{1}{\circ} \hspace{-.18cm}\sn\hspace{-.18cm}\stackrel{2}{\circ} 
\dots 
\stackrel{\text{n-1}}{\circ}\hspace{-.18cm}\dnl\hspace{-.07cm}\stackrel{\text{n}}{\circ}$$}

\item $E_6^{(2)}$ : {\large \vspace{-.15cm} 
$$\hspace{-.4cm}\stackrel{0}{\circ}\hspace{-.18cm}\sn\hspace{-.18cm}\stackrel{1}{\circ}\hspace{-.18cm}\sn\hspace{-.18cm}\stackrel{2}{\circ} 
\hspace{-.08cm}\dnr\hspace{-.08cm}\stackrel{3}{\circ} 
\hspace{-.18cm}\sn\hspace{-.18cm}\stackrel{4}{\circ}$$}

\item $D_4^{(3)}$ : {\large \vspace{-.15cm} 
$$\hspace{-.4cm}\stackrel{0}{\circ}\hspace{-.18cm}\sn\hspace{-.18cm}\stackrel{1}{\circ}\hspace{-.12cm}\equiv\hspace{-.18cm}<\hspace{-.18cm}\equiv 
\hspace{-.1cm}\stackrel{2}{\circ}$$

}
\end{itemize}

\end{multicols}

We can choose for $i\in \hat{I}^\sigma$, $r_i = d_i$. Here we have an exception to the convention $r_i\in\NN$ as for type $A_{2n}^{(2)}$, $r_0 = 1/2$; this is not a problem as $q^{\frac{1}{2}}$ makes sense.

Remark : for $i\in I_{\sigma}$, we have 
\begin{equation*}
\begin{split}
C_{i,\sigma(i)} = 0 &\Rightarrow d_i = 1,
\\C_{i,\sigma(i)} = 2 &\Rightarrow d_i = M,
\\C_{i,\sigma(i)} = -1 &\Rightarrow d_i = 1/2.
\end{split}
\end{equation*}

Let $\omega\in\CC$ be a primitive $M^{th}$-root of $1$. We denote $q_i = q^{d_i}$ for $i\in \hat{I}_{\sigma}$.

For $i,j\in I$, define $d_{i,j}\in\QQ$ and $P_{i,j}^\pm(u_1,u_2)\in\QQ[u_1,u_2]$ by :

\noindent if $C_{i,\sigma(i)} =2$, then $d_{i,j}=\frac{1}{2}$ and $P_{i,j}^\pm(u_1,u_2)=1$,

\noindent if $C_{i,\sigma(i)}=0$ and $\sigma(j)\ne j$, then $d_{i,j}=\frac{1}{4M}$ and $P_{i,j}^\pm(u_1,u_2)=1$,

\noindent if $C_{i,\sigma(i)}=0$ and $\sigma(j)=j$, then $d_{i,j}=\frac{1}{2}$ and $P_{i,j}^\pm(u_1,u_2)=\frac{u_1^M q^{\pm 2M}-u_2^M}{u_1q^{\pm 2}-u_2}$,

\noindent if $C_{i,\sigma(i)}=-1$, then $d_{i,j}=\frac{1}{8}$ and $P_{i,j}^\pm(u_1,u_2)=u_1q^{\pm 1}+u_2$.

Note that $\U_q(\Lo\Glie^\sigma)$ has a $\ZZ$-grading defined by $\text{deg}(X_i^+) = \text{deg}(X_i^-) = \text{deg}(K_i) = 0$ for $i\in I_\sigma$, $\text{deg}(X_\epsilon^+) = 1$, $\text{deg}(X_\epsilon^-) = -1$, $\text{deg}(K_\epsilon^{\pm 1}) = 0$. Then we have $$\text{deg}(x_{i,r}^{\pm}) = r\text{ , }\text{deg}(h_{i,m}) = m\text{ , }\text{deg}(k_i^{\pm 1}) = 0.$$

\begin{defi}\label{phid} For $m\in\ZZ, i\in I$, $\phi_{i,m}^{\pm}\in \U_q(\Lo\Glie)$ is defined by the formal power series in $\U_q(\Lo\Glie^{\sigma})[[z]]$ (resp. in $\U_q(\Lo\Glie^{\sigma})[[z^{-1}]]$): 
$$\underset{m\geq 0}{\sum}\phi_{i,\pm m}^{\pm}z^{\pm m}=k_i^{\pm 1}\text{exp}(\pm(q_{\overline{i}}-q_{\overline{i}}^{-1})\underset{m'\geq
1}{\sum}h_{i,\pm m'}z^{\pm m'})$$ 
and $\phi_{i,m}^+=0$ for $m < 0$, $\phi_{i,m}^- = 0$ for $m > 0$.
\end{defi}

Let $N_+ = \sum_{i\in I, m\in\ZZ} \U_q(\Lo\Glie^\sigma) x_{i,m}^+$.

\begin{thm}\label{mult}
For $i\in I$, $k\geq 0$, $m\in\ZZ$ we have :
$$\Delta(\Psi_{i,\pm k}^{\pm}) = \sum_{k'=0\cdots k}\Psi_{i,\pm k'}^{\pm}\otimes \Psi_{i,\pm (k-k')}^{\pm}\text{ mod }(N_+ \otimes \U_q(\Lo\Glie^\sigma) + \U_q(\Lo\Glie^\sigma)\otimes N_+),$$
$$\Delta(x_{i,m}^+) \in (N_+ \otimes \U_q(\Lo\Glie^\sigma) + \U_q(\Lo\Glie^\sigma)\otimes N_+).$$
\end{thm}
\noindent For the case of $A_2^{(2)}$ see \cite{Cha5}, in general see \cite[Proposition 7.1.2]{da}, \cite[Proposition 7.1.5]{da} and \cite[Theorem 2.2]{jm}.

\noindent Example : let $\U_q^\tau = \U_q(\Lo sl_3^{\tau})$ where $\tau$ is the unique non trivial automorphism of the Dynkin diagram of $sl_3$. As there is only one orbit $|I_{\tau}|=1$, we can consider generators : $k^{\pm 1}, x_m^{\pm}, h_r$ (in the presentation of $\U_q^\tau$ in \cite{Cha5}, we have to replace $q$ by $q^{\frac{1}{2}}$ as we choose the normalization of \cite{jin}). As we will use it in the next sections, let us describe more explicitly the algebra $\U_q^\tau$. 

We have $C_{1,0}^\tau = -1$, $C_{0,1}^\tau = -4$, $q_1 = q^2$, $q_0 = q^{\frac{1}{2}}$. $\U_q^\tau$ is the algebra with generators $X_0^+$, $X_1^+$, $X_0^-$, $X_1^-$, $K_0^{\pm 1}$, $K_1^{\pm 1}$ and relations :
$$K_1 = K_0^{-2},$$
$$K_0X^+_1K_0^{-1} = q^{-2}X^+_1\text{ , }K_0X^+_0K_0^{-1} = qX^+_0,$$
$$K_0X^-_1K_0^{-1} = q^2X^-_1\text{ , }K_0X^-_0K_0^{-1} = q^{-1}X^-_0,$$
$$[X^+_1,X^-_0] = [X^+_0,X^-_1] = 0\text{ , }[X^+_0,X^-_0] = \frac{K_0 - K_0^{-1}}{q^{\frac{1}{2}} - q^{-\frac{1}{2}}}\text{ , }[X^+_1,X^-_1] = \frac{K_1 - K_1^{-1}}{q^2 - q^{-2}},$$
$$X_0(X_1)^{(2)} - X_1X_0X_1 + (X_1)^{(2)} X_0 = 0,$$
$$X_1(X_0)^{(5)} - X_0X_1X_0^{(4)} + X_0^{(2)}X_1X_0^{(3)} - X_0^{(3)}X_1X_0^{(2)} + X_0^{(4)}X_1X_0 - X_0^{(5)}X_1 = 0.$$
where $X = X^+$ or $X = X^-$.

Consider the Lusztig automorphisms $T_0$, $T_1$ of $\U_q^\tau$ defined by ($j\neq i$) :
$$T_i(X^+_i) = -X^-_iK_i\text{ , }T_i(X^-_i) = -K_i^{-1}X^+_i \text{ , }T_i(K_j) = K_i^{-C_{i,j}}K_j,$$
$$T_i(X^+_j) = \sum_{r = 0\cdots -C_{i,j}}(-1)^rq_i^{-r}(X^+_i)^{(-C_{i,j}-r)}X^+_j(X^+_i)^{(r)},$$
$$T_i(X^-_j) = \sum_{r = 0\cdots -C_{i,j}}(-1)^rq_i^r(X^-_i)^{(r)}X^-_j(X^-_i)^{(-C_{i,j}-r)}.$$
We have :
$$T_i^{-1}(X^+_i) = -K_i^{-1}X^-_i\text{ , }T_i^{-1}(X^-_i) = - X^+_iK_i \text{ , }T_i^{-1}(K_j) = K_i^{-C_{i,j}}K_j,$$
$$T_i^{-1}(X^+_j) = \sum_{r = 0\cdots -C_{i,j}}(-1)^rq_i^{-r}(X^+_i)^{(r)}X^+_j(X^+_i)^{(-C_{i,j}-r)},$$
$$T_i^{-1}(X^-_j) = \sum_{r = 0\cdots -C_{i,j}}(-1)^rq_i^r(X^-_i)^{(-C_{i,j}-r)}X^-_j(X^-_i)^{(r)}.$$
Let $T = T_1\circ T_0$. Then we have by definition :
$$x_n^+ = T^{-n}(X^+_0)\text{ , }x_n^- = T^n(X^-_0)\text{ , }k^{\pm 1} = K_0^{\pm 1}.$$

\begin{lem}\label{back} We have 
$$K_0 = k\text{ , }X^+_0 = x_0^+\text{ , }X^-_0 = x_0^-,$$
$$K_1 = k^{-2}\text{ , }X^+_1 = [4]_{q^{\frac{1}{2}}}^{-1}k^{-2}[x_0^-,x_1^-]_q\text{ , }X^-_1 = [4]_{q^{\frac{1}{2}}}^{-1}q[x_{-1}^+,x_0^+]_{q^{-1}}k^2.$$
\end{lem}

These relations are stated in \cite{a} (note that the Cartan matrix of type $A_2^{(2)}$ used in \cite{a} is the transposed of the one in \cite{k}). An an illustration let us write the proof as it is analog for all types :

\demo By definition we have $x_0^+ = X^+_0$ and 
$$x_{-1}^+ = T(X^+_0) = T_1T_0(X^+_0) = T_1(-X^-_0K_0) = - (X^-_0X^-_1 - q^2 X^-_1X^-_0)K_1K_0.$$
Moreover from the relations 
$$[X^+_0,X^-_0] = (K_0 - K_0^{-1})/(q^{1/2} - q^{-1/2}) = 0 \text{ and }[X^+_0,X^-_1] = 0$$ 
we get that $x_{-1}^+x_0^+ - q^{-1} x_0^+x_{-1}^+$ is equal to
\begin{equation*}
\begin{split}
 &- (X^-_0X^-_1 - q^2 X^-_1X^-_0)K_1K_0 X^+_0 + q^{-1} X_0^+(X^-_0X^-_1 - q^2 X^-_1X^-_0)K_1K_0
\\=& q^{-1}[4]_{q^{1/2}} X^-_1 K_0^2.
\end{split}
\end{equation*}
By definition we have 
$$x_1^- = -K_1^{-1}K_0^{-1}(X^+_1X^+_0 - q^{-2}X^+_0X^+_1)\text{ and }x_0^- = X^-_0.$$ 
So by using the relation $[X^+_0,X^-_0] = (K_0 - K_0^{-1})/(q^{1/2} - q^{-1/2})$ and $[X^-_0,X^+_1] = 0$ we get that $x_0^-x_1^- - q x_1^- x_0^-$ is equal to
\begin{equation*}
\begin{split}
 &-X^-_0 K_1^{-1}K_0^{-1}(X^+_1X^+_0 - q^{-2}X^+_0X^+_1) + q K_1^{-1}K_0^{-1}(X^+_1X^+_0 - q^{-2}X^+_0X^+_1) X^-_0 
\\= &K_0^2X^+_1 q[4]_{q^{1/2}}.
\end{split}
\end{equation*}
\qed

Let us go back to the general case. 

It is known that the Drinfeld generators generate the algebra $\U_q(\hat{\Glie}^\sigma)$. The proof is word by word the same as for the untwisted case which is given in \cite{bec} in the last remark in the proof of \cite[Theorem 4.7]{bec}. Let $W$ be the Weyl group of $\hat{\Glie}^\sigma$ and $W_0\subset W$ the Weyl group of $\Glie^\sigma$. The Lusztig automorphisms are defined for any element of $W$. For $i\in I_\sigma$, we have $X_i^{\pm}\in \CC^* x_{i,0}^{\pm}$, and so it suffices to check that $X_\epsilon^{\pm}$ can be expressed in terms of Drinfeld generators. Let $\theta = \delta - \alpha_{\epsilon}$. Let $i\neq \epsilon$ and $s_{\theta_i}\in W_0$ such that $s_{\theta_i}(\alpha_i) = \theta$ in the weight lattice of $\Glie^\sigma$. Then we have $X^+_{\epsilon}\in\CC[K_j^{\pm 1}]_{j\neq\epsilon}T_{\theta_i}(x_i^-(1))$ (see the proof of \cite[Theorem 4.7]{bec}; this is analog for $X^-_\epsilon$). We get explicit formulas for $X_\epsilon^\pm$, see the formulas in \cite{Dri2, jm}. For illustration of this result see Lemma \ref{back} for type $A_2^{(2)}$, and for other types :

Type $A_{2n}^{(2)}$ ($n\geq 2$) : $\Glie^\sigma$ is of type $B_n$ and
$$\theta = 2(\alpha_{n-1} + \alpha_{n-2} + \cdots + \alpha_0) = (s_{n-1}s_{n-2}\cdots s_1s_0s_1\cdots s_{n-2})(\alpha_{n-1}),$$

Type $A_{2n-1}^{(2)}$ ($n\geq 3$) : $\Glie^\sigma$ is of type $C_n$ and 
$$\theta = \alpha_1 + 2(\alpha_2+\cdots + \alpha_{n-1}) + \alpha_n = (s_2s_3\cdots s_{n-1}s_ns_{n-1}\cdots s_2)(\alpha_1),$$

Type $D_{n+1}^{(2)}$ ($n\geq 2$) : $\Glie^\sigma$ is of type $B_n$ and
$$\theta = \alpha_1 +\alpha_2+\cdots +\alpha_n = (s_1s_2\cdots s_{n-1})(\alpha_n),$$

Type $E_6^{(2)}$ : $\Glie^\sigma$ is of type $F_4$ and
$$\theta = 2\alpha_1 + 3\alpha_2 + 2\alpha_3 + \alpha_4 = (s_1 s_2 s_3 s_2 s_4 s_3 s_2)(\alpha_1),$$

Type $D_4^{(3)}$ : $\Glie^\sigma$ is of type $G_2$ and
$$\theta = 2\alpha_1 + \alpha_2 = (s_1s_2)(\alpha_1).$$

(Note that for type $A_{2n}^{(2)}$ the convention $\alpha_0 = \Lambda_0 - \Lambda_1$ is used). 

The following is proved in \cite{da} (see Section 3, 4, Lemma 5.1 and Theorem 5.3.2 in \cite{da}) :

\begin{thm}\cite{da} We have the relations ($i,j\in I, r,r'\in\ZZ, m\in\ZZ-\{0\}$) :
\begin{equation}\label{eqdun}x_{\sigma(i),r}^{\pm} = \omega^r x_{i,r}^{\pm}\text{ , }h_{\sigma(i),m} = \omega^mh_{i,m}\text{ , }k_{\sigma(i)}^{\pm 1} = k_i^{\pm 1},\end{equation}
\begin{equation}\label{eqddeux}[k_i,k_j]=[k_i,h_{j,m}]=[h_{i,m},h_{j,m'}]=0,\end{equation}
\begin{equation}\label{eqdtrois}k_ix_{j,r}^{\pm}=q^{\pm \sum_{k = 1,\cdots,M}C_{i,\sigma^k(j)}}x_{j,r}^{\pm}k_i,\end{equation}
\begin{equation}\label{eqdquatre}[h_{i,m},x_{j,r}^{\pm}]=\pm \frac{1}{m}(\sum_{k=1\cdots M}[mC_{i,\sigma^k(j)}/d_{\overline{i}}]_{q_{\overline{i}}}\omega^{mk})x_{j,m+r}^{\pm},\end{equation}
\begin{equation}\label{eqdcinq}[x_{i,r}^+,x_{j,r'}^-]=\sum_{k=1\cdots M}\delta_{\sigma^k(i),j}\omega^{kr'}\frac{\phi^+_{i,r+r'}-\phi^-_{i,r+r'}}{q_{\overline{i}}-q_{\overline{i}}^{-1}},\end{equation}
where the $\phi^\pm_{i,r}$ are given in Definition \ref{phid}.
\end{thm}

As a consequence we have a surjective map :
$$\U_q^+(\Lo\Glie^\sigma)\otimes \U_q(\Lo\Hlie^\sigma)\otimes \U_q^-(\Lo\Glie^\sigma)\rightarrow \U_q(\Lo\Glie^\sigma),$$
where $\U_q^\pm(\Lo\Glie^\sigma)$ (resp. $\U_q(\Lo\Hlie^\sigma)$) is the subalgebra generated by the $x_{i,m}^\pm$ (resp. the $h_{i,r}$ and the $k_i^{\pm 1}$). 

%Note that $\U_q(\Lo\Glie^\sigma)$ can be presented by the Drinfeld generators and the above relations with in addition  the relations (\ref{aidcont}), (\ref{aidcontd}) where the $X^\pm_i$ are expressed with Drinfeld generators (we get all relations of Definition \ref{dkm}).

The following Drinfeld-Serre relations are stated in \cite{Dri2, jin} :
\begin{equation}\label{eqdsix}(\prod_{k=1\cdots M} (u_1 - \omega^k q^{\pm C_{i,\sigma^k(j)}}u_2)) x_i^{\pm}(u_1)x_j^{\pm}(u_2)\end{equation}
\begin{equation*}= (\prod_{k=1\cdots M}(u_1q^{\pm C_{i,\sigma^k(j)}} - \omega^ku_2))x_j^{\pm}(u_2)x_i^{\pm}(u_1),\end{equation*}
if $C_{i,j}=-1$ and $\sigma(i)\neq j$ then (Sym denotes the symmetrization over $u_1,u_2$) :
\begin{equation}\label{eqdsept}\text{Sym}\{P_{ij}^\pm(u_1,u_2)(x_j^\pm(v)x_i^\pm(u_1)x_i^\pm(u_2)\end{equation} \begin{equation*}-(q^{2md_{ij}}+q^{-2md_{ij}})x_i^\pm(u_1)x_j^\pm(v)x_i^\pm(u_2)
+x_i^\pm(u_1)x_i^\pm(u_2)x_j^\pm(v))\}=0,\end{equation*}
if $C_{i,\sigma(i)}=-1$ then (Sym denotes the symmetrization over $u_1,u_2,u_3$) :
\begin{equation}\label{eqdhuit}\text{Sym}\{(q^{\frac{3}{2}}u_1^{\mp 1}-(q^{\frac{1}{2}}+q^{-\frac{1}{2}})u_2^{\mp 1}+q^{-\frac{3}{2}}u_3^{\mp 1})x_i^\pm(u_1)x_i^\pm(u_2)x_i^\pm(u_3)\}=0,\end{equation}
\begin{equation}\label{eqdneuf}\text{Sym}\{(q^{-\frac{3}{2}}u_1^{\pm 1}-(q^{\frac{1}{2}}+q^{-\frac{1}{2}})u_2^{\pm 1}+q^{\frac{3}{2}}u_3^{\pm 1})x_i^\pm(u_1)x_i^\pm(u_2)x_i^\pm(u_3)\}=0,\end{equation}
where for $i\in I$, $x_i^{\pm}(u) = \sum_{l\in\ZZ}x_{i,l}^{\pm}u^{-l}$.
At the time he wrote the paper, the author could not complete the proof of these relations sketched in \cite{Dri2, jin}, but he noticed that these additional relations are not needed for the results of this paper and so are not used. 

In the case of $\U_q^\tau$, the following relations between Drinfeld generators $k^{\pm 1}$, $h_r$, $x_m^{\pm}$ ($r\neq 0$, $m\in\ZZ$) are proved in \cite[Appendix B]{a} :
\begin{equation}\label{taucart} k k^{-1} = k^{-1} k =1\text{ , }[k,h_r] = [h_r,h_{r'}] = 0\text{ , }kx_m^{\pm} k^{-1} = q^{\pm 1}x_m^{\pm},\end{equation}
\begin{equation}\label{plusmoins}[x_r^+,x_{r'}^-] = \frac{\phi^+_{r + r'} - \phi^-_{r+ r'}}{q^{1/2} - q^{-1/2}},\end{equation}
\begin{equation}\label{carp}[h_r,x_m^\pm] = \pm \frac{q^r - q^{-r}}{r(q^{1/2}-q^{-1/2})} (q^r + q^{-r} + (-1)^{r+1})x_{r+m}^\pm,\end{equation}
\begin{equation}\label{int}x_{p}^{\pm}x_{m+2}^\pm + (q^{\mp 1} - q^{\pm 2})x_{m+1}^\pm x_{p+1}^\pm - q^{\pm 1}x_m^\pm x_{p+2}^\pm \end{equation}
\begin{equation*}= q^{\pm 1} x_p^\pm x_{m+2}^\pm + (q^{\pm 2} - q^{\mp 1})x_{p+1}^\pm x_{m+1}^\pm - x_{p+2}^\pm x_m^\pm,\end{equation*}
The following (conjectural) Drinfeld-Serre relations (\ref{serreun}), (\ref{serredeux}) are not used in the present paper :
\begin{equation}\label{serreun} \text{Sym}(q^{3/2} x_{k\mp 1}^\pm x_l^\pm x_m^\pm - (q^{1/2} + q^{-1/2})x_k^\pm x_{l\mp 1}^\pm x_m^\pm + q^{-3/2} x_k^\pm x_l^\pm x_{m\mp 1}^\pm) = 0,
\end{equation}
\begin{equation}\label{serredeux} \text{Sym}(q^{-3/2} x_{k\pm 1}^\pm x_l^\pm x_m^\pm - (q^{1/2} + q^{-1/2})x_k^\pm x_{l\pm 1}^\pm x_m^\pm + q^{3/2} x_k^\pm x_l^\pm x_{m\pm 1}^\pm) = 0.
\end{equation}

\subsection{Finite dimensional representations of twisted quantum loop algebras}\label{fdeptq} Denote by 
$\text{Rep}(\U_q(\Lo\Glie^{\sigma}))$ the Grothendieck ring of finite 
dimensional representations of $\U_q(\Lo\Glie^{\sigma})$ (we work with type $1$ representations).

\begin{defi}\label{lhigh} A $\U_q(\Lo\Glie^\sigma)$-module $V$ is said to be of $l$-highest weight if there is $v\in V$ eigenvector of all $\phi_{i,m}^{\pm}$ such that $V=\U_q(\Lo\Glie^\sigma)^-.v$ and $\forall i\in I, m\in\ZZ, x_{i,m}^+.v=0$.\end{defi}

For such a vector, denote by $\gamma = (\gamma_{i,\pm m}^{\pm})_{i\in I, m\geq 0}$ the corresponding eigenvalues of the $\phi_{i,\pm m}^{\pm}$. Note that we have necessarily $\gamma_{\sigma(i),\pm m}^{\pm} = \omega^{\pm m}\gamma_{i,\pm m}^{\pm}$ and $\gamma_{i,0}^+\gamma_{i,0}^- = 1$. 

Moreover if these conditions are satisfied, a unique simple module $L(\gamma)$ of $l$-highest weight $\gamma$ is considered in \cite{Cha5}. 

The following result was first stated in \cite{Cha5} :

\begin{thm}\label{fint} $L(\gamma)$ is finite dimensional (and non trivial) if and only if there exist polynomials $(P_i)_{i\in I}$, $P_i(u)\in\CC[u]$, such that $P_i(0) = 1$ and $\gamma$ satisfies in $\CC[[u]]$ (resp. in $\CC[[u^{-1}]]$) :
$$\underset{m\geq 0}{\sum}\gamma_{i,\pm m}^{\pm}u^{\pm 
m}=q^{M\text{deg}(P_i)}\frac{P_i(u^Mq^{-M})}{P_i(u^Mq^M)}\text{ if $i = \sigma(i)$,}$$
$$\underset{m\geq 0}{\sum}\gamma_{i,\pm m}^{\pm}u^{\pm 
m}=q^{\text{deg}(P_i)}\frac{P_i(uq^{-1})}{P_i(uq)}\text{ if $i\neq \sigma(i)$.}$$
Moreover all (type $1$) simple finite dimensional representations of $\U_q(\Lo\Glie^{\sigma})$ are of this form.\end{thm}

Remarks : 

We have necessarily $\sigma(i)\neq i\Rightarrow P_{\sigma(i)}(u) = P_i(\omega u)$. Thus $\gamma$ is determined by a set of polynomials indexed by $I_{\sigma}$.

The case $C_{i,\sigma(i)} = 0$ is different than the case $C_{i,\sigma(i)} = -1$ in the statement of \cite{Cha5}; here we choose the normalization of \cite{jin} $d_{\overline{i}} = 1/2$ if $C_{i,\sigma(i)} = -1$ which allows to unify both cases.

A priori the fact that $L(\gamma)$ is not trivial is not clear without the triangular decomposition of $\U_q(\Lo\Glie^{\sigma})$ (although this triangular decomposition should follow from the PBW basis of \cite{da, a, bn}). We partly rewrite the proof of \cite{Cha5} with references or proof of the relations between Drinfeld generators that are used (see the introduction).

The following relations were first stated in \cite{Cha5} :

\begin{prop}\label{relpr} We have the following relations in $\U_q^\tau$ for $r\geq 0$:
\begin{equation}\label{csept}(x_0^+)^{(r)}x_1^+ = -q^{-3r/2}[r-1]_{q^{1/2}}x_1^+(x_0^+)^{(r)} + q^{(-3r+3)/2} x_0^+x_1^+(x_0^+)^{(r-1)}\end{equation}
\begin{equation}\label{chuit}
\begin{split}
[h_1,(x_0^+)^{(r)}] [3]_{q^{1/2}}^{-1} = &(\frac{q^{3(1-r)/2} + q^{(3-r)/2} - q^{(1-r)/2} - q^{-(r+1)/2}}{q^{1/2} - q^{-1/2}})x_1^+(x_0^+)^{(r-1)} 
\\&+ q^{-r+2} x_0^+x_1^+(x_0^+)^{(r-2)}
\end{split}
\end{equation}
\begin{equation}\label{cneuf}
\begin{split}
[(x_0)^+,x_1^-] = &q^{(1-r)/2}kh_1(x_0^+)^{(r-1)} 
\\&+ \frac{q^{-r+1/2} + q^{-r-1/2} -q^{-r+5/2} - q^{-2r+5/2}}{q^{1/2} - q^{-1/2}}kx_1^+(x_0^+)^{(r-2)} 
\\&- q^{(-3r+5)/2}kx_0^+x_1^+(x_0^+)^{(r-3)}
\end{split}
\end{equation}
\begin{equation}\label{cdix}
\begin{split}
[(x_0^+)^{(r)},\tilde{e}_0] = &q^{(3-r)/2}K_1x_1^- (x_0^+)^{(r-1)} + q^{2-r}(q+q^{-1})k^2h_1(x_0^+)^{(r-2)}  
\\&+ q^{(6-3r)/2}\frac{q^{-5/2} + q^{-3/2} - q^{-3/2} - q^{-r+3/2}}{q^{1/2} - q^{-1/2}} k^2x_1^+ (x_0^+)^{(r-3)} 
\\&- q^{-2r+4} k^2 x_0^+x_1^+(x_0^+)^{(r-4)}
\end{split}
\end{equation}
\begin{equation}\label{conze}
\tilde{e}_0x_1^- - q^2x_1^-\tilde{e}_0 = 0
\end{equation}
\begin{equation}\label{cdouze}
[h_1,\tilde{e}_0^{(r)}] = q^{-2r + 5/2}(q^{1/2} - q^{-1/2})[3]_{q^{1/2}}[4]_{q^{1/2}}\tilde{e}_0^{(r-1)}(x_1^-)^2
\end{equation}
\begin{equation}\label{ctreize}[x_0^+,\tilde{e}_0^{(r)}] = q^{-2r+2}[4]_{q^{1/2}}\tilde{e}_0^{(r-1)}x_1^-k\end{equation}
\end{prop}

\demo By definition we have $x_0^+ = X^+_1$ and 
$$x_{-1}^+ = T(X^+_0) = T_1T_0(X^+_0) = T_1(-X^-_0K_0) = - (X^-_0X^-_1 - q^2 X^-_1X^-_0)K_1K_0.$$
Moreover from the relation $[X^+_0,X^-_0] = (K_0 - K_0^{-1})/(q^{1/2} - q^{-1/2}) = 0$ and $[X^+_0,X^-_1] = 0$ we get that $x_{-1}^+x_0^+ - q^{-1} x_0^+x_{-1}^+$ is equal to
\begin{equation*}
\begin{split}
 &- (X^-_0X^-_1 - q^2 X^-_1X^-_0)K_1K_0 X^+_0 + q^{-1} X_0^+(X^-_0X^-_1 - q^2 X^-_1X^-_0)K_1K_0 
\\= &q^{-1}[4]_{q^{1/2}} X^-_1 K_0^2.
\end{split}
\end{equation*}

From the quantum Serre relation $(X^-_1)^2X^-_0 - (q^2 + q^{-2}) X^-_1X^-_0X^-_1 + X^-_0(X^-_1)^2 = 0$ we get 
$$X^-_1(X^-_1X^-_0 - q^{-2}X^-_0X^-_1) + (X^-_0X^-_1 - q^2X^-_1X^-_0)X^-_1 = 0$$
$$-[x_{-1}^+,x_0^+]_{q^{-1}} K_0^2q^{-2}x_{-1}^+K_0 + x_{-1}^+K_0[x_{-1}^+,x_0^+]_{q^{-1}} K_0^2 = 0$$ 
$$-q^{-2}[x_{-1}^+,x_0^+]_{q^{-1}}x_{-1}^+ + x_{-1}^+ [x_{-1}^+,x_0^+]_{q^{-1}} = 0$$
$$(x_{-1}^+)^2 x_0^+ - (q^{-1} + q^{-2}) x_{-1}^+ x_0^+ x_{-1}^+ + q^{-3} x_0^+ (x_{-1}^+)^2 = 0$$
Now by applying the automorphism $T^{-1}$ we get :
$$(x_0^+)^2 x_1^+ - (q^{-1} + q^{-2}) x_0^+ x_1^+ x_0^+ + q^{-3} x_1^+ (x_0^+)^2 = 0$$
This is the relation (\ref{csept}) with $r = 2$. Then the relation (\ref{csept}) for general $r\geq 2$ follows from the case $r=2$ by induction on $r$ as stated in \cite{Cha5}.

We have $[h_1,x_0^+] =[2]_{q^{1/2}} x_1^+$. This is the relation (\ref{chuit}) with $r=1$. Then by using the relation (\ref{csept}) this relation (\ref{chuit}) is proved by induction for general $r\geq 1$ as stated in \cite{Cha5}.

We have $[x_0^+,x_1^-] = k h_1$. This is the relation (\ref{cneuf}) with $r=1$. Then by using the relation (\ref{csept}) and $[x_0^+,x_1^-] = k h_1$, $[h_1,x_0^+] = [2]_{q^{1/2}} x_1^+$ this relation (\ref{cneuf}) is proved by induction for general $r\geq 1$ as stated in \cite{Cha5}.

Consider $\tilde{e}_0 = x_0^-x_1^- - q x_1^- x_0^-$. By definition we have 
$$x_1^- = -K_1^{-1}K_0^{-1}(X^+_1X^+_0 - q^{-2}X^+_0X^+_1)\text{ and }x_0^- = X^-_0.$$ 
So by using the relation $[X^+_0,X^-_0] = (K_0 - K_0^{-1})/(q^{1/2} - q^{-1/2})$ and $[X^-_0,X^+_1] = 0$ we get
\begin{equation*}
\begin{split}
\tilde{e}_0 &= -X^-_0 K_1^{-1}K_0^{-1}(X^+_1X^+_0 - q^{-2}X^+_0X^+_1) + q K_1^{-1}K_0^{-1}(X^+_1X^+_0 - q^{-2}X^+_0X^+_1) X^-_0 
\\&= K_0^2X^+_1 q[4]_{q^{1/2}}.
\end{split}
\end{equation*}
So we have 
$$x_1^- = -K_1^{-1}K_0^{-1}\frac{q^{-1}}{[4]_{q^{1/2}}}(K_0^{-2}\tilde{e}_0X^+_0 - q^{-2}X^+_0K_0^{-2}\tilde{e}_0)$$
and so $[x_0^+,\tilde{e}_0] = q [4]_{q^{1/2}}K_0 x_1^-$. This is exactly the relation (\ref{cdix}) with $r = 1$. Then by using the relation (\ref{csept}) and $[x_0^+,x_1^-] = k h_1$ this relation (\ref{cdix}) is proved by induction for general $r\geq 1$ as stated in \cite{Cha5}.

From the quantum Serre relation $(X^+_1)^2X^+_0 - (q^2 + q^{-2})X^+_1X^+_0X^+_1 + X^+_0(X^+_1)^2 = 0$, we get 
\begin{equation*}
\begin{split}
\tilde{e}_0x_1^- - q^2x_1^-\tilde{e}_0 = &q[4]_{q^{1/2}} (-K_0^2X^+_1 K_1^{-1}K_0^{-1}(X^+_1X^+_0 - q^{-2}X^+_0X^+_1) 
\\&+ q^2 K_1^{-1}K_0^{-1}(X^+_1X^+_0 - q^{-2}X^+_0X^+_1)K_0^2X^+_1) = 0.
\end{split}
\end{equation*}
This is precisely the relation (\ref{conze}). 

We have 
$$[h_1,\tilde{e}_0] = [2]_{q^{1/2}}(1+q+q^{-1})(-(x_1^-)^2 - x_0x_2^- + qx_2^-x_0^- + q (x_1^-)^2).$$
So from the relation (\ref{int}) with $p=m=0$, we get
$$[h_1,\tilde{e}_0] = -q^{1/2} (q^{1/2} - q^{-1/2}) [3]_{q^{1/2}}[4]_{q^{1/2}} (x_1^-)^2.$$
This is precisely the relation (\ref{cdouze}) with $r=1$. Then by using the relation (\ref{conze}) this relation (\ref{cdouze}) is proved by induction for general $r\geq 1$ as stated in \cite{Cha5}.

The relation (\ref{cdix}) with $r=1$ is the relation (\ref{ctreize}) with $r=1$. Then by using the relation (\ref{conze}) this relation (\ref{ctreize}) is proved by induction for general $r\geq 1$ as stated in \cite{Cha5}.
\qed

We can now end the proof of Theorem \ref{fint} :

\demo First let us prove the "if" part. A simple finite dimensional representation corresponding to $(P_i)_{i\in I_\sigma}$ is said to be fundamental if $\sum_{i\in I_\sigma}\text{deg}(P_i) = 1$. For each $i\in I$ a fundamental representation such that $\text{deg}(P_i)=1$ has been constructed in \cite{kas} from level $0$ fundamental extremal modules (see \cite{naw} for the computation of the corresponding Drinfeld polynomials). Then we get all fundamental representations by twisting with the algebra automorphism $\tau_a : \U_q(\Lo\Glie^\sigma)\rightarrow\U_q(\Lo\Glie^\sigma)$ defined by $\tau_a(X_i^+) = X_i^+$, $\tau_a(X_i^-) = X_i^-$, $\tau_a(K_i) = K_i$ for $i\in I_\sigma$, $\tau_a(X_\epsilon^+) = a X_\epsilon^+$, $\tau_a(X_\epsilon^-) = a^{-1}X_\epsilon^-$, $\tau_a(K_\epsilon) = K_\epsilon$, that is to say $\tau_a(x_{i,m}^{\pm}) = a^m x_{i,m}^{\pm}$, $\tau_a(h_{i,m}) = a^m h_{i,m}$ and $\tau_a(k_i) = k_i$. The existence for general $n$-uplet of polynomials follows from Theorem \ref{mult}.

So it suffices to prove the "only if" part for $\U_q^\tau$. We see how to prove it without using relations (\ref{serreun}) and (\ref{serredeux}) (other relations were proved in \cite{a}). We follow the proof of \cite{Cha5} without using these relations. The crucial point for this is \cite[Proposition 4.2]{Cha5}. The result is a consequence of the relations (7 - 16) in the proof of \cite[Proposition 4.2]{Cha5}. So it suffices to prove that these relations hold. This is a consequence of Proposition \ref{relpr} (the relation (14) of \cite{Cha5} is the relation (\ref{cdix}) modified by using relations (\ref{chuit}) and (\ref{cneuf}, the relations (15) and (16) of \cite{Cha5} follow from the previous relations as stated in \cite{Cha5}).\qed

\subsection{Subalgebras of $\U_q(\Lo\Glie^{\sigma})$}\label{subalg}

The type of $\Glie^\sigma$ corresponding to the type of $\hat{\Glie}^\sigma$ can be read in the following table (here $n\geq 2$) :
\begin{equation*}
\begin{array}{l|l|l|l|l|l|l}
   \hat{\Glie}^\sigma  &A_2^{(2)} & A_{2n}^{(2)} & A_{2n-1}^{(2)} & D_{n+1}^{(2)} & E_6^{(2)} & D_4^{(3)}
\\ \Glie^\sigma        & A_1      &  B_n         & C_n            & B_n           & F_4       & G_2
\end{array}
\end{equation*}

Let $\overline{\U}_q(\Glie^\sigma)$ be the subalgebra of $\U_q(\Lo\Glie^{\sigma})$ generated by the $x_{i,0}^{\pm}, k_i^{\pm 1}$ where $i\in I_\sigma$. Let $\tilde{\U}_q(\Glie^\sigma)$ be the subalgebra of of $\U_q(\Lo\Glie^\sigma)$ generated by the $X_i^{\pm}, K_i^{\pm 1}$, $1\leq i\leq n$.

If $\hat{\Glie}^\sigma$ is not of type $A_{2n}^{(2)}$, we have $\overline{\U}_q(\Glie^\sigma) = \tilde{\U}_q(\Glie^\sigma)\simeq \U_q(\Glie^\sigma)$.

If $\hat{\Glie}^\sigma$ is of type $A_{2n}^{(2)}$ where $n\geq 2$, we have $\overline{\U}_q(\Glie^\sigma)\simeq \U_{q^{\frac{1}{2}}}(\Glie^\sigma) \simeq \U_{q^{\frac{1}{2}}}(B_n)$ and $\tilde{\U}_q(\Glie^\sigma)\simeq \U_q(C_n)$.

If $\hat{\Glie}^\sigma$ is of type $A_2^{(2)}$, we have $\overline{\U}_q(\Glie^\sigma)\simeq \U_{q^{\frac{1}{2}}}(\Glie^\sigma) = \U_{q^{\frac{1}{2}}}(sl_2)$ and $\tilde{\U}_q(\Glie^\sigma)\simeq \U_{q^2}(sl_2)$.

Note that we have a natural grading of $\U_q(\Lo\Glie^\sigma)$ by the weight lattice of $\overline{\U}_q(\Glie^\sigma)$ and by the weight lattice of $\tilde{\U}_q(\Glie^\sigma)$.

For $i\in I$, denote by $\hat{\U}_i$ the subalgebra of $\U_q(\Lo\Glie^{\sigma})$ generated by the $x_{i,m}^{\pm}$ ($m\in\ZZ$), $h_{i,r}$ ($r\in\ZZ-\{0\}$) and $k_i^{\pm 1}$. We have different cases \cite{da} :

$\bullet$ $C_{i,\sigma(i)} = 2$ : for $m\neq 0 [M]$, we have $ 0 = (1 - \omega^m)x_{i,m}^{\pm} = (1 - \omega^m)h_{i,m}$ and so $x_{i,m}^{\pm} = h_{i,m} = 0$. This implies $\phi_{i,m}^{\pm} = 0$. So $\U_i$ is generated by the $x_{i,Mm}^{\pm}, \phi_{i,Mm}$ where $m\in\ZZ$. Moreover we have an algebra isomorphism : 
$$\U_{q^M}(\Lo sl_2)\rightarrow \hat{\U}_i$$
$$x_k^+\mapsto x_{i,kM}^+/M\text{ , }x_k^-\mapsto x_{i,kM}^-\text{ , }\phi_{m}^{\pm}\mapsto\phi_{i,Mm}^{\pm}.$$ 

$\bullet$ $C_{i,\sigma(i)} = 0$ : we have an algebra isomorphism
$$\U_q(\Lo sl_2)\rightarrow \hat{\U}_i$$
$$x_k^+\mapsto x_{i,k}^+\text{ , }x_k^-\mapsto x_{i,k}^-\text{ , }\phi_{m}^{\pm}\mapsto\phi_{i,m}^{\pm}.$$

$\bullet$ $C_{i,\sigma(i)} = -1$ : we have an algebra morphism 
$$\U_q^\tau\rightarrow \hat{\U}_i$$
$$x_k^+\mapsto x_{i,k}^+\text{ , }x_k^-\mapsto x_{i,k}^-\text{ , }\phi_{m}^{\pm}\mapsto\phi_{i,m}^{\pm}.$$ 

In particular the "elementary" subalgebras $\hat{\U}_i$ may be of type $A_1^{(1)}$ or of type $A_2^{(2)}$. This is an important difference with the theory of untwisted quantum affine algebras where all "elementary" subalgebras are of type $A_1^{(1)}$ and this is one of the reason why several new technical developments are considered in the present paper.

Remark : the last situation $C_{i,\sigma(i)} = -1$ appears only for types $A_{2n}^{(2)}$. 

\subsection{Invariant subalgebra of the quantum loop algebra} In \cite{jin}, the following subalgebra $\U_q'(\Lo\Glie^\sigma)$ of $\U_q(\Lo\Glie)$ is introduced. $\U_q'(\Lo\Glie^\sigma)$ is generated by the :
$$x_{i,r}'^{\pm} = \frac{1}{[d_{\overline{i}}]_qM^{\frac{1}{2}}}\sum_{s = 0\cdots M-1}\tilde{x}_{\sigma^s(i),m}^{\pm}\omega^{-ms},$$
$$h_{i,m}' = \frac{1}{[d_{\overline{i}}]_q}\sum_{s = 0\cdots M-1}\tilde{h}_{\sigma^s(i),m}\omega^{-ms}\text{ , }k_i' = \prod_{s=0\cdots M-1}\tilde{k}_{\sigma^s(i)},$$
where the $\tilde{x}_{i,m}^{\pm}$, $\tilde{h}_{i,m}$, $\tilde{k}_i^{\pm 1}$ are the generators of $\U_q(\Lo\Glie)$.
 
In particular $\U'_q(\Lo\Glie^\sigma)$ is a subalgebra invariant by the automorphism $\sigma$ of $\U_q(\Lo\Glie)$ defined by 
\begin{equation*}
\sigma(\tilde{k}_i) = \tilde{k}_{\sigma(i)}\text{ , }\sigma(\tilde{x}_{i,r}^{\pm}) = \omega^{-r} \tilde{x}_{\sigma(i),r}^{\pm}\text{ , }\sigma(\tilde{h}_{i,r}) = \omega^{-r} \tilde{h}_{\sigma(i),r}. 
\end{equation*}
But there are invariant elements which are not in this subalgebra as for example $\prod_{s=0\cdots M - 1}\tilde{h}_{\sigma^s(i),r}$.

The identification $\theta(x_{i,r}^\pm) = x_{i,r}'^\pm$, $\theta(h_{i,m}) = h_{i,m}'$, $\theta(k_i) = k'_i$ considered in \cite{jin} does not give an isomorphism between $\U_q(\Lo\Glie^\sigma)$ and $\U'_q(\Lo\Glie^\sigma)$. For example for $\U_q^\tau$ we should have $\theta(k) = \tilde{k}_1\tilde{k}_2$, $\theta([x_0^+,x_0^-])=\frac{\tilde{k}_1\tilde{k}_2 - \tilde{k}_1^{-1}\tilde{k}_2^{-1}}{q^{\frac{1}{2}} - q^{-\frac{1}{2}}}$, but also 
\begin{equation*}
\begin{split}
\theta([x_0^+ , x_0^-])  &= \frac{1}{2}([\tilde{x}_{1,0}^+,\tilde{x}_{1,0}^-] + [\tilde{x}_{2,0}^+,\tilde{x}_{2,0}^-])
\\               &= \frac{1}{2(q^{\frac{1}{2}}-q^{-\frac{1}{2}})}(\tilde{k}_1 - \tilde{k}_1^{-1} + \tilde{k}_2 - \tilde{k}_2^{-1}),
\end{split}
\end{equation*}
but we do not have 
$$2(\tilde{k}_1\tilde{k}_2 - \tilde{k}_1^{-1}\tilde{k}_2^{-1}) = \tilde{k}_1+\tilde{k}_2 - \tilde{k}_1^{-1} - \tilde{k}_2^{-1}$$ 
as $\tilde{k}_1$ and $\tilde{k}_2$ are algebraically independent in $\U_q(\Lo sl_3)$.

\begin{rem}\label{sitpart} Note that in the particular situations where we specialize 
$$\tilde{k}_1 = 1\text{ or }\tilde{k}_2 = 1\text{ or }\tilde{k}_1 = \tilde{k}_2^{-1},$$ 
the above relation is satisfied (but is does not imply that the morphism is well-defined under this condition). This situation will appear in the next sections.\end{rem}

A priori it seems difficult to get direct information from the representation theory in the untwisted case to the twisted case as in general : 

\begin{prop} There is no algebra morphism $\theta:\U_q^\tau \rightarrow \U_q(\Lo sl_3)$ satisfying one of the following properties :

1) $\theta(x_0^{\pm}) = \alpha^\pm (\tilde{x}_{1,0}^\pm + \beta^\pm \tilde{x}_{2,0}^\pm )$, where $\alpha^+,\alpha^-,\beta^+,\beta^- \in\CC^*$,

2) $\theta(k) = \tilde{k}_1\tilde{k}_2$ and $\theta(x_0^{\pm})\in\U_q^{\pm}(\Lo sl_3)$.

\end{prop}

\demo 

1) We have :
\begin{equation*}
\begin{split}
\theta([x_0^+,x_0^-]) &= \alpha^+\alpha^-  ([\tilde{x}_{1,0}^+,\tilde{x}_{1,0}^-] + \beta^+\beta^- [\tilde{x}_{2,0}^+,\tilde{x}_{2,0}^-])
\\&=\frac{\alpha^+\alpha^-}{(q - q^{-1})}(\tilde{k}_1 - \tilde{k}_1^{-1} + \beta^+\beta^-(\tilde{k}_2 - \tilde{k}_2^{-1})).
\end{split}
\end{equation*}
So 
$$\theta(k - k^{-1}) = \alpha^+\alpha^- \frac{\tilde{k}_1 - \tilde{k}_1^{-1} + \beta^+\beta^-(\tilde{k}_2 - \tilde{k}_2^{-1})}{[\frac{1}{2}]_q}.$$
We have the relation $(k - k^{-1}) x_0^+ = x_0^+ (q k - q^{-1} k)$. But $\theta((k - k^{-1})x_0^+)$ is equal to :
\begin{equation*}
\begin{split}
 &\alpha^+\alpha^-\frac{\tilde{k}_1 - \tilde{k}_1^{-1} + \beta^+\beta^-(\tilde{k}_2 - \tilde{k}_2^{-1})}{[\frac{1}{2}]_q}\alpha^+(\tilde{x}_1^+ +\beta^+ \tilde{x}_2^+)
\\= &\frac{(\alpha^+)^2\alpha^-}{[\frac{1}{2}]_q} (\tilde{x}_1^+ (q^2\tilde{k}_1 - q^{-2}\tilde{k}_1^{-1} + \beta^+\beta^-(q^{-1}\tilde{k}_2 - q\tilde{k}_2^{-1})) 
\\&+ \beta^+\tilde{x}_{2,r}^+ (q^{-1}\tilde{k}_1 - q\tilde{k}_1^{-1} + \beta^+\beta^-(q^2\tilde{k}_2 - q^{-2}\tilde{k}_2^{-1}))),
\end{split}
\end{equation*}
and 
$$\theta(x_0^+ (qk - q^{-1}k^{-1})) = \alpha^+ (\tilde{x}_1^+ + \beta^+ \tilde{x}_2^+)(q\theta(k) - q^{-1}\theta(k)^{-1}).$$
So :
$$q^2\tilde{k}_1 - q^{-2}\tilde{k}_1^{-1} + \beta^+\beta^-(q^{-1}\tilde{k}_2 - q\tilde{k}_2^{-1})$$
$$= q^{-1}\tilde{k}_1 - q\tilde{k}_1^{-1} + \beta^+\beta^-(q^2\tilde{k}_2 - q^{-2}\tilde{k}_2^{-1}),$$
contradiction.

2) As $kx_0^+ = q^{\pm 1} x_0^+k$, $\theta(x_0^+)$ is a sum of elements with weight of form $r_1\alpha_1 + r_2\alpha_2$ where $r_1 + r_2 = 1$. By hypothesis $r_1,r_2\geq 0$. We get $\theta(x_0^+) = a \tilde{x}_1^+ + b \tilde{x}_2^+$ with $a,b\in\CC$. In the same way $\theta(x_0^-) = c \tilde{x}_1^- + d \tilde{x}_2^-$. Contradiction as $$[\theta(x_0^+),\theta(x_0^-)] = \frac{\tilde{k}_1\tilde{k}_2 - \tilde{k}_2^{-1}\tilde{k}_1^{-1}}{q^{\frac{1}{2}} - q^{-\frac{1}{2}}}.$$\qed

So a priori it is not clear if $\U_q'(\Lo\Glie^\sigma)$ has a Hopf algebra structure or if its representation theory is related to the one of $\U_q(\Lo\Glie^\sigma)$. As the main subject of this paper is to study twisted quantum affine algebras, we will not discuss $\U_q'(\Lo\Glie^\sigma)$ outside this section.

First from the triangular decomposition of $\U_q(\Lo\Glie)$ we get a triangular decomposition of $\U_q'(\Lo\Glie^\sigma)$ :
$$\U_q'^+(\Lo\Glie^\sigma)\otimes \U_q'(\Lo\Hlie^\sigma)\otimes \U_q'^-(\Lo\Glie^\sigma)\simeq \U_q'(\Lo\Glie^\sigma),$$
where $\U_q'^\pm(\Lo\Glie^\sigma)$ (resp. $\U_q'(\Lo\Hlie^\sigma)$) is the subalgebra generated by the $x_{i,m}'^\pm$ (resp. the $h_{i,r}'$ and the $k_i'^{\pm 1}$).

So we can define the simple $\U_q'(\Lo\Glie^\sigma)$-modules $\tilde{L}(\gamma)$ as for $\U_q(\Lo\Glie^\sigma)$. 

\begin{prop}Suppose that the conditions of Theorem \ref{fint} are satisfied by $\gamma$. Then $\tilde{L}(\gamma)$ is finite dimensional.\end{prop}

\demo As $\U_q'(\Lo\Glie^\sigma)$ is a subalgebra of $\U_q(\Lo\Glie)$ the result is a direct consequence of Theorem \ref{finu}.\qed

\section{Twisted $q$-characters}\label{tw}

In this section we develop the theory of twisted Frenkel-Reshetikhin $q$-characters which is a crucial tool in the proofs of the conjectures in this paper. It is a generalization of the theory for untwisted quantum affine algebras developed in \cite{Fre, Fre2}. For many points the generalization is not straightforward, as for example the characterization of the image. The existence of such a theory for twisted cases was expected in \cite{Fre, Fre2}.

\subsection{Definition}

First let us explain the construction of twisted $q$-characters.

\subsubsection{$l$-weight spaces}

From the relation (\ref{eqddeux}), the subalgebra $\U_q(\Lo\Hlie^{\sigma})\subset\U_q(\Lo\Glie^{\sigma})$ is commutative. So for $V$ a finite dimensional $\U_q(\Lo\Glie^{\gamma})$-module we have :
$$V=\underset{\gamma=(\gamma_{i,\pm m}^{\pm})_{i\in I, m\geq 0}}{\bigoplus}V_{\gamma},$$
where :
$$V_{\gamma}=\{v\in V | \exists p\geq 0, \forall i\in I, m\geq 0, (\phi_{i,\pm m}^{\pm}-\gamma_{i,\pm m}^{\pm})^p.v=0\}.$$
The $V_{\gamma}$ are called the $l$-weight spaces of $V$. We prove in the following that the $\gamma$ satisfying $V_{\gamma}\neq \{0\}$ have a particular form. For an $l$-weight $\gamma$ we consider the generating series $$\gamma_i^+(z) = \sum_{m\geq 0}\gamma_{i,m}^+z^m\in\CC[[z]]\text{ , }\gamma_i^-(z) = \sum_{m\leq 0}\gamma_{i,-m}^-z^{-m}\in\CC[[z^{-1}]].$$

\subsubsection{Example : case of $\U_q^\tau$}\label{exfund} We study $\U_q^\tau$.

For $a\in\CC^*$, the simple finite dimensional representation corresponding to $P(u) = 1 - aq^{-1} u$ is of dimension $3$ and was explicitly constructed in \cite{Cha5} with Drinfeld generators. Let us give this construction in terms of the Drinfeld-Jimbo generators.

Let $V = \CC v_0\oplus \CC v_1\oplus \CC v_2$. The action of $\U_q^\tau$ is defined as follows 
\begin{equation*}
\begin{array}{l|l|l|l}
       &v_0                                       & v_1                         & v_2
\\ X^+_0 & 0                                        & [2]_{q^{\frac{1}{2}}} v_0  & [2]_{q^{\frac{1}{2}}} v_1
\\ X^-_0 & v_1                                      & v_2                        & 0
\\ K_0 & qv_0                                     & v_1                        & q^{-1}v_2
\\ X^+_1 & [4]_{q^{\frac{1}{2}}}^{-1} q(a(1+q^2))v_2& 0                          & 0
\\ X^-_1 & 0                                        & 0                          & [4]_{q^{\frac{1}{2}}}^{-1} q^{-1}[2]_{q^{\frac{1}{2}}}^2a^{-1}(1+q^{-2})v_0
\\ K_1 & q^{-2}v_0                                & v_1                        & q^2 v_2
\end{array}
\end{equation*}
We check directly that the relations and in particular the quantum Serre relations are satisfied.

\demo

The quantum Serre relations with $i = 0$ and $j = 1$ are trivial as $(X^+_0)^3 = (X^-_0)^3 = 0$ on the representation. 

The quantum Serre relations with $i = 1$ and $j = 0$ are also clear as $(X^+_1)^2 = (X^-_1)^2 = 0$ on the representation, and $\text{Im}(X^+_0X^+_1) = \CC v_1\subset \text{Ker}(X^+_1)$ so $X^+_1X^+_0X^+_1 = 0$ on the representation, and $\text{Im}(X^-_0X^-_1) = \CC v_1\subset \text{Ker}(X^-_1)$ so $X^-_1X^-_0X^-_1 = 0$ on the representation.

The relations between the $E,F$ and the $K$ are clear by construction.

We have $[X^+_0,X^-_1] = 0$ as $\text{Im}(X^-_1) = \CC v_0\subset\text{Ker}(X^+_0)$, and $\text{Im}(X^+_0) = \CC v_0\oplus \CC v_1\subset\text{Ker}(X^-_1)$ and so $X^-_1X^+_0 = X^+_0X^-_1 = 0$ on the representation. In the same way we have $[X^+_1,X^-_0]=0$.

The action of $[X^+_0,X^-_0]$ is diagonal with eigenvectors $(v_0,v_1,v_2)$ with respective eigenvalues $([2]_{q^{\frac{1}{2}}},0,-[2]_{q^{\frac{1}{2}}})$ and so is equal to the action of $\frac{K_0 - K_0^{-1}}{q^{\frac{1}{2}} - q^{-\frac{1}{2}}}$.

The action of $[X^+_1,X^-_1]$ is diagonal with eigenvectors $(v_0,v_1,v_2)$ with respective eigenvalues $(-\frac{[2]_{q^{\frac{1}{2}}}(1+q^2)(1+q^{-2})}{[4]_{q^{\frac{1}{2}}}},0,\frac{[2]_{q^{\frac{1}{2}}}(1+q^2)(1+q^{-2})}{[4]_{q^{\frac{1}{2}}}}) = (-1,0,1)$ and so is equal to the action of $\frac{K_1 - K_1^{-1}}{q^2 - q^{-2}}$.
\qed

We can compute the action in terms of the Drinfeld generators and we recover the formulas of \cite{Cha5} ($m\in\ZZ$, $r\in\ZZ-\{0\}$) :
\begin{equation*}
\begin{array}{l|l|l|l}
&v_0&v_1&v_2
\\ x_m^+                         & 0              & a^m[2]_{q^{\frac{1}{2}}}v_0  & [2]_{q^{\frac{1}{2}}}(-aq)^{m}v_1
\\ x_m^-                         & a^m v_1        &(-aq)^m v_2           & 0
\\ k                             &qv_0            &v_1                   &q^{-1}v_2
\\ \frac{r h_r }{[2r]_{q^{\frac{1}{2}}}} &(aq^{-1})^r v_0 &((-a)^r - (aq)^r) v_1 &- (-aq^2)^rv_2
\end{array}
\end{equation*}
In particular for $m > 0$ :
\begin{equation*}
\begin{split}
\phi_{\pm m}^{\pm}.v_0 &= \pm a^{\pm } (q-q^{-1})v_0
\\\phi_{\pm m}^{\pm}.v_1 &= \pm (q-q^{-1})((-aq)^{\pm m} - a^{\pm m}) v_1,
\\\phi_{\pm m}^\pm.v_2 &= \mp(q - q^{-1})(-aq)^{\pm m}v_2.
\end{split}
\end{equation*}
Remark : we are in the situation of remark \ref{sitpart}. 

It is checked in \cite{Cha5} that the relations between Drinfeld generators are satisfied.

We have the $l$-weight spaces $\CC.v_0 = V_{\gamma_0}$, $\CC.v_1 = V_{\gamma_1}$, $\CC.v_2 = V_{\gamma_2}$ where $\gamma_0,\gamma_1,\gamma_2$ are given by the above formulas. Let us compute the corresponding generating series. First in $\CC[[z]]$ (and in $\CC[[z^{-1}]]$ for the second identity):
$$\gamma_0^+(z) = q + \sum_{m\geq 1}a^m (q - q^{-1})z^m = q + (q - q^{-1})\frac{az}{1-az} = q\frac{1 - q^{-2}az}{1-az},$$
$$\gamma_0^-(z) = q^{-1} - \sum_{m\geq 1}a^{-m} (q - q^{-1})z^{-m} = q^{-1} + (-q + q^{-1})\frac{(az)^{-1}}{1 - (az)^{-1}}
= q\frac{1-aq^{-2}z}{1 - az}.
$$
In fact $\gamma_0^+(z),\gamma_0^-(z)$ where already computed in \cite{Cha5} and correspond to the $l$-highest weight. In the present paper we are also interested in the lower $l$-weights.
\begin{equation*}
\begin{split}
\gamma_1^+(z) &= 1 + \sum_{m\geq 1}(q - q^{-1})((-aq)^m - a^m)z^m
\\&= 1 + (q - q^{-1})(\frac{-aqz}{1 + aqz} -\frac{az}{1-az}) 
= \frac{(1-zaq^2)(1+zaq^{-1})}{(1+aqz)(1-az)},
\\\gamma_1^-(z) &= 1 - \sum_{m\geq 1}(q - q^{-1}) ((-aq)^{-m} - a^{-m})z^{-m}
\\&= 1 + (-q + q^{-1})(\frac{-(aqz)^{-1}}{1+(aqz)^{-1}}-\frac{(az)^{-1}}{1 - (az)^{-1}} ) 
\\&= \frac{(1-(zaq^2)^{-1})(1+(zaq^{-1})^{-1})}{(1+(aqz)^{-1})(1-(az)^{-1})}
= \frac{(1-zaq^2)(1+zaq^{-1})}{(1+aqz)(1-az)},
\\\gamma_2^+(z) &= q^{-1} - \sum_{m\geq 1}(-aq)^m (q - q^{-1})z^m = q^{-1} + (-q + q^{-1})\frac{-aqz}{1 + aqz} 
= q^{-1}\frac{1 + q^3az}{1 + aqz},
\\\gamma_2^-(z) &= q + \sum_{m\geq 1}(-aq)^{-m} (q - q^{-1})z^{-m} = q + (q - q^{-1})\frac{(-aqz)^{-1}}{1 - (-aqz)^{-1}} 
= q^{-1}\frac{1 + aq^3 z}{1 + aqz}.
\end{split}
\end{equation*}

We have :
\begin{equation*}
\begin{split}
\gamma_0^{\pm}(z) &= q^{\text{deg}(P_0)} \frac{P_0(zq^{-1})}{P_0(zq)},
\\\gamma_1^{\pm}(z) &= q^{(\text{deg}(P_1) - \text{deg}(Q_1))}\frac{P_1(zq^{-1})Q_1(zq)}{P_1(zq)Q_1(zq^{-1})},
\\\gamma_2^{\pm}(z) &= q^{-\text{deg}(Q_2)}\frac{Q_2(zq)}{Q_2(zq^{-1})},
\end{split}
\end{equation*}
where $P_0(z) = (1 - azq^{-1})$, $P_1(z) = (1 + az)$, $Q_1(z) = (1- azq)$, $Q_2(z) = (1 + zaq^2)$.

In fact we can deduce from this observation the general statement :

\begin{lem}\label{casfac} Let $V$ be a finite dimensional representation of $\U_q^\tau$ and $\gamma$ be an $l$-weight of $V$. Then there are polynomials $P,Q\in\CC[z]$ such that $P(0) = Q(0) = 1$ and in $\CC[[z^{\pm 1}]]$ : 
$$\gamma^{\pm}(z) = q^{(\text{deg}(P) - \text{deg}(Q))}\frac{P(zq^{-1})Q(zq)}{P(zq)Q(zq^{-1})}.$$
\end{lem}

\demo From the above computation and Theorem \ref{mult}, Theorem \ref{fint}, the result is true for simple finite dimensional representations, and so for any finite dimensional representations.\qed

\subsubsection{Definition of twisted $q$-characters}\label{defimono} 

\begin{thm}\label{deftw} Let $V$ be a finite dimensional representation of $\U_q(\Lo\Glie^\sigma)$ and $\gamma$ an $l$-weight of $V$. Then there exists polynomials $(P_i)_{i\in I}, (Q_i)_{i\in I}$, $P_i(u),Q_i(u)\in\CC[u]$, such that $P_i(0) = Q_i(0) = 1$ and $\gamma$ satisfies in $\CC[[u]]$ (resp. in $\CC[[u^{-1}]]$) :
$$\gamma_i^{\pm}(u) = q^{M(\text{deg}(P_i) - \text{deg}(Q_i))}\frac{P_i(u^Mq^{-M})Q_i(u^Mq^M)}{P_i(u^Mq^M)Q_i(u^Mq^{-M})}\text{ if $i = \sigma(i)$,}$$
$$\gamma_i^{\pm}(u) = q^{\text{deg}(P_i)-\text{deg}(Q_i)}\frac{P_i(uq^{-1})Q_i(uq)}{P_i(uq)Q_i(uq^{-1})}\text{ if $i\neq \sigma(i)$.}$$
\end{thm}

\demo It suffices to look at the subalgebras $\hat{\U}_i$ considered in section \ref{subalg}, and then the result follows from Lemma \ref{casfac} and Section \ref{unt}.
\qed

Remark : we have necessarily $\sigma(i)\neq i\Rightarrow P_{\sigma(i)}(u) = P_i(\omega u)$, $Q_{\sigma(i)}(u) = Q_i(\omega u)$.

First we can define the twisted version of the Frenkel-Reshetikhin $q$-characters morphism $\chi_q$ \cite{Fre} (see also \cite{kn}) considered in section \ref{unt}.

Let $\text{Rep}(\U_q(\Lo\Glie^\sigma))$ be the Grothendieck group of (type $1$) finite dimensional representations of $\U_q(\Lo\Glie^\sigma)$. 

\begin{defi} The twisted $q$-character morphism is the group morphism
$$\chi_q^\sigma:\text{Rep}(\U_q(\Lo\Glie^\sigma))\rightarrow \ZZ[Y_{i,a}^{\pm}]_{i\in I, a\in\CC^*}\text{ , }\chi_q^\sigma(V)=\underset{\gamma}{\sum}\text{dim}(V_{\gamma})m_{\gamma}$$
where 
$$m_{\gamma}=\underset{i\in I, a\in\CC^*}{\prod}Y_{i,a}^{q_{i,a}-r_{i,a}}\text{ , }P_i(u)=\underset{a\in\CC^*}{\prod}(1-ua)^{q_{i,a}}\text{ , }Q_i(u)=\underset{a\in\CC^*}{\prod}(1-ua)^{r_{i,a}},$$
where $P_i$, $Q_i$ are the polynomials considered in Theorem \ref{deftw}.
\end{defi}

The $m_{\gamma}$ are called monomials or $l$-weight (they are analogs of weight) and we denote $V_{\gamma}=V_{m_{\gamma}}$. Note that $\chi_q^\sigma$ also makes sense for finite dimensional representations of $\U_q(\Lo\Hlie^\sigma)$ which are a sub $\U_q(\Lo\Hlie^\sigma)$ of a finite dimensional representation of $\U_q(\Lo\Glie^\sigma)$.

\subsection{First properties}

For each $\overline{i}\in I_{\sigma}$ choose a representative $i\in I_{\sigma}$. We can make this choice such that $$(C_{i,j},C_{i,\sigma(j)},\cdots,C_{i,\sigma^{M-1}(j)})\neq (0,\cdots,0)\Rightarrow C_{i,j} = -1.$$ 
This choice is fixed in the following. Moreover we identify $i$ and $\overline{i}$.

Let $i\in I$ such that $i\neq \sigma(i)$. As above for any $a\in\CC^*$, $r\in\ZZ$, $q_{\sigma^r(i),a} = q_{i,a\omega^r}$ and $r_{\sigma^r(i),a} = r_{i,a\omega^r}$.  For $a\in\CC^*$, we put $Z_{\overline{i},a} = \prod_{k=1\cdots M} Y_{\sigma^k(i),a\omega^k}$. 

Let $i\in I$ such that $i = \sigma(i)$. For $a\in\CC^*$, we put $Z_{\overline{i},a} = Y_{i,a}$.

As a conclusion, we have $\text{Im}(\chi_q^{\sigma})\subset \ZZ[Z_{i,a}^{\pm 1}]_{a\in\CC^*,i\in I_{\sigma}}$. In the following we will consider the twisted $q$-character morphism :
$$\chi_q^\sigma:\text{Rep}(\U_q(\Lo\Glie^\sigma))\rightarrow \ZZ[Z_{i,a}^{\pm}]_{i\in I_\sigma, a\in\CC^*}.$$

\begin{thm} $\chi_q^\sigma$ is an injective ring morphism.
\end{thm}

\demo 
The additivity is clear. As the $\{\chi_q^\sigma(L)|\text{ $L$ is simple}\}$ are linearly independent, the map is injective. The multiplicity follows from Theorem \ref{mult}.\qed

As in the untwisted case, we have the following consequence :

\begin{cor} 
The Grothendieck ring $\text{Rep}(\U_q(\Lo\Glie^\sigma))$ is commutative.
\end{cor}

\noindent For $J\subset I_\sigma$, a monomial $m=\underset{i\in I_\sigma, a\in\CC^*}{\prod}Z_{i,a}^{z_{i,a}(m)}$ is said to be 
$J$-dominant if for all $j\in J, a\in\CC^*$ we have $z_{j,a}(m)\geq 0$. An $I_\sigma$-dominant monomials is said to be dominant. 

\noindent In the following for $M$ a finite dimensional $\U_q(\Lo\Glie^\sigma)$-module, we denote by $\mathcal{M}(M)$ the set of monomials occurring in $\chi_q^\sigma(M)$. 

\subsection{Restriction maps}\label{fprop}

The compatibility with the restriction functors requires additional work.

Let $\chi^\sigma:\text{Rep}(\tilde{U}_q(\Glie^{\sigma}))\rightarrow \ZZ[z_i^{\pm 1}]_{1\leq i\leq n}$ be the usual character map (here we denote $z_i = e^{\Lambda_i}$). Let 
$$\beta : \ZZ[Z_{i,a}^{\pm 1}]_{a\in\CC^*,i\in I_{\sigma}}\rightarrow \ZZ[z_i^{\pm 1}]_{1\leq i\leq n},$$ 
be the ring morphism such that for any $i\in I_\sigma, a\in\CC^*$ :
\begin{equation*}
\begin{split}
\beta(Z_{i,a}) = \begin{cases}
z_i\text{ if $\hat{\Glie}^\sigma \neq A_{2n}^{(2)}$,}
\\z_{n-i}\text{ if $\hat{\Glie}^\sigma = A_{2n}^{(2)}$.}
\end{cases}
\end{split}
\end{equation*}
Eventually let $res^\sigma$ be the restriction map from $\text{Rep}(\U_q(\hat{\Glie}^{\sigma}))$ to $\text{Rep}(\tilde{\U}_q(\Glie^{\sigma}))$. 

\begin{prop}\label{dag} The following diagram is commutative :
$$\begin{CD}
\text{Rep}(\U_q(\hat{\Glie}^{\sigma}))  @>{\chi_q^\sigma}>>  \ZZ[Z_{i,a}^{\pm 1}]_{i\in I_\sigma,a\in\CC^*}\\
@VV{\text{res}^\sigma}V     @VV{\beta}V\\
\text{Rep}(\U_q(\Glie^{\sigma}))  @>{\chi^\sigma}>>  \ZZ[z_i^{\pm 1}]_{1\leq i\leq n}
\end{CD}.$$
\end{prop}

The analog result for untwisted quantum affine algebras was proved in \cite{Fre}. In general the proof is modified, as for example for type $A_{2n}^{(2)}$, $\overline{\U}_q(\Glie^\sigma)$ and $\tilde{\U}_q(\Glie^\sigma)$ are not isomorphic.

\demo Let $m$ be a monomial corresponding to $((P_i,Q_i))_{i\in I_\sigma}$ and let $i\in I_\sigma$.

If $i = \sigma(i)$ : the multiplicity $z_i(m)$ of $z_i$ in $\beta(m)$ is the sum of the multiplicities of the $Z_{i,a}$ in $m$ for $a\in\CC^*$, that is to say 
$$z_i(m) = \text{deg}(P_i) - \text{deg}(Q_i).$$ 
The corresponding eigenvalue of $k_i$ is $q^{M(\text{deg}(P_i)-\text{deg}(Q_i))} = q^{d_i z_i(m)}$.

If $i\neq \sigma(i)$ : the multiplicity $z_i(m)$ of $z_i$ in $\beta(m)$ is the sum of the multiplicities of the $Z_{j,a}$ in $m$ for $a\in\CC^*$, $j\in\overline{i}$, that is to say 
$$z_i(m) = M(\text{deg}(P_i) - \text{deg}(Q_i))/M = \text{deg}(P_i) - \text{deg}(Q_i).$$ 
The corresponding eigenvalue of $k_i$ is $q^{(\text{deg}(P_i)-\text{deg}(Q_i))} = q^{d_i z_i(m)}$.

If $\hat{\Glie}^\sigma$ is not of type $A_{2n}^{(2)}$, for $i\in I_\sigma$, $q^{d_i z_i(m)}$ is the eigenvalue of $k_i = K_i\in\U_q(\hat{\Glie}^\sigma)$ corresponding to $z_i^{z_i(m)}$. So the result is clear.

For $\hat{\Glie}^\sigma$ of type $A_{2n}^{(2)}$ and $0\leq i\leq n-1$ the eigenvalue of $k_i$ is $q^{z_i(m)}$. For $0\leq i\leq n$, we have by definition :
\begin{equation*}
\begin{split}
K_i = \begin{cases} k_i &  \text{ if $0\leq i\leq n-1$,}
\\ (k_0k_1\cdots k_{n-1})^{-2} & \text{ if $i = n$.}
\end{cases}
\end{split}
\end{equation*}
So the eigenvalue of $K_n$ corresponding to $m$ is $q^{-2(z_0(m)+\cdots + z_{n-2}(m) + z_{n-1}(m))}$.
Let $V$ be a representation of $\U_q(\Lo\Glie^\sigma)$ in $\text{Rep}(\U_q(\Lo\Glie^\sigma))$. For $m$ a monomial in the $z_i^{\pm 1}$, let $n_m$ be the multiplicity of $m$ in $\beta(\chi_q^\sigma(V))$. We have :
$$\chi^\sigma(\text{res}^\sigma(V)) = \sum_m n_m z_1^{z_1(m)}z_2^{z_2(m)}\cdots z_{n-1}^{z_{n-1}(m)} z_n^{-z_0(m)-\cdots -z_{n-2}(m) - z_{n-1}(m)}.$$
By the usual invariance of characters by the Weyl group, we get that $\chi^\sigma(\text{res}^\sigma(V))$ is equal to :
\begin{equation*}
\begin{split}
 &\sum_m n_m z_1^{z_1(m)}z_2^{z_2(m)}\cdots z_{n-2}^{z_{n-2}(m)+z_{n-1}(m)} z_{n-1}^{-z_{n-1}(m)}z_n^{-z_1(m)-\cdots -z_{n-2}(m)}
\\=& \sum_m n_m z_1^{z_1(m)}z_2^{z_2(m)}\cdots z_{n-3}^{z_{n-3}(m)+z_{n-2}(m)+z_{n-3}(m)}z_{n-2}^{-z_{n-1}(m)} z_{n-1}^{-z_{n-2}(m)}z_n^{-z_1(m)-\cdots -z_{n-3}(m)}
\\=& \cdots = \sum_m n_m z_1^{-z_{n-1}(m)}z_2^{-z_{n-2}(m)}\cdots z_{n-1}^{-z_1(m)}z_n^{-z_0(m)}.
\end{split}
\end{equation*}
As the Dynkin diagram of type $C$ has no non trivial Dynkin automorphism, the character of a $\U_q(\Glie^\sigma)$-module is invariant by the transformation satisfying for any $1\leq i\leq n$, $z_i\mapsto z_i^{-1}$. We can conclude that $\beta(\chi_q^\sigma(V)) = \chi^\sigma(\text{res}^\sigma(V))$.
\qed

Let $\overline{\chi}^\sigma:\text{Rep}(\overline{\U}_q(\Glie^{\sigma}))\rightarrow \ZZ[z_i^{\pm 1}]_{i\in I_{\sigma}}$ and
$$\overline{\beta} : \ZZ[Z_{i,a}^{\pm 1}]_{a\in\CC^*,i\in I_{\sigma}}\rightarrow \ZZ[z_i^{\pm 1}]_{i\in I_{\sigma}},$$ 
be the ring morphism such that for any $a\in\CC^*$, $i\in I_\sigma$ :
\begin{equation*}
\begin{split}
\overline{\beta}(Z_{i,a}) = \begin{cases}
z_i\text{ if $(\hat{\Glie}^\sigma,i)\neq (A_{2n}^{(2)},0)$,}
\\z_0^2\text{ if $(\hat{\Glie}^\sigma,i) = (A_{2n}^{(2)},0)$.}
\end{cases}
\end{split}
\end{equation*}
It is somewhat analog to a such a morphism considered in \cite{hkott}. Eventually let $\overline{res}^\sigma$ be the restriction map from $\text{Rep}(\U_q(\hat{\Glie}^{\sigma}))$ to $\text{Rep}(\overline{\U}_q(\Glie^{\sigma}))$. 

\begin{prop}\label{dagdeux} The following diagram is commutative :
$$\begin{CD}
\text{Rep}(\U_q(\hat{\Glie}^{\sigma}))  @>{\chi_q^\sigma}>>  \ZZ[Z_{i,a}^{\pm 1}]_{i\in I_\sigma,a\in\CC^*}\\
@VV{\overline{res}^\sigma}V     @VV{\overline{\beta}}V\\
\text{Rep}(\overline{\U}_q(\Glie^{\sigma}))  @>{\overline{\chi}^\sigma}>>  \ZZ[z_i^{\pm 1}]_{i\in I_\sigma}
\end{CD}.$$
\end{prop}

\demo We follows the proof of Proposition \ref{dag}. In the case $A_{2n}^{(2)}$, the eigenvalue of $K_0$ is $q^{z_0(m)} = (q^{\frac{1}{2}})^{2z_0(m)}$, that is why $z_0$ has to be replaced by $z_0^2$ in the definition of $\overline{\beta}$.\qed

\subsection{Examples}\label{exaide}

Let us look at two examples which will be crucial for the following. 

For $\U_q(\Lo sl_2)$, we have $\chi_q(V_a) = Y_a + Y_{aq^2}^{-1}$, where $V_a$ is a two dimensional fundamental representation of $\U_q(\Lo sl_2)$. We have (in fact it is a particular case of Theorem \ref{imchiqu}) : 
\begin{prop}\cite{Fre} We have $Im(\chi_q) = \ZZ[(Y_a + Y_{aq^2}^{-1})]_{a\in\CC^*}$.\end{prop}

For $\U_q^\tau$, we have 
\begin{equation*}
\begin{split}
\chi_q^\tau (V_a) &= Y_{1,a} Y_{2,-a} + Y_{1,-aq^2}Y_{1,aq^4}^{-1}Y_{2,aq^2}Y_{2,-aq^4}^{-1} + Y_{1,- aq^6}^{-1}Y_{2,aq^6}^{-1}
\\&= Z_a + Z_{-aq^2}Z_{aq^4}^{-1} + Z_{-aq^6}^{-1}
\end{split}
\end{equation*}
where $V_a$ is the representation described in section \ref{exfund}. In particular 
\begin{prop} We have $\text{Im}(\chi_q^\tau) = \ZZ[Z_a + Z_{-aq^2}Z_{aq^4}^{-1} + Z_{aq^6}^{-1}]_{a\in\CC^*}$.\end{prop}

Remark : consider the fundamental representations $V_1(a)$, $V_2(a)$ of the untwisted algebra quantum affine algebra $\U_{q^2}(\Lo sl_3)$. We have :
$$\chi_q(V_1(a)) = Y_{1,a} + Y_{1,aq^4}^{-1}Y_{2,aq^2} + Y_{2,aq^6}^{-1},$$
$$\chi_q(V_2(a)) = Y_{2,a} + Y_{2,aq^4}^{-1}Y_{1,aq^2} + Y_{1,aq^6}^{-1}.$$
Consider the ring morphism $\pi : \ZZ[Y_{1,a}^{\pm},Y_{2,a}^{\pm}]_{a\in\CC^*}\rightarrow \ZZ[Z_a^{\pm}]_{a\in\CC^*}$ defined by $\pi(Y_{1,a}) = Z_a$ and $\pi(Y_{2,a}) = Z_{-a}$. Then :
$$\pi(\chi_q(V_1(a))) = \chi_q^\tau (V_a)\text{ , }\pi(\chi_q(V_2(a))) = \chi_q^\tau(V_{-a}).$$
So $\text{Im}(\chi_q^\tau) = \pi(\text{Im}(\chi_q))$. This a particular case of a more general relation between twisted and untwisted case that we will study below (in section \ref{expform}).

Let us look at another example : consider the type $A_4^{(2)}$ and the fundamental representation $V_1(a)$. We have :
$$\chi_q^\sigma(V_1(a)) = Z_{1,a} + Z_{1,aq^2}^{-1}Z_{0,aq} + Z_{0,aq^3}^{-1}Z_{0,-aq^2} + Z_{0,-aq^4}^{-1}Z_{1,-aq^3} + Z_{1,-aq^5}^{-1}.$$
In particular we have 
$$\chi^\sigma(\text{res}^\sigma(V_1(a))) = (z_1 + z_2 z_1^{-1} + z_1 z_2^{-1} + z_1^{-1})  + 1 = \chi^\sigma(V(\Lambda_1)\oplus V(0)),$$
$$\overline{\chi}^\sigma(\overline{\text{res}}^\sigma(V_1(a))) = z_1 + z_1^{-1}z_0^2 + 1 + z_1z_0^{-2} + z_1^{-1} = \overline{\chi}^\sigma(\overline{V}(\Lambda_1)),$$
where the $V(\lambda)$ (resp. $\overline{V}(\lambda)$) are the simple representations of $\U_q(\Glie^\sigma)$ (resp. $\overline{\U}_q(\Glie^\sigma)$). In particular $\text{res}^\sigma(V_1(a))$ is not simple, but $\overline{\text{res}}^\sigma(V_1(a))$ is simple. These branching rules were known as $V_1(a)$ is a fundamental representation, but more general branching rules will be proved in Section \ref{branch}.

This representation $V_1(a)$ has a also crystal basis \cite{kas} and the crystal graph can be computed (see \cite{kkmmnn}) :

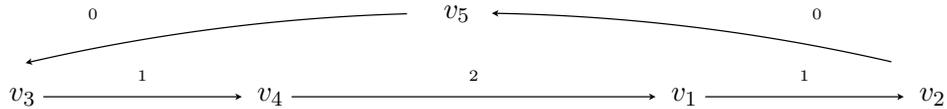
\begin{figure}[htbp]
\centering
\psset{xunit=.55mm,yunit=.55mm,runit=.55mm}
\psset{linewidth=0.3,dotsep=1,hatchwidth=0.3,hatchsep=1.5,shadowsize=1}
\psset{dotsize=0.7 2.5,dotscale=1 1,fillcolor=black}
\psset{arrowsize=1 2,arrowlength=1,arrowinset=0.25,tbarsize=0.7 5,bracketlength=0.15,rbracketlength=0.15}
\begin{pspicture}(0,0)(220,50)
\rput(105,40){$v_5$}
\rput{0}(106,-500){\parametricplot[arrows=->]{78}{89.15}{ t cos 500.26 mul t sin 540.26 mul }}
\rput{0}(106,-500){\parametricplot[arrows=<-]{102.15}{91.5}{ t cos 500.26 mul t sin 540.26 mul }}
\rput(192,40){$\scriptscriptstyle 0$}
\rput(17,40){$\scriptscriptstyle 0$}
\rput(0,20){$v_3$}
\rput(60,20){$v_4$}
\psline{->}(5,20)(53,20)
\psline{->}(65,20)(153,20)
\rput(160,20){$v_1$}
\psline{->}(165,20)(213,20)
\rput(220,20){$v_2$}
\rput(29,25){$\scriptscriptstyle 1$}
\rput(109,25){$\scriptscriptstyle 2$}
\rput(189,25){$\scriptscriptstyle 1$}
\end{pspicture}
\caption{(Type $A_4^{(2)}$) the crystal of $V_1(a)$}
\label{fig:A22}
\end{figure}

By erasing the $0$-arrows we get the $\tilde{\U}_q(\Glie^\sigma)$-crystal graph which is not connected.

In fact we can have the identification (see the discussion in the last section of \cite{hn}) : 
$$(V_1(a))_{Z_{1,a}} = \CC.v_1\text{ , }(V_1(a))_{Z_{1,aq^2}^{-1}Z_{2,aq}} = \CC.v_2\text{ , }(V_1(a))_{Z_{2,aq^3}^{-1}Z_{2,-aq^2}} = \CC.v_5,$$ 
$$(V_1(a))_{Z_{2,-aq^4}^{-1}Z_{1,-aq^3}} = \CC.v_3\text{ , }(V_1(a))_{Z_{1,-aq^5}^{-1}} = \CC.v_4.$$
To get the graph associated with twisted $q$-characters as in \cite{Fre}, we should add a $2$-arrow from $v_2$ to $v_5$ and from $v_5$ to $v_2$ (see section \ref{fundhe} for the general definition of such a graph).

Note that by renumbering the nodes $(0,1,2)\leftrightarrow (2,1,0)$ we get the crystal graph of $V_1(a)$ viewed as a representation of type $A_4^{(2)^\dagger}$ (see \cite{hn} for example). By erasing the $0$-arrows we get the $\overline{\U}_q(\Glie^\sigma)$-crystal graph which is connected.

\begin{figure}[htbp]
\centering
\psset{xunit=.55mm,yunit=.55mm,runit=.55mm}
\psset{linewidth=0.3,dotsep=1,hatchwidth=0.3,hatchsep=1.5,shadowsize=1}
\psset{dotsize=0.7 2.5,dotscale=1 1,fillcolor=black}
\psset{arrowsize=1 2,arrowlength=1,arrowinset=0.25,tbarsize=0.7 5,bracketlength=0.15,rbracketlength=0.15}
\begin{pspicture}(0,0)(220,30)
\rput(0,20){$v_1$}
\rput(60,20){$v_2$}
\psline{->}(5,20)(53,20)
\psline{->}(65,20)(103,20)
\rput(110,20){$v_5$}
\psline{->}(115,20)(153,20)
\rput(160,20){$v_3$}
\psline{->}(165,20)(213,20)
\rput(220,20){$v_4$}
\rput{0}(110.26,520){\parametricplot[arrows=<-]{-102.15}{-77.85}{ t cos 520.26 mul t sin 520.26 mul }}
\rput(29,25){$\scriptscriptstyle 1$}
\rput(84,25){$\scriptscriptstyle 2$}
\rput(134,25){$\scriptscriptstyle 2$}
\rput(189,25){$\scriptscriptstyle 1$}
\rput(110,3){$\scriptscriptstyle 0$}
\end{pspicture}
\caption{(Type $A_4^{(2)^\dagger}$) the crystal of $V_1(a)$}
\end{figure}
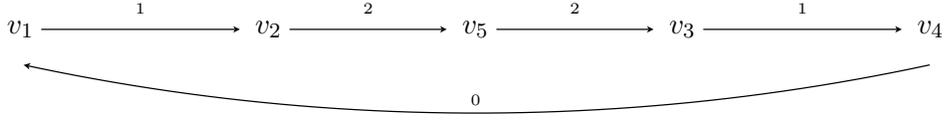

\subsection{Combinatorics of twisted $q$-characters} The aim of this section is to compute $\text{Im}(\chi_q^\sigma)$, and so to prove an analog of Theorem \ref{imchiqu} for twisted $q$-characters (Theorem \ref{imchiq}).

Let $\epsilon_i$ defined by $\epsilon_i=M$ if $i=\sigma(i)$, and $\epsilon_i=1$ if $i\neq \sigma(i)$.

For $i,j\in I$ and $m\neq 0$ denote 
\begin{equation*}
F_{i,j}(m) = \sum_{k = 1\cdots M}[\epsilon_i mC_{i,\sigma^k(j)}/d_{\overline{i}}]_{q_{\overline{i}}}\omega^{\epsilon_i mk}.
\end{equation*} 
In fact this definition is set such that : 
$$[h_{i,\epsilon_i m},x_{j,r}^{\pm}]=\pm \frac{1}{\epsilon_i m} F_{i,j}(m) x_{j,\epsilon_im+r}^{\pm}.$$
We remind that this is the relation (\ref{eqdquatre}) proved in \cite[Theorem 5.3.2]{da}.

As we choose a class of representative of $I_{\sigma}$, it makes sense to consider the matrix $F(m) = (F_{i,j}(m))_{i,j\in I_{\sigma}}$. 

\begin{lem}\label{condgen} 
For a generic parameter $q$, for any $m\in\ZZ-\{0\}$ the matrix $F(m)$ is invertible and the diagonal coefficients of the inverse matrix $\tilde{F}(m)$ of $F(m)$ are non zero.
\end{lem}

Here $q$ generic means that there is a finite set of algebraic equations such that for any $m\in\ZZ$, $q^m$ is not a solution of one of these equations.

\demo Let $D_m(q)$ be the determinant of 
$$F'(m) = (\prod_{i\in I_\sigma}(q_{\overline{i}} - q_{\overline{i}}^{-1}))F(m)$$ 
and for $i\in I_{\sigma}$, let $D_{i,m}(q)$ be the determinant of 
$$(\prod_{j\in I_\sigma,j\neq i}(q_{\overline{j}} - q_{\overline{j}}^{-1}))(F_{j,k}(m))_{j,k\in I_{\sigma},j\neq i,k\neq i}.$$ 
Then there are $D(z)\in \CC[z^{\pm 1}]$ and $D_i(z)\in\CC[z^{\pm 1}]$ such that for any $m\in\ZZ$ :
$$D_m(q) = D(q^m)\text{ and } D_{i,m}(q) = D_i(q^m).$$ 
So it suffices to prove that $D(z)\neq 0$ and $D_i(z)\neq 0$. By looking at the degree in $q$ of the coefficients of $F'(m)$, we notice that for each row $i$ the maximal degree is on the diagonal where $q^{\epsilon_i 2m}$ appears. The coefficient of $q^{\epsilon_i 2m}$ is $1$ if $i\neq \sigma(i)$ and $M$ if $i = \sigma(i)$. So the maximal degree in $P(z)$ is $2\sum_{i\in I_\sigma}\epsilon_i$ and the coefficient of $z^{2\sum_{i\in I_\sigma}\epsilon_i}$ in $P(z)$ is $\prod_{i\in I_\sigma}\epsilon_i$, so $P(z)\neq 0$. 

A similar argument gives the result for $D_i(z)$.\qed

In the untwisted case, it suffices to suppose that $q$ is not a root of unity to have a similar result (see \cite{Fre, Fre2}). For the untwisted case it is not true in general : indeed for type $D_4^{(3)}$ we have 
$$D_m(q) = 3(q^m-q^{-m})(q^{3m}-q^{-3m})(q^{4m}+q^{-4m}+q^{2m}+q^{-2m}-1-j^m-j^{2m})$$
and so for $m\in\ZZ$, we have $D_{3m}(q)=0$ if $Q = q^{6m}+q^{-6m}$ satisfies the equation $Q^2+Q-5=0$. We get for example $q = (2(21)^{\frac{1}{2}}-6)^{\frac{1}{2}}/2$ which is obviously not a root of one. 

\begin{lem} 
For types $A_{2n}^{(2)}$, the statement of Lemma \ref{condgen} it satisfied if $q$ is not a root of unity. 
\end{lem}

\demo
Let us denote $D_m^{(n)}(q)=D_m(q)$ and $D_{i,m}^{(n)}(q)=D_{i,m}(q)$ for type $A_{2n}^{(2)}$. From classical result we have 
$$D_{n,m}^{(n)}(q) = \frac{q^{mn} - q^{-mn}}{q^m - q^{-m}}.$$ 
Moreover by developing the determinant we have :
\begin{equation*}
\begin{split}
D_m^{(n)}(q) &= (q^m+q^{-m}- (-1)^m)D_{n,m}^{(n)}(q) - D_{n-1,m}^{(n-1)}(q) 
\\&= \frac{(1+(-q^{-2n-1})^m)((q^{n+1})^m - (-q^n)^m)}{q^m-q^{-m}},
\end{split}
\end{equation*}
whose roots are roots of unity. We can conclude as 
$$D_{i,m}^{(n)}(q) = D_{i,m}^{(i)}(q)D_m^{(n-i)}(q).$$
\qed

Note that $\epsilon_i$ is crucial in the definition of $F(m)$ : indeed if $\sigma(i) = i$ and $m\neq 0 \text{ mod }M$, then the determinant of the matrix defined with $1$ instead of $\epsilon_i$ is $0$. In fact this follows from the fact that in this case $h_{i,m}  = 0$. In particular we do not loose information by using the $\epsilon_i$.

In the following we suppose that $q$ is generic, that is to say that the conditions of Lemma \ref{condgen} are satisfied.

$\tilde{F}(m) = (\tilde{F}_{i,j}(m))_{i,j\in I_{\sigma}}$ is the inverse matrix of $F(m)$. For $m\neq 0$ and $i\in I_{\sigma}$ let us define 
$$\tilde{h}_{i,m} = \sum_{j\in I_\sigma}\tilde{F}_{i,j}(m) h_{j,\epsilon_j m}.$$ 
Remark : as for $i,j\in I$, $F_{\sigma(i), j}(m) = \omega^{\epsilon_i m} F_{i,j}(m)$, we have $\tilde{F}_{i,\sigma(j)}(m) = \omega^{-\epsilon_j m} \tilde{F}_{i,j}(m)$. Moreover we have $\tilde{h}_{\sigma(j),\epsilon_j m} = \omega^{\epsilon_j m} \tilde{h}_{j,\epsilon_j m}$, so $\tilde{F}_{i,j}(m) h_{j,\epsilon_j m}$ does not depend of the choice of the representative $j$. Moreover $\CC.\tilde{h}_{i,m}$ does not depend of the choice of the representatives of $I_{\sigma}$.

 Let $i\in I_{\sigma}$ and
$$\U_q(\Hlie^{\sigma})_i^\perp = \{K\in \CC[k_i^{\pm}]_{i\in I} | [\hat{\U}_i,K] = 0\}.$$ 
Let $\U_q(\Lo\Hlie^\sigma)_i^\perp$ be the subalgebra of $\U_q(\Lo\Hlie)$ generated by $\U_q(\Hlie^{\sigma})_i^\perp$ and the $\tilde{h}_{j,m}$, $j\in I_\sigma$, $j\neq i$, $m\neq 0$.

\noindent For $\gamma$ an $l$-weight, we denote by $\tau_i(\gamma)$ the couple $(\gamma_i(z), \gamma_i')$ where $\gamma_i' : \U_q(\Lo\Hlie^\sigma)_i^\perp\rightarrow\CC$ is the algebra morphism induced by $\gamma$. $\tau_i$ can be extended to a map : 
$$\tau_i : \ZZ[Z_{j,a}]_{j\in I_{\sigma},a\in\CC^*}\rightarrow \ZZ[Z_{i,a}]_{a\in\CC^*}\otimes (\U_q(\Lo\Hlie^\sigma)_i^\perp)^*.$$ 
Note that $\tau_i$ is clearly injective. We set $u_i(\tau_i(m)) = z_i(m)$, so $u_i$ is well-defined on $\text{Im}(\tau_i)$.

\begin{lem}\label{aidedeux} Let $V\in\text{Rep}(\U_q(\Lo\Glie))$ and consider a decomposition : $$\tau_i(\chi_q^\sigma(V))=\underset{r}{\sum} P_r\otimes Q_r,$$ 
where all monomials $Q_r\in (\U_q(\Lo\Hlie^\sigma)_i^\perp)^*$ are distinct. Then the
$\hat{\U}_i$-module $V$ is isomorphic to a direct sum $V = \underset{r}{\bigoplus}V_r$ where the twisted $q$-character of the $\hat{\U}_i$-module $V_r$ is $P_r$.\end{lem}

\noindent There are previous analog results for untwisted cases (\cite[Lemma 3.4]{Fre2}) and untwisted general quantum affinizations (\cite[Lemma 5.10]{her04}).

\demo  For $i\in I$, we have 
$$h_{i,\epsilon_i m} = \sum_{j\in I_{\sigma}}F_{i,j}(m)\tilde{h}_{j, m}.$$ 
As moreover for $i\in I_\sigma$, $\tilde{F}_{i,i}(m)\neq 0$, we have :
$$\bigoplus_{j\in I_{\sigma}}\CC h_{i,\epsilon_i m} = \bigoplus_{j\in I_{\sigma}}\CC \tilde{h}_{j, m} = \CC h_{i,\epsilon_i m} \oplus\bigoplus_{j\in I_{\sigma},j\neq i}\CC\tilde{h}_{j,m}.$$
From the relation (\ref{eqdquatre}) proved in \cite[Theorem 5.3.2]{da}, for $j\neq \overline{i}$ we have $[\hat{\U}_i,\tilde{h}_{j,m}] = 0$. We have $[\hat{\U}_i,\U_q(\Lo\Hlie^\sigma)_i^\perp] = 0$ and $\U_q(\Lo\Hlie)$ is generated by $\U_q(\Lo\Hlie^\sigma)_i^\perp$ and the $k_i^{\pm 1}, h_{i,r}$, $r\neq 0$. As $\U_q(\Lo\Hlie^\sigma)_i^\perp$ is commutative, it suffices to decompose $V$ in common eigensubspaces for all elements of $\U_q(\Lo\Hlie^\sigma)_i^\perp$.
\qed

For $i\in I_{\sigma}$, we denote $j\sim i$ if $j\in I_{\sigma}$ and $C_{i,j} = -1$.

\noindent For $i\in I_{\sigma}, a\in\CC^*$, let us define elements $A_{i,a}\in\ZZ[Z_{j,b}]_{j\in I_{\sigma},b\in\CC^*}$ analogs of simple root for monomials. In \cite{Fre} the definition of the $A_{i,a}$ in the untwisted case was  given. For the twisted case we have a new definition with modifications. Indeed for $a\in\CC^*,i\in I_\sigma$ we set : 
\\If $C_{i,\sigma(i)} = 2$ :
$$A_{i,a} = Z_{i,aq^M}Z_{i,aq^{-M}}\times \prod_{j\sim i|j = \sigma(j)}Z_{j,a}^{-1}\times \prod_{j\sim i| j\neq \sigma(j)}(\prod_{a'\in\CC^*|(a')^M = a}Z_{j,a'}^{-1}).$$
If $C_{i,\sigma(i)} = 0$ :
$$A_{i,a} = Z_{i,aq}Z_{i,aq^{-1}}\times \prod_{j\sim i|j = \sigma(j)}Z_{j,a^M}^{-1}\times \prod_{j\sim i| j\neq \sigma(j)}Z_{j,a}^{-1}.$$
If $C_{i,\sigma(i)} = -1$ :
$$A_{i,a} = Z_{i,aq}Z_{i,aq^{-1}}Z_{i,-a}^{-1}\times \prod_{j\sim i}Z_{j,a}^{-1}.$$
The motivation for the definition of the $A_{i,a}$ will appear in the proof of Theorem \ref{imchiq}.
\\For $i\in I_\sigma$ such that $C_{i,\sigma(i)} = 2$ define
$$\mathcal{K}_i = \ZZ[Z_{i,a}(1 + A_{i,aq^M}^{-1})]_{a\in\CC^*}\times\ZZ[Z_{j,a}^{\pm 1}]_{j\in I_{\sigma}, j\neq i,a\in\CC^*}.$$
For $i\in I_\sigma$ such that $C_{i,\sigma(i)} = 0$ define
$$\mathcal{K}_i = \ZZ[Z_{i,a}(1 + A_{i,aq}^{-1})]_{a\in\CC^*}\times\ZZ[Z_{j,a}^{\pm 1}]_{j\in I_{\sigma}, j\neq i,a\in\CC^*}.$$
For $i\in I_\sigma$ such that $C_{i,\sigma(i)} = -1$ define
$$\mathcal{K}_i = \ZZ[Z_{i,a}(1 + A_{i,aq}^{-1} + A_{i,aq}^{-1}A_{i,-aq^2}^{-1})]_{a\in\CC^*}\times\ZZ[Z_{j,a}^{\pm 1}]_{j\in I_{\sigma}, j\neq i,a\in\CC^*}.$$

Let us state and prove the analog of Theorem \ref{imchiqu} for the twisted cases, which is the main result of this section :

\begin{thm}\label{imchiq} We have $\text{Im}(\chi_q^\sigma) = \bigcap_{i\in I_{\sigma}} \mathcal{K}_i$.\end{thm}

\demo First let us prove that $\text{Im}(\chi_q^\sigma) \subset \bigcap_{i\in I_{\sigma}} \mathcal{K}_i$.

We start with the decomposition of Lemma \ref{aidedeux}. On each component $V_r$, $\gamma_i'$ is constant. For $\gamma$ an $l$-weight satisfying $\gamma_i' = 0$, we have for $l\in I_{\sigma}$ and $m\neq 0$, $$F_{i,i}(m)\gamma(h_{l,\epsilon_lm}) = F_{l,i}(m) \gamma(h_{i,\epsilon_im}).$$ 
So if we have a given $\gamma_i(z)$ we can determine uniquely and explicitly a corresponding $\gamma(z)$.  So it suffices to prove that the monomials $A_{i,a}$ defined above are the monomials corresponding by this process to the $A_{i,a}$ for the subalgebra $\hat{\U}_i$. Let us check it case by case.

Let $i,j\in i_{\sigma}$ such that $\overline{i}\neq\overline{j}$ and $C_{i,j} = - 1$. 

Remark : we can not have simultaneously $C_{i,\sigma(i)} = 2$ and $C_{j,\sigma(j)} = -1$.

We study the different cases corresponding to the values of $C_{i,\sigma(i)}$ in $\{-1, 0, 2\}$.

$\bullet$ $C_{i,\sigma(i)} = 2$. We have $F_{i,i}(m) = M \frac{q^{2mM} - q^{-2mM}}{q_{\overline{i}} - q_{\overline{i}}^{-1}}$. From the $sl_2$-case we set for $m\in\ZZ - \{0\}$ :
$$\gamma(h_{i,mM}) = \frac{q^{mM} - q^{-mM}}{m (q^M - q^{-M})} (- (aq^{-M})^m - (aq^M)^m).$$
If $C_{j,\sigma(j)} = 2$ : we have in this situation $d_{\overline{i}} = d_{\overline{j}}$, and :
$$F_{j,i} (m) = M \frac{q^{-mM} - q^mM}{q_{\overline{j}} - q_{\overline{j}}^{-1}}.$$
So $\gamma(h_{j,mM})$ is equal to :
\begin{equation*}
\begin{split}
 &\frac{q_{\overline{i}} - q_{\overline{i}}^{-1}}{M(q^{2mM} - q^{-2mM})}M \frac{q^{-mM} - q^{mM}}{q_{\overline{j}} - q_{\overline{j}}^{-1}} 
\times\frac{q^{mM} - q^{-mM}}{m (q^M - q^{-M})} (- (aq^{-M})^m - (aq^M)^m)
\\= &\frac{q^{-mM} - q^{mM}}{q^{2mM} - q^{-2mM}} \frac{q^{mM} - q^{-mM}}{m (q^M - q^{-M})} a^m(- q^{-mM} - q^{mM})
=\frac{q^{mM} - q^{-mM}}{m (q^M - q^{-M})} a^m.
\end{split}
\end{equation*}
This corresponds to $Z_{j,a}$.

If $C_{j,\sigma(j)} = 0$ : we have in this situation $d_{\overline{i}} = M$ and $d_{\overline{j}} = 1$. We denote $\delta_m^{[M]}$ equal to $1$ if $m=0\text{ mod}[M]$ and equal to $0$ otherwise. We have :
$$F_{j,i} (m) = M \frac{q^{-m} - q^m}{q - q^{-1}}\delta_{m}^{[M]}.$$
Let $a'\in\CC^*$ such that $(a')^M = a$. We have
\begin{equation*}
\begin{split}
\gamma(h_{j,m}) &= \frac{q^M - q^{-M}}{M(q^{2m} - q^{-2m})}M \frac{q^{-m} - q^{m}}{q - q^{-1}} 
\times \frac{M(q^{m} - q^{-m})}{m (q^M - q^{-M})} (a')^m(- q^m - q^m)\delta_{m}^{[M]}
\\&=  \frac{(q^{m} - q^{-m})}{m (q - q^{-1})} M(a')^m \delta_{m}^{[M]} = \frac{(q^{m} - q^{-m})}{m (q - q^{-1})} \sum_{b\in\CC^*| b^M = a}b^m .
\end{split}
\end{equation*}
This corresponds to $\prod_{\{b\in\CC^*| b^M = a\}} Z_{j,b}$.

The case $C_{j,\sigma(j)} = -1$ can not occur.

$\bullet$ $C_{i,\sigma(i)} = 0$ : We have $F_{i,i}(m) = \frac{q^{2m} - q^{-2m}}{q_{\overline{i}} - q_{\overline{i}}^{-1}}$. From the $sl_2$-cases we set for $m\in\ZZ-\{0\}$ : 
$$\gamma(h_{i,m}) = \frac{q^{m} - q^{-m}}{m (q - q^{-1})} (- (aq^{-1})^m - (aq)^m).$$

If $C_{j,\sigma(j)} = 2$ : we have in this situation $d_{\overline{i}} = 1$ and $d_{\overline{j}} = M$. So we have
\begin{equation*}
\begin{split}
F_{j,i} (m) &= \frac{q^{-mM} - q^{mM}}{q^M - q^{-M}}M, 
\\\gamma(h_{j,mM}) &= M\frac{q^{-mM} - q^{mM}}{q^M - q^{-M}} \frac{q - q^{-1}}{q^{2mM} - q^{-2mM}}
\times\frac{q^{mM} - q^{-mM}}{mM (q - q^{-1})} (- (aq^{-1})^{mM} - (aq)^{mM})
\\&=\frac{q^{mM} - q^{-mM}}{m(q^M - q^{-M})}(a^M)^m.
\end{split}
\end{equation*}
This corresponds to $Z_{j,a^M}$.

If $C_{j,\sigma(j)} = 0$ : we have in this situation $d_{\overline{i}} = d_{\overline{j}} = 1$. We have :
\begin{equation*}
\begin{split}
F_{j,i} (m) &= \frac{q^{-m} - q^m}{q - q^{-1}}.
\\\gamma(h_{i,m}) &= \frac{q^{-m} - q^m}{q - q^{-1}}
\frac{q^{m} - q^{-m}}{m (q - q^{-1})} (- (aq^{-1})^m - (aq)^m)\frac{q - q^{-1}}{q^{2m} - q^{-2m}}
=a^m \frac{q^m - q^{-m}}{m(q - q^{-1})}.
\end{split}
\end{equation*}
This corresponds to $Z_{j,a}$.

If $C_{j,\sigma(j)} = -1$ : we have in this situation $d_{\overline{j}} = \frac{1}{2}$ and $d_{\overline{i}} = 1$. So we have
\begin{equation*}
\begin{split}
F_{j,i} (m) &= \frac{q^{-m} - q^m}{q^{\frac{1}{2}} - q^{-\frac{1}{2}}}.
\\\gamma(h_{i,m}) &= \frac{q^{-m} - q^m}{q^{\frac{1}{2}} - q^{-\frac{1}{2}}}
\frac{q^{m} - q^{-m}}{m (q - q^{-1})} (- (aq^{-1})^m - (aq)^m)\frac{q - q^{-1}}{q^{2m} - q^{-2m}}
\\&=a^m \frac{(q^{\frac{1}{2}})^m - (q^{\frac{1}{2}})^{-m}}{m(q^{\frac{1}{2}} - q^{-\frac{1}{2}})}.
\end{split}
\end{equation*}
This corresponds to $Z_{j,a}$.

$\bullet$ $C_{i,\sigma(i)} = - 1$. The cases $C_{j,\sigma(j)} = 2$ and $C_{j,\sigma(j)} = -1$ do not occur. So $C_{j,\sigma(j)} = 0$. We have $d_{\overline{i}} = \frac{1}{2}$ and $d_{\overline{j}} = 1$. We have 
\begin{equation*}
\begin{split}
F_{i,i}(m) &= [4m]_{q^{\frac{1}{2}}} + (-1)^m [-2m]_{q^{\frac{1}{2}}} 
\\&= \frac{q^{2m} - q^{-2m} + (-1)^mq^{-m} - (-1)^m q^m}{q^{\frac{1}{2}} - q^{-\frac{1}{2}}}
= \frac{ (q^m - q^{-m}) (q^{m} + q^{-m} - (-1)^m)}{q^{\frac{1}{2}} - q^{-\frac{1}{2}}}.
\end{split}
\end{equation*}
From the case of $\U_q^\tau$, we set 
\begin{equation*}
\begin{split}
\gamma(h_{i,m}) &= \frac{(q^{\frac{1}{2}})^{2m} - (q^{\frac{1}{2}})^{-2m}}{m (q^{\frac{1}{2}} - q^{-\frac{1}{2}})} (- (aq)^m - (aq^{-1})^m + (-a)^m) 
\\&= \frac{q^m - q^{-m}}{m (q^{\frac{1}{2}} - q^{-\frac{1}{2}})} a^m (- q^m - q^{-m} + (-1)^m).
\end{split}
\end{equation*}
So $F_{j,i} (m) = \frac{q^{-m} - q^m}{q - q^{-1}}$ and 
\begin{equation*}
\begin{split}
\gamma(h_{j,m}) = &\frac{q^m - q^{-m}}{m (q^{\frac{1}{2}} - q^{-\frac{1}{2}})} a^m (- q^m - q^{-m} + (-1)^m) 
\\&\times\frac{q^{-m} - q^m}{q - q^{-1}}\frac{q^{\frac{1}{2}} - q^{-\frac{1}{2}}}{(q^m - q^{-m}) (q^{m} + q^{-m} - (-1)^m)}
= \frac{q^m - q^{-m}}{m(q - q^{-1})} a^m.
\end{split}
\end{equation*}
This corresponds to $Z_{j,a}$.

Now we prove the other inclusion. Here we consider the usual order on the integral weight lattice of $\Glie^\sigma$. For $\chi\in\bigcap_{i\in I_\sigma}\mathcal{K}_i$ non equal to $0$, an highest weight element is clearly dominant in the sense of monomials. So a non zero element of $\bigcap_{i\in I_\sigma}\mathcal{K}_i$ has at least one dominant monomial. As the twisted $q$-characters of simple modules give elements in $\text{Im}(\chi_q^\sigma)$ of the form 
$$m \text{ + elements of lower weight},$$ 
it suffices to prove that there is a finite number of possible weights corresponding to dominant monomials lower than $m$. For such an $m'$, we have 
$$\beta(m) = \beta(m') + \sum_{i\in I_\sigma} a_i \alpha_i\text{ where for any $i\in I_\sigma$, $a_i\geq 0$.}$$ 
As the Cartan matrix of $\Glie^\sigma$ is of finite type and so invertible, there are $b_i\in\QQ$ such that $\beta(m)=\sum_{i\in I_\sigma} b_i\alpha_i$. So we have 
$$\beta(m') = \sum_{i\in I_\sigma} (b_i - a_i)\alpha_i.$$ 
From \cite[Theorem 4.3]{kac}, we have $b_i - a_i\geq 0$ for any $i$, and so $0\leq a_i\leq b_i$. As the $b_i$ are fixed, there is only a finite number of possible $a_i\in\ZZ$.\qed

It is possible to define corresponding twisted screening operators such that $\mathcal{K}_i = \text{Ker}(S_i)$ (see \cite{Fre} for the untwisted case).

As explained in the proof of Theorem \ref{imchiq}, we have :

\begin{cor} A non-zero element of $\text{Im}(\chi_q^\sigma)$ has at least one dominant monomial. \end{cor}

In particular an element of $\text{Im}(\chi_q^\sigma)$ with a unique dominant monomial is uniquely determined by this monomial.

\noindent For $j\in I_\sigma$ and $m\in A$ denote $m^{(j)}=\underset{a\in\CC^*}{\prod}Z_{j,a}^{z_{j,a}(m)}$. For $a\in\CC^*$ consider $A_{j,a}^{j, \pm}=(A_{j,a}^{\pm})^{(j)}$. Define 
$$\mu_j^I:\ZZ[A_{j,a}^{J, \pm}]_{a\in\CC^*}\rightarrow \ZZ[A_{j,a}^{\pm}]_{a\in\CC^*}$$ 
as the ring morphism such that $\mu_j^I(A_{j,a}^{j, \pm})=A_{j,a}^{\pm}$. For $m\in B_j$, denote $L^j(m^{(j)})$ the (twisted) $q$-characters defined for the sub (twisted) quantum affine algebra $\hat{\U}_j$. Define :
$$L_j(m)=m \mu_j^I((m^{(j)})^{-1}L_j(m^{(j)})).$$
We have :

\begin{prop}\label{jdecomp} 
For a module $V\in\text{Rep}(\U_q(\Lo \Glie^\sigma))$ and $j\in I_\sigma$, there is unique
decomposition in a finite sum : 
$$\chi_q^\sigma(V)=\underset{m'\in B_j}{\sum}\lambda_j(m')L_j(m').$$
Moreover for all $m'$, $\lambda_j(m')\geq 0$.
\end{prop}

See \cite[Proposition 3.9]{her05} for the analog result in the untwisted case. The proof in the twisted case is analog by using Theorem \ref{aidedeux} and Lemma \ref{imchiq}.

\subsection{Additional definitions}\label{addef}

As the $A_{i,a}^{-1}$ are algebraically independent (as the matrices $F(m)$ are invertible), for $M$ a 
product of $A_{i,a}^{-1}$ we can define $v_{i,a}(M)\geq 0$ by $M=\underset{i\in I_\sigma, 
a\in\CC^*}{\prod}A_{i,a}^{-v_{i,a}(m)}$. We put $v(M)=\underset{i\in I_\sigma, a\in\CC^*}{\sum}v_{i,a}(m)$. 

\noindent We denote $m\leq m'$ if $m'm^{-1}$ is a product of $A_{i,a}$ ($i\in I, a\in\CC^*$). This partial ordering is called the partial ordering in the sense of monomials. 

\noindent We remind the maps $\beta$, $\overline{\beta}$ defined in Section \ref{fprop}. For $m$ a monomial, we can consider $\beta(m)$ (resp. $\overline{\beta}(m)$) as an element of the weight lattice of $\U_q(\Glie^\sigma)$ (resp. $\overline{\U}_q(\Glie^\sigma)$.

\noindent For $\lambda$ in the weight lattice of $\overline{\U}_q(\Glie^\sigma)$, we set 
$$v(\lambda) = -\lambda(\Lambda_1^{\vee}+...+\Lambda_n^{\vee}).$$ 
For a product $\mathcal{P}$ of $A_{i,a}^{-1}$, we have $v(\mathcal{P})=v(\overline{\beta}(\mathcal{P}))$. So the map $v$ can be extended to any monomial.

For the untwisted case, the notion of right-negative monomial was introduced in \cite{Fre2}. The definition in the untwisted case is modified :

\begin{defi}\label{monomrn} Suppose that $\hat{\Glie}^\sigma$ is not of type $A_{2n}^{(2)}$. A monomial $m\neq 1$ is said to be right-negative if for all $a\in\CC^*$, 
for 
$$L=\text{max}\{l\in\ZZ | \exists i\in I_\sigma,\exists r\in\ZZ, z_{i,(a\omega^rq^L)^{d_i}}(m)\neq 0\},$$ 
we have $\forall j\in I_\sigma,\forall r\in\ZZ$, $z_{j,(a\omega^rq^L)^{d_i}}(m) \leq 0$.

Suppose that $\hat{\Glie}^\sigma$ is of type $A_{2n}^{(2)}$. A monomial $m\neq 1$ is said to be right-negative if for all $a\in\CC^*$, 
for 
$$L=\text{max}\{l\in\ZZ | \exists i\in I_\sigma, z_{i,a q^L}(m)\neq 0\text{ or }z_{i,-a q^L}(m)\neq 0\},$$ 
we have $\forall j\in I_\sigma$, $z_{j,a q^L}(m)\leq 0$ and $z_{j,-a q^L}(m)\leq 0$.\end{defi}

\noindent Note that a right-negative monomial is not dominant. As in \cite{Fre2}, we have :

\begin{lem}\label{rn} 1) For $i\in I_\sigma, a\in\CC^*$, $A_{i,a}^{-1}$ is right-negative.

2) A product of right-negative monomials is right-negative.

3) If $m$ is right-negative, then $m'\leq m$ implies that $m'$ is right-negative.\end{lem}

Let $k\geq 0, a\in\CC^*, i\in I_{\sigma}$. We set
\begin{equation*}
m_{k,a}^{(i)} = 
\begin{cases}
\underset{s=1\cdots k}{\prod}Z_{i,aq_i^{2s-2}}\text{ for $\hat{\Glie}^\sigma$ not of type $A_{2n}^{(2)}$,}
\\\underset{s=1\cdots k}{\prod}Z_{i,aq^{2s-2}}\text{ for $\hat{\Glie}^\sigma$ of type $A_{2n}^{(2)}$.}
\end{cases}
\end{equation*}

\begin{defi} 
For $k\geq 0, a\in\CC^*, i\in I_{\sigma}$, the Kirillov-Reshetikhin module $W_{k,a}^{(i)}$ is the simple module corresponding to the monomial $m_{k,a}^{(i)}$.
\end{defi}

For the untwisted case, the definition of Kirillov-Reshetikhin modules is analog where the highest monomials $\underset{s=1\cdots k}{\prod}Y_{i,aq^{r_i(2s-2)}}$ are used.

\noindent For $i\in I_\sigma$ and $a\in\CC^*$, $W_{1,a}^{(i)} = V_{i,a}$ is called a fundamental representation (this coincides with the definition of Section \ref{fdeptq}). 

\noindent Two dominant monomials $m_1=m_{k_1,a_1}^{(i)}$, $m_2=m_{k_2,a_2}^{(i)}$ are said to be in special position if the monomial $m_3=\underset{a\in\CC^*}{\prod}Z_{i,a}^{\text{max}(z_{i,a}(m_1),z_{i,a}(m_2))}$ is of the form $m_3=m_{k_3,a_3}^{(i)}$ and $m_3\neq m_1, m_3\neq m_2$.

\noindent A normal writing of a dominant monomial $m$ is a product decomposition $$m=\underset{i=1,...,L}{\prod}m_{k_l,a_l}^{(i_l)}$$ 
such that for $l\neq l'$, if $i_l=i_{l'}$ then $m_{k_l,a_l}^{(i_l)}$, $m_{k_{l'},a_{l'}}^{(i_{l'})}$ are not in special position. Any dominant monomial has a unique normal writing up to permuting the monomials (see \cite{Cha2}).

\noindent It follows from the study of the representations of $\U_q(\Lo sl_2)$ in \cite{Cha0, Cha, Fre} that :

\begin{prop}\label{aidesldeux} Suppose that $\Glie=sl_2$. 

(1) $W_{k,a}$ is of dimension $k+1$ and :
$$\chi_q(W_{k,a})=m_{k,a}(1+A_{aq^{2k-1}}^{-1}(1+A_{aq^{2(k-1)-1}}^{-1}(1+...(1+A_{aq^{2-1}}^{-1}))...).$$

(2) For $m$ a dominant monomial and $m=m_{k_1,a_1}...m_{k_l,a_l}$ a normal writing we have :
$$L(m)\simeq W_{k_1,a_1}\otimes ...\otimes W_{k_l,a_l}.$$
\end{prop}

We extend the notion of special modules \cite{Nab} to the twisted case :

\begin{defi} A $\U_q(\Lo\Glie^\sigma)$-module is said to be special if his twisted $q$-character has a unique dominant 
monomial.\end{defi}

\noindent Note that a special module is a simple $l$-highest weight module. But in general all simple $l$-highest
weight module are not special.

\noindent For example the following result was proved in \cite{Fre} for the $sl_2$-case, \cite{Fre2} for fundamental representation in the untwisted cases, \cite{Nab, Nad} in the simply-laced untwisted case and in \cite{her06} in the general untwisted case :

\begin{thm} The Kirillov-Reshetikhin modules of an untwisted quantum affine algebra are special.\end{thm}

In the present paper we will prove that this statement also holds for twisted quantum affine algebras. This is a crucial point for the results of this paper.

\section{Twisted $T$-systems and main results}\label{deux}

In this section we state the main results about twisted $Q$-systems and twisted $T$-systems.

\subsection{Twisted $Q$-systems}

For $i\in I_\sigma, k\geq 1$ consider the Kirillov-Reshetikhin module restricted to $\U_q(\Glie^\sigma)$ : $Q_k^{(i)}=Res(W_{k,a}^{(i)})$ (it is independent of $a\in\CC^*$). 

\noindent For $i\in I_{\sigma}, k\geq 1$ define the $\U_q(\Glie^\sigma)$-module $R_k^{(i)}$ by : 
\\If $C_{i,\sigma(i)} = 2$ :
$$R_k^{(i)} = (\bigotimes_{\{j\in I_{\sigma}|C_{i,j} = -1,\sigma(j) = j\}}Q_k^{(j)})\otimes (\bigotimes_{\{j\in I_\sigma|C_{i,j} = -1, \sigma(j)\neq j\}}(Q_k^{(j)})^{\otimes m}).$$
If $C_{i,\sigma(i)} = 0$ :
$$R_k^{(i)} = (\bigotimes_{\{j\in I_{\sigma}|C_{i,j} = -1,\sigma(j) = j\}}Q_k^{(j)})\otimes (\bigotimes_{\{j\in I_\sigma|C_{i,j} = -1, \sigma(j)\neq j\}}Q_k^{(j)}).$$
If $C_{i,\sigma(i)} = -1$ :
$$R_k^{(i)} = Q_k^{(i)}\otimes (\bigotimes_{\{j\in I_{\sigma}|C_{i,j} = -1\}}Q_k^{(j)}).$$

Note that $\otimes$ is commutative in the category of finite dimensional representations of $\U_q(\Glie^\sigma)$ which is semi-simple.

\begin{thm}[The twisted $Q$-system]\label{conjdeux} Let $a\in\CC^*, k\geq 1, i\in I_\sigma$. We have :
$$Q_k^{(i)}\otimes Q_k^{(i)}=Q_{k+1}^{(i)}\otimes Q_{k-1}^{(i)}\oplus R_k^{(i)}.$$
\end{thm}

Moreover in type $A_{2n}^{(2)}$, the same relations hold between the restrictions $\overline{res}^\sigma(W_{k,a}^{(i)})$ to the subalgebra $\overline{\U}_q(\Glie^\sigma)$.

 This $Q$-system appeared in \cite{hkoty}. As a consequence we prove the conjecture of \cite{hkoty} that there is a solution which is the character of a representation (it is a purely combinatorial statement). Note that this is a particular case of the general Laurent phenomena described in another context in \cite{fz, fzd} : a priori the solutions could be rational fractions in the variables $z_i = e^{\Lambda_i}$, but in fact they are Laurent polynomials (with positive coefficients) as they correspond to characters of representations.

\noindent We will prove a stronger version of Theorem \ref{conjdeux} called twisted $T$-system (Theorem \ref{tsyst}). 

\subsection{Twisted $T$-system} For untwisted types, the $T$-system was introduced in \cite{kns} as a system of functional relations associated with solvable lattice models. Motivated by results of \cite{Fre}, it was conjectured in \cite{kosy} that the $q$-characters of Kirillov-Reshetikhin modules solve the $T$-system. This was proved in \cite{Nab, Nad} for simply-laced types and in \cite{her06} for untwisted non simply-laced types.

For twisted types, the twisted $T$-system were defined in \cite{ks}.

For $i\in I_{\sigma}, k\geq 1, a\in\CC^*$ define the $\U_q(\Lo\Glie^\sigma)$-module $S_{k,a}^{(i)}$ by : 
\\If $C_{i,\sigma(i)} = 2$ :
$$S_{k,a}^{(i)} = (\bigotimes_{\{j\in I_{\sigma}|C_{i,j} = -1,\sigma(j) = j\}}W_{k,aq_i}^{(j)})\otimes (\bigotimes_{\{j\in I_\sigma,a'\in\CC|C_{i,j} = -1, \sigma(j)\neq j,(a')^M = aq_i\}}W_{k,a'}^{(j)}).$$
If $C_{i,\sigma(i)} = 0$ :
$$S_{k,a}^{(i)} = (\bigotimes_{\{j\in I_{\sigma}|C_{i,j} = -1,\sigma(j) = j\}}W_{k,(aq)^m}^{(j)})\otimes (\bigotimes_{\{j\in I_\sigma|C_{i,j} = -1, \sigma(j)\neq j\}}W_{k,aq}^{(j)}).$$
If $C_{i,\sigma(i)} = -1$ :
$$S_{k,a}^{(i)} = W_{k,-aq}^{(i)}\otimes (\bigotimes_{\{j\in I_{\sigma}|C_{i,j} = -1\}}W_{k,aq}^{(j)}).$$

\noindent Remark : we will see later, in view of Lemma \ref{commute} and Proposition \ref{xspecial}, that in all cases the tensor products of the modules involved in the definition of $S_{k,a}^{(i)}$ commute for $\otimes$, and so $S_{k,a}^{(i)}$ is well-defined. However we only consider $\chi_q^\sigma(S_{k,a}^{(i)})$ which is clearly well-defined.

\begin{thm}[The twisted $T$-system]\label{tsyst} Let $a\in\CC^*, k\geq 1, i\in I_\sigma$. We have :
$$\chi_q^\sigma(W_{k,a}^{(i)})\chi_q^\sigma(W_{k,aq_i^2}^{(i)})=\chi_q^\sigma(W_{k+1,a}^{(i)})\chi_q^\sigma(W_{k-1,aq_i^2}^{(i)})
+\chi_q^\sigma(S_{k,a}^{(i)}).$$ 
\end{thm}

\noindent By Proposition \ref{dag}, Theorem \ref{tsyst} implies Theorem \ref{conjdeux} because $\text{res}^\sigma(W_{k,a}^{(i)})=Q_k^{(i)}$, $\text{res}^\sigma(S_{k,a}^{(i)})=R_k^{(i)}$ and the category of finite dimensional representations of $\U_q(\Glie^\sigma)$ is semi-simple. By Proposition \ref{dagdeux}, this is the same for the restrictions to $\overline{\U}_q(\Glie^\sigma)$. 

\noindent Note that as the category of finite dimensional representations of $\U_q(\Lo\Glie^\sigma)$ is not semi-simple, the twisted $T$-system can not a priori be directly translated in the category, as we did for the twisted $Q$-system.

\noindent Theorem \ref{tsyst} will follow from the following result proved in section \ref{special} :

\begin{thm}\label{domkr} The Kirillov-Reshetikhin modules of a twisted quantum affine algebra are special.\end{thm}

\subsection{Formulas for the twisted $T$-systems}\label{exptsyst} In this section we give explicit formulas for the twisted $T$-system of the theorem \ref{tsyst} (these formulas appeared in \cite{ks}). 

We denote by $X_{k,a}^{(i)}$ the representative of $W_{k,a}^{(i)}$ in the Grothendieck ring. 

\noindent Type $A_{2}^{(2)}$ :
$$X_{k,a}X_{k,aq^2} = X_{k+1,a} X_{k-1,aq^2} + X_{k,-aq}.$$

\noindent Type $A_{2n}^{(2)}$ ($n\geq 2$) : for $1\leq i\leq n - 1$ :
\begin{equation*}
\begin{split}
X_{k,a}^{(i)}X_{k,aq^2}^{(i)} &= X_{k+1,a}^{(i)} X_{k-1,aq^2}^{(i)} + X_{k,aq}^{(i - 1)}X_{k,aq}^{(i + 1)},
\\X_{k,a}^{(n-1)}X_{k,aq^2}^{(n-1)} &= X_{k+1,a}^{(n-1)} X_{k-1,aq^2}^{(n-1)} + X_{k,aq}^{(n-2)},
\\X_{k,a}^{(0)}X_{k,aq^2}^{(0)} &= X_{k+1,a}^{(0)} X_{k-1,aq^2}^{(0)} + X_{k,aq}^{(1)}X_{k,-aq}^{(0)}.
\end{split}
\end{equation*}

\noindent Type $A_{2n-1}^{(2)}$ ($n\geq 3$) : for $2\leq i\leq n - 2$ :
\begin{equation*}
\begin{split}
X_{k,a}^{(i)}X_{k,aq^2}^{(i)} &= X_{k+1,a}^{(i)} X_{k-1,aq^2}^{(i)} + X_{k,aq}^{(i - 1)}X_{k,aq}^{(i + 1)},
\\X_{k,a}^{(1)}X_{k,aq^2}^{(1)} &= X_{k+1,a}^{(1)} X_{k-1,aq^2}^{(1)} + X_{k,aq}^{(2)},
\\X_{k,a}^{(n-1)}X_{k,aq^2}^{(n-1)} &= X_{k+1,a}^{(n-1)} X_{k-1,aq^2}^{(n-1)} + X_{k,aq}^{(n - 2)}X_{k,a^2q^2}^{(n)},
\\X_{k,a}^{(n)}X_{k,aq^4}^{(n)} &= X_{k+1,a}^{(n)} X_{k-1,aq^4}^{(n)} + X_{k,a'q}^{(n-1)}X_{k,-a'q}^{(n-1)}.
\end{split}
\end{equation*} 
where $a'$ satisfies $(a')^2 = a$.

\noindent Type $D_{n+1}^{(2)}$ ($n\geq 2$) : for $2\leq i\leq n - 2$ :
\begin{equation*}
\begin{split}
X_{k,a}^{(i)}X_{k,aq^4}^{(i)} &= X_{k+1,a}^{(i)} X_{k-1,aq^4}^{(i)} + X_{k,aq^2}^{(i - 1)}X_{k,aq^2}^{(i + 1)},
\\X_{k,a}^{(1)}X_{k,aq^4}^{(1)} &= X_{k+1,a}^{(1)} X_{k-1,aq^4}^{(1)} + X_{k,aq^2}^{(2)},
\\X_{k,a}^{(n-1)}X_{k,aq^4}^{(n-1)} &= X_{k+1,a}^{(n-1)} X_{k-1,aq^4}^{(n-1)} + X_{k,aq^2}^{(n - 2)}X_{k,a'q}^{(n)}X_{k,a'q}^{(n)},
\\X_{k,a}^{(n)}X_{k,aq^2}^{(n)} &= X_{k+1,a}^{(n)} X_{k-1,aq^2}^{(n)} + X_{k,a^2q^2}^{(n - 1)},
\end{split}
\end{equation*} 
where $a'$ satisfies $(a')^2 = a$.

\noindent Type $E_6^{(2)}$ :
\begin{equation*}
\begin{split}
X_{k,a}^{(1)}X_{k,aq^2}^{(1)} &= X_{k+1,a}^{(1)} X_{k-1,aq^2}^{(1)} + X_{k,aq}^{(2)},
\\X_{k,a}^{(2)}X_{k,aq^2}^{(2)} &= X_{k+1,a}^{(2)} X_{k-1,aq^2}^{(2)} + X_{k,aq}^{(1)}X_{k,a^2q^2}^{(3)},
\\X_{k,a}^{(3)}X_{k,aq^4}^{(3)} &= X_{k+1,a}^{(3)} X_{k-1,aq^4}^{(3)} + X_{k,a'q}^{(2)}X_{k,-a'q}^{(2)}X_{k,aq^2}^{(4)},
\\X_{k,a}^{(4)}X_{k,aq^4}^{(4)} &= X_{k+1,a}^{(4)} X_{k-1,aq^4}^{(4)} + X_{k,aq^2}^{(3)},
\end{split}
\end{equation*}
where $a'$ satisfies $(a')^2 = a$.

\noindent Type $D_4^{(3)}$ :
\begin{equation*}
\begin{split}
X_{k,a}^{(1)}X_{k,aq^2}^{(1)} &= X_{k+1,a}^{(1)} X_{k-1,aq^2}^{(1)} + X_{k,a^3q^3}^{(2)},
\\X_{k,a}^{(2)}X_{k,aq^6}^{(2)} &= X_{k+1,a}^{(2)} X_{k-1,aq^6}^{(2)} + X_{k,a''q}^{(1)}X_{k,a''jq}^{(1)}X_{k,a''j^2q}^{(1)},
\end{split}
\end{equation*}
where $a''$ satisfies $(a'')^3 = a$ and $j = \text{exp}(2i\pi/3)$.
  
The specializations of these $T$-systems give the $Q$-systems of Theorem \ref{conjdeux}.

\section{Case of $\U_q^\tau$ (type $A_2^{(2)}$)}\label{adde}

First we have to study the case of $\U_q^\tau$ which is crucial for the proof in the twisted cases. Indeed twisted quantum affine algebras may have "elementary" subalgebras $\hat{\U}_i$ not only of type $A_1^{(1)}$, but also of type $A_2^{(2)}$ (see Section \ref{subalg}). So we have to treat partly the type $A_2^{(2)}$ "by hand".

\subsection{The representation $W_{2,a}$ of $\U_q^\tau$} 

The crucial result for our purposes is to describe the Kirillov-Reshetikhin modules $W_{k,a}$ of $\U_q^\tau$ for $k = 2$ as they will provide important informations for the following. 

\begin{prop}\label{kdeux} The simple representation $V = L(Z_aZ_{aq^2})$ of $\U_q^\tau$ has dimension $6$ and twisted $q$-character 
\begin{equation*}
\begin{split}
\chi_q^\tau (V) = &Z_aZ_{aq^2} + Z_aZ_{-aq^3} Z_{aq^4}^{-1} + Z_a Z_{-aq^5}^{-1} + Z_{-aq}Z_{aq^2}^{-1}Z_{-aq^3}Z_{aq^4}^{-1} 
\\&+ Z_{-aq}Z_{aq^2}^{-1}Z_{-aq^5}^{-1} + Z_{-aq^3}^{-1}Z_{-aq^5}^{-1}.
\end{split}
\end{equation*}
\end{prop}

\demo The above formula is first conjecturally given in the following way : $V$ is a subquotient of $W_{1,a}\otimes W_{1,aq^2}$, and $\chi_q^\tau(W_{1,a})\chi_q^\tau(W_{1,aq^2})$ has only two dominant monomials $Z_{a}Z_{aq^2}$ and $Z_{-aq}$. So $\chi_q^\tau(V)$ may have one dominant monomial $Z_{a}Z_{aq^2}$ and $V$ is of dimension $6$, or two dominant monomials $Z_{a}Z_{aq^2}$, $Z_{-aq}$ and $V$ is of dimension $9$. A priori it is not clear if the dimension of $V$ is equal to $9$ or $6$. If it is $6$ then $\chi_q^\tau(V)$ is equal to $\chi_q^\tau(W_{1,a})\chi_q^\tau(W_{1,aq^2}) - \chi_q^\tau(W_{1,-aq})$ which is equal to the above formula.

By the above discussion it suffices to construct a representation of highest monomial $Z_1Z_{q^2}$ of dimension $6$. To prove this result we prove that $V_{aq^2}\otimes V_a$ has a proper submodule isomorphic to $V_{-aq}$. Let $(v_0,v_1,v_2)$ be a basis of $V_{aq^2}$ and $(v_0',v_1',v_2')$ be a basis of $V_a$ as in the definition of these fundamental representations. Let 
$$\tilde{v}_0 = v_0\otimes v_1' - q v_1\otimes v_0',$$ 
$$\tilde{v}_1 = X^-_0.\tilde{v}_0 = (1-q)v_1\otimes v_1' + v_0\otimes v_2' - v_2\otimes v_0'$$ 
and 
$$\tilde{v}_2 = X^-_0.\tilde{v}_1 = v_1\otimes v_2' - q v_2\otimes v_1'.$$ 
Then $\CC. \tilde{v}_0\oplus \CC.\tilde{v}_1\oplus \CC.\tilde{v}_2$ is isomorphic to $V_{-aq}$ by identifying with the corresponding vector of the definition. It suffices to check that the action of the Drinfeld-Jimbo generators satisfy the correct relations. By construction it is clear for $K_0$ and $K_1$. For $X^-_0$ we only have to check that 
$$X^-_0\tilde{v}_2 = v_2\otimes qv_2' - q v_2\otimes v_2' = 0.$$
As $[X^+_0,X^-_0] = (K_0 - K_0^{-1})/(q^{\frac{1}{2}} - q^{-\frac{1}{2}})$, for $X^+_0$ we only have to check that
$$X^+_0.\tilde{v}_0 = q(v_0\otimes [2]_{q^{\frac{1}{2}}} v_0') -q [2]_{q^{\frac{1}{2}}}(v_0\otimes v_0') = 0.$$
As $[X^+_1,X^-_0] = 0$, for $X^+_0$ we only have to check
$$X^+_1.\tilde{v}_0 = aq^2q(1+q^2)[4]_{q^{\frac{1}{2}}}^{-1} v_2\otimes v_1' - q(v_1\otimes v_2') [4]_{q^{\frac{1}{2}}}^{-1}(1+q^2) a q$$ 
$$= (- aq)[4]_{q^{\frac{1}{2}}}^{-1} q (1+q^2)\tilde{v}_2.$$
As $[X^+_1,X^-_1] = (K_1 - K_1^{-1})/(q^2 - q^{-2})$, for $X^-_1$ we only have to check 
$$X^-_1\tilde{v}_0 = 0 \text{ as $X^-_1v_0 = X^-_1v_1 = X^-_1\tilde{v}_0 = X^-_1\tilde{v}_1 = 0$},$$
and 
$$X^-_1\tilde{v}_1 = [4]_{q^{\frac{1}{2}}}^{-1}q^{-1}[2]_{q^{\frac{1}{2}}}^2a^{-1}(1+q^{-2}) (v_0\otimes v_0') (1 - q^{-2} q^2)= 0.$$\qed

Remark : $V_a\otimes V_{aq^2}$ has no proper submodule isomorphic to $V_{-aq}$. Indeed $\tilde{v}_0$ would be represented by a vector of respective eigenvalues $(q,q^{-2})$ for $(K_0,K_1)$ and so would be of the form $\alpha = \lambda (v_1'\otimes v_0) + \mu(v_0'\oplus v_1)$. The condition $X^+_0 \alpha = 0$ gives $\mu = -q\lambda$. We can compute the vector corresponding to $\tilde{v}_2$ which is 
$$\beta = (X^-_0)^2\alpha = \lambda (v_1'\otimes v_2 - q v_2'\otimes v_1).$$
But then $X^+_1.\alpha$ should be equal to $(- aq)[4]_{q^{\frac{1}{2}}}^{-1} q (1+q^2) \beta$, contradiction. In particular $V_{aq_2}\otimes V_a$ is not semi-simple, and although $V_{aq^2}\otimes V_a$ and $V_a\otimes V_{aq^2}$ have the terms in their Jordan-Hölder decomposition, they are not isomorphic.

We can describe explicitly the representation. Let $V$ be a $6$-dimensional vector space with a basis $(v_0,v_1,v_2,v_3,v_4,v_5)$. The action of Drinfeld-Jimbo generators is given in the following tables :
$$
\begin{array}{l|l|l|l|l}
      & X^+_0                                    & X^-_0                      & K_0     & K_1
\\v_0 & 0                                      & [4]_{q^{\frac{1}{2}}}v_1 & q^2v_0  & q^{-4} v_0
\\v_1 & v_0                                    & v_2 +  v_3               & qv_1    & q^{-2} v_1
\\v_2 & \nu v_1                                & \nu v_4                  & v_2     & v_2
\\v_3 & \mu v_1                                & \mu v_4                  & v_3     & v_3
\\v_4 & v_3 + v_2                              & v_5                      & q^{-1}v_4& q^2 v_4
\\v_5 & [4]_{q^{\frac{1}{2}}} v_4              & 0                        & q^{-2} v_5& q^4 v_5
\end{array}
$$
$$
\begin{array}{l|l|l}
      & X^+_1                                                  & X^-_1                   
\\v_0 & [4]_q^{-1}q^{-1}((q^3+q^5)v_2 +(q^3 - q^2)v_3)       & 0                                  
\\v_1 & v_4                                                  & 0                     
\\v_2 & \nu [4]_{q^{\frac{1}{2}}}^{-1} q^3 (q^3 + q^{-1})v_5 & [4]_{q^{\frac{1}{2}}}^{-1}q^5\nu(q^{-3}+q^{-5})v_0    
\\v_3 & \mu [4]_{q^{\frac{1}{2}}}^{-1}q^3(q^3 - q^4)v_5      & \mu[4]_{q^{\frac{1}{2}}}q^5(q^{-3} - q^{-2}) v_0     
\\v_4 & 0 & [4]_{q^{\frac{1}{2}}}^{-1}q^{-2}(\mu 1 - \nu q^{-3} + \mu q^{-4} - \nu q^{-1}) v_1         
\\v_5 & 0                       &  q^{-3}((q^{-1} + q^{-2}) v_2 + ( q^{-2} - q^{-4}) v_3)                 
\end{array}
$$
where :
$$\nu = q^{-3/2} \frac{(1+q)(1 - q + q^2 - q^3 + q^4)}{1 - q + q^2}\text{ , }\mu = q^{-\frac{1}{2}}\frac{(1+q^2)(1+q)}{1-q+q^2}.$$
Note that we have $\nu$ and $\mu$ are invariant by $q\mapsto q^{-1}$.

So we can give the action of the Drinfeld generators : the action of the operators $x_r^{\pm}$, $k$ is given in the following table :
$$
\begin{array}{l|l|l|l}
&x_r^+&x_r^-&k
\\v_0&0                                                 &q^{3r}[4]_{q^{\frac{1}{2}}}v_1              &q^2v_0
\\v_1&q^{3r} v_0                                        &(-q^4)^r v_2 + q^r v_3 &qv_1
\\v_2&(-q^4)^r \nu v_1                                  &\nu q^r v_4           &v_2
\\v_3&q^r \mu      v_1                                  &\mu(-q^4)^r v_4      &v_3
\\v_4&(-q^4)^r v_3 +  q^r v_2   &(-q^2)^r v_5           &q^{-1}v_4
\\v_5&[4]_{q^{\frac{1}{2}}} (-q^2)^rv_4                 & 0                     &q^{-2} v_5
\end{array}
$$
Remark : we are in the situation of remark \ref{sitpart}.

And so we get the twisted $q$-character by computing the action of $h_r = \frac{q^r - q^{-r}}{r(q^{\frac{1}{2}} - q^{-\frac{1}{2}})} H_r$ for $r\neq 0$ :
\begin{equation*}
\begin{split}
  H_r.v_0 &= (1 + q^{2r})v_0,
\\H_r.v_1 &= (1 + (-q^{3})^r - q^{4r})v_1,
\\H_r.v_2 &= (1 - (-q^5)^r)v_2,
\\H_r.v_3 &= ((-q)^r - q^{2r} + (-q^3)^r - q^{4r})v_3,
\\H_r.v_4 &= ((-q)^r - q^{2r} - (-q^5)^r)v_4,
\\H_r.v_5 &= (- (-q^3)^r - (-q^5)^r)v_5.
\end{split}
\end{equation*}

\noindent Here we can prove directly that all relations between Drinfeld generators are satisfied for this representation. We give the complete proof of this point here as this is a new evidence that the Drinfeld relations hold and as this results could be used to prove them.

By the notation $v_i$ -> $v_j$, we mean the coefficient on $v_j$ of the vector obtained from $v_i$ by the action of the considered element of $\U_q^\tau$.

Relation (\ref{carp}) : by symmetry we only check that 
$$[H_r,x_m^-] = - (q^r + q^{-r} + (-1)^{r+1})x_{r+m}^-.$$

$v_0$ -> $v_1$ : 
\\$q^{3m}(1 + (-q^{3})^r - q^{4r}) - q^{3m} (1 + q^{2r}) 
= - q^{3(r+m)}(q^r +(-1)^{r+1} + q^{-r})$,

$v_1$ -> $v_2$ : 
\\$(-q^4)^m(1 - (-q^5)^r) - (-q^4)^m(1 + (-q^{3})^r - q^{4r}) 
= - (-q^4)^{m+r}(q^r +(-1)^{r+1} + q^{-r})$

$v_1$ -> $v_3$ : 
\\$q^m ((-q)^r - q^{2r} + (-q^3)^r - q^{4r})- q^m (1 + (-q^{3})^r - q^{4r})
= - q^{m+r}(q^r +(-1)^{r+1} + q^{-r})$

$v_2$ -> $v_4$ : 
\\$q^m ((-q)^r - q^{2r} - (-q^5)^r) - q^m (1 - (-q^5)^r) 
= - q^{m+r} (q^r +(-1)^{r+1} + q^{-r})$

$v_3$ -> $v_4$ : 
\\$(-q^4)^m ((-q)^r - q^{2r} - (-q^5)^r) - (-q^4)^m ((-q)^r - q^{2r} + (-q^3)^r - q^{4r}) 
\\= - (-q^4)^{m+r} (q^r +(-1)^{r+1} + q^{-r})$

$v_4$ -> $v_5$ : 
\\$(-q^2)^m (- (-q^3)^r - (-q^5)^r) - (-q^2)^m ((-q)^r - q^{2r} - (-q^5)^r)
\\= - (-q^2)^{m+r} (q^r +(-1)^{r+1} + q^{-r})$

Relation (\ref{int}) :

$v_0$ -> $v_2$ : 
\\$(-q^4)^{m+2}q^{3p} + (q -q^{-2})(-q^4)^{m+1}q^{3(p+1)} - q^{-1}(-q^4)^{m}q^{3(p+2)} 
\\= (-q^4)^m q^{3p} (q^8 - (q - q^{-2}) q^7 - q^5) = 0$
\\$q^{-1}(-q^4)^p q^{3(m+2)} + (q^{-2} - q)(-q^4)^{p+1} q^{3(m+1)} - (-q^4)^{p+2} q^{3m} 
\\= (-q^4)^pq^{3m}(q^5 - (q^{-2} -q) q^7 - q^8) = 0$

$v_0$ -> $v_3$ :
\\$q^{m+2}q^{3p} + (q -q^{-2})q^{m+1}q^{3(p+1)} - q^{-1}q^{m}q^{3(p+2)} 
= q^{m+3p}(q^2 + (q - q^{-2}) q^4 - q^5) = 0$
\\$q^{-1}q^p q^{3(m+2)} + (q^{-2} - q)q^{p+1} q^{3(m+1)} - q^{p+2} q^{3m} 
= q^pq^{3m}(q^5 + (q^{-2} -q) q^4 - q^2) = 0$

$v_1$ -> $v_4$ :
\\$\nu(q^{m+2}(-q^4)^p + (q -q^{-2})q^{m+1}(-q^4)^{p+1} - q^{-1}q^{m}(-q^4)^{p+2}) 
\\+ \mu((-q^4)^{m+2}q^p + (q -q^{-2})(-q^4)^{m+1}q^{p+1} - q^{-1}(-q^4)^m q^{p+2}) 
\\= \nu q^m(-q^4)^p(q^2 - q^6 + q^3 - q^7)
+ \mu q^p(-q^4)^m(q^8 - q^6 + q^3 - q)
\\= \nu q^p(-q^4)^m (q^7 - q^3 + q^6 - q^2) 
+\mu (-q^4)^p q^m (q - q^3 + q^6 - q^8)
\\=\nu (q^{-1}q^p (-q^4)^{m+2} 
+ (q^{-2} - q)q^{p+1} (-q^4)^{m+1} - q^{p+2} (-q^4)^m) 
\\+ \mu(q^{-1}(-q^4)^p q^{m+2} 
+ (q^{-2} -q)(-q^4)^{p+1} q^{m+1} - (-q^4)^{p+2} q^m)$

$v_2$ -> $v_5$
\\$(-q^2)^{m+2}q^p + (q -q^{-2})(-q^2)^{m+1}q^{p+1} - q^{-1}(-q^2)^{m}q^{p+2} 
\\= (-q^2)^mq^p(q^4 - (q - q^{-2}) q^3 - q) = 0$
\\$q^{-1}(-q^2)^p q^{m+2} + (q^{-2} - q)(-q^2)^{p+1} q^{m+1} - (-q^2)^{p+2} q^m 
\\= (-q^2)^p q^m (q - (q^{-2} -q) q^3 - q^4) = 0$

$v_3$ -> $v_5$
\\$(-q^2)^{m+2}(-q^4)^p + (q -q^{-2})(-q^2)^{m+1}(-q^4)^{p+1} - q^{-1}(-q^2)^{m}(-q^4)^{p+2} 
\\= (-q^2)^m(-q^4)^p(q^4 + (q - q^{-2}) q^6 - 7) = 0$
\\$q^{-1}(-q^2)^p (-q^4)^{m+2} + (q^{-2} - q)(-q^2)^{p+1} (-q^4)^{m+1} - (-q^2)^{p+2} (-q^4)^m 
\\= (-q^2)^p (-q^4)^m (q^7 + (q^{-2} -q) q^6 - q^4) = 0$

Relation (\ref{serreun})

$v_0$ -> $v_4$ : the left member has $6. 6 = 36$ terms. We associate the terms in $6$ sums corresponding to $q^{3k}(-q^4)^l q^m$ (and permutations of $k,l,m$). We get :
\\$q^{3k}(-q^4)^l q^m (\nu (q^{3/2} q - (q^{\frac{1}{2}} + q^{-\frac{1}{2}})(-q^4) + q^{-3/2}q^3) 
\\+ \mu (q^{3/2}(-q^4) - (q^{\frac{1}{2}} + q^{-\frac{1}{2}})q + q^{-3/2}q^3))
\\= q^{3k}(-q^4)^l q^m q^{\frac{1}{2}}(\nu (q + q^2 + q^3 + q^4) + \mu (-q^5 -1 ))
\\= q^{3k}(-q^4)^l q^m q^{\frac{1}{2}}(q+1) (\nu q(1 + q^2) - \mu (1 - q + q^2 - q^3 + q^4)) = 0.$

$v_1$ -> $v_5$ :
\\$q^{3k}(-q^4)^l q^m (\nu (q^{3/2} (-q^2) - (q^{\frac{1}{2}} + q^{-\frac{1}{2}}) q + q^{-3/2}(-q^4)) 
\\+ \mu (q^{3/2}(-q^2) - (q^{\frac{1}{2}} + q^{-\frac{1}{2}}) (-q^4) + q^{-3/2}q))
\\= q^{3k}(-q^4)^l q^m q^{\frac{1}{2}}(\nu (-q^3 -q - 1 - q^2) + \mu (q^4 + q^{-1}))
\\= q^{3k}(-q^4)^l q^m q^{\frac{1}{2}}(q+1) q^{-1}(-\nu q(1 + q^2) + \mu (1 - q + q^2 - q^3 + q^4)) = 0.$

$v_4$ -> $v_0$ :
\\$\nu (q^{3/2} q^{-3} - (q^{\frac{1}{2}} + q^{-\frac{1}{2}})(-q^{-4}) + q^{-3/2}q{-1}) 
\\+ \mu (q^{3/2}q^{-3} - (q^{\frac{1}{2}} + q^{-\frac{1}{2}})q^{-1} + q^{-3/2} (-q^{-4}))
\\= q^{-\frac{1}{2}}(\nu (q^{-1} + q^{-2} + q^{-3} + q^{-4}) +\mu (-q^{-5} - 1 ))
\\= q^{-\frac{1}{2}}(q^{-1}+1) (\nu q^{-1}(1 + q^{-2}) - \mu (1 - q^{-1} + q^{-2} - q^{-3} + q^{-4})) = 0.$

$v_5$ -> $v_1$ :
\\$\nu (q^{3/2} (-q^4)^{-1} - (q^{\frac{1}{2}} + q^{-\frac{1}{2}}) q^{-1} + q^{-3/2}(-q^{-2})^{-1}) 
\\+ \mu (q^{3/2}q^{-1} - (q^{\frac{1}{2}} + q^{-\frac{1}{2}}) (-q^{-4}) + q^{-3/2}(-q^2)^{-1})
\\= q^{-\frac{1}{2}}(\nu (-q^{-3} - q^{-1} - 1 - q^{-2}) + \mu (q^{-4} + q))
\\= q^{-\frac{1}{2}}(q^{-1}+1) q(-\nu q^{-1}(1 + q^{-2}) + \mu (1 - q^{-1} + q^{-2} - q^{-3} + q^{-4})) = 0.$

Relation (\ref{serredeux})

$v_0$ -> $v_4$ :
\\$q^{3k}(-q^4)^l q^m (\nu (q^{-3/2} q^{-1} - (q^{\frac{1}{2}} + q^{-\frac{1}{2}})(-q^4)^{-1} + q^{3/2}q^{-3}) 
\\+ \mu (q^{-3/2}(-q^4)^{-1} - (q^{\frac{1}{2}} + q^{-\frac{1}{2}})q^{-1} + q^{3/2}q^{-3}))
\\= q^{3k}(-q^4)^l q^m q^{-\frac{1}{2}}(\nu (q^{-1} + q^{-2} + q^{-3} + q^{-4}) + \mu (- q^{-5} - 1))
\\= q^{3k}(-q^4)^l q^m q^{-1\frac{1}{2}}(q+1) (\nu q(1 + q^2) - \mu (1 - q + q^2 - q^3 + q^4)) = 0.$

$v_1$ -> $v_5$ : 
\\$q^{3k}(-q^4)^l q^m (\nu (q^{-3/2} (-q^2)^{-1} - (q^{\frac{1}{2}} + q^{-\frac{1}{2}}) q^{-1} + q^{3/2}(-q^4)^{-1}) 
\\+ \mu (q^{-3/2}(-q^2)^{-1} - (q^{\frac{1}{2}} + q^{-\frac{1}{2}}) (-q^4)^{-1} + q^{3/2}q^{-1}))
\\= q^{3k}(-q^4)^l q^m q^{-\frac{1}{2}}(\nu (-q^{-3} - 1 - q^{-1} - q^{-2}) + \mu (q^{-4} + q ))
\\= q^{3k}(-q^4)^l q^m q^{-9/2}(q+1) (-\nu q(1 + q^2) + \mu (1 - q + q^2 - q^3 + q^4)) = 0.$

$v_4$ -> $v_0$ :
\\$\nu (q^{-3/2} q^3 - (q^{\frac{1}{2}} + q^{-\frac{1}{2}})(-q^4) + q^{3/2}q) + \mu (q^{-3/2}q^3 - (q^{\frac{1}{2}} + q^{-\frac{1}{2}})q + q^{3/2}(-q^4))
\\= q^{\frac{1}{2}}(\nu (q + q^2 + q^3 + q^4) + \mu (- q^5 - 1))
\\= q^{1\frac{1}{2}}(q^{-1}+1) (\nu q^{-1}(1 + q^{-2}) - \mu (1 - q^{-1} + q^{-2} - q^{-3} + q^{-4})) = 0.$

$v_5$ -> $v_1$ : 
\\$\nu (q^{-3/2} (-q^4) - (q^{\frac{1}{2}} + q^{-\frac{1}{2}}) q + q^{3/2}(-q^2)) + \mu (q^{-3/2}q - (q^{\frac{1}{2}} + q^{-\frac{1}{2}}) (-q^4) + q^{3/2}(-q^2))
\\=  q^{\frac{1}{2}}(\nu (-q^2 - 1 - q - q^3) + \mu (q^4 + q^{-1}))
\\=  q^{9/2}(q^{-1}+1) (-\nu q^{-1}(1 + q^{-2}) + \mu (1 - q^{-1} + q^{-2} - q^{-3} + q^{-4})) = 0.$

Relation (\ref{plusmoins}) : it suffices to prove for $r'\in\ZZ$ the two relations : 
$$\sum_{r\geq 0} [x_{r + r'}^+,x_{-r'}^-] z^r = \frac{\phi^+(z) - k^{-1}}{q^{\frac{1}{2}} - q^{-\frac{1}{2}}}\text{ , }\sum_{r\leq 0} [x_{r + r'}^+,x_{-r'}^-] z^r = \frac{k - \phi^-(z)}{q^{\frac{1}{2}} - q^{-\frac{1}{2}}}.$$

$v_0$ : 
$\sum_{r\geq 0} [4]_{q^{\frac{1}{2}}}q^{3(r+r')}q^{-3r'} z^r 
= \frac{q^2 - q^{-2}}{(q^{\frac{1}{2}} - q^{-\frac{1}{2}})(1-q^3 z)} 
= \frac{q^2 - qz - q^{-2}  + qz}{(q^{\frac{1}{2}} - q^{-\frac{1}{2}})(1 - q^3z)}
\\=\frac{1}{q^{\frac{1}{2}} - q^{-\frac{1}{2}}}(q^2\frac{(1-q^{-1}z)(1-qz)}{(1-qz)(1-q^3z)} - q^{-2}).$ 
\\$\sum_{r\leq 0}  [4]_{q^{\frac{1}{2}}}q^{3(r+r')}q^{-3r'} z^r 
= \frac{q^2 - q^{-2}}{(q^{\frac{1}{2}} - q^{-\frac{1}{2}})(1-q^{-3} z^{-1})} 
= \frac{-q^2 + qz + q^2  - q^5z}{(q^{\frac{1}{2}} - q^{-\frac{1}{2}})(1 - q^3z)}
\\= \frac{1}{q^{\frac{1}{2}}-q^{-\frac{1}{2}}}(-q^2\frac{(1-q^{-1}z)(1-qz)}{(1-qz)(1-q^3z)} + q^2).$ 

$v_1$ : 
$\sum_{r\geq 0} ( (-q^4)^{r+r'}\nu (-q^4)^{-r'}  + q^{r+r'}\mu q^{-r'} - q^{-3r'}[4]_{q^{\frac{1}{2}}}q^{3(r+r')})z^r 
\\= \frac{\nu}{1 + q^4z}+ \frac{\mu}{1 - qz} + \frac{q^{-2} - q^2}{(q^{\frac{1}{2}} - q^{-\frac{1}{2}})(1-q^3 z)} 
\\= \frac{\nu(q^{\frac{1}{2}} - q^{-\frac{1}{2}})(1-q^3z)(1 - qz) + \mu(q^{\frac{1}{2}} - q^{-\frac{1}{2}})(1-q^3z)(1+q^4z)+(q^{-2} - q^2)(1-qz)(1+q^4z))}{(q^{\frac{1}{2}} - q^{-\frac{1}{2}})(1-q^3z)(1-qz)(1+q^4z)} 
\\= \frac{1}{q^{\frac{1}{2}}-q^{-\frac{1}{2}}} \frac{z^2(-q^8 + q^6 + q^4 -q^3 -q^2) + z(q^2 - q^6) + q - q^{-1}}{(1-qz)(1-q^3z)(1+q^4z)}
\\= \frac{1}{q^{\frac{1}{2}}-q^{-\frac{1}{2}}}(q \frac{(1-q^{-1}z)(1+q^2z)(1-q^5z)}{(1-qz)(1-q^3z)(1+q^4z)} - q^{-1}).$ 
\\$\sum_{r\leq 0} ( (-q^4)^{r+r'}\nu (-q^4)^{-r'}  + q^{r+r'}\mu q^{-r'} - q^{-3r'}[4]_{q^{\frac{1}{2}}}q^{3(r+r')})z^r 
\\= \frac{\nu}{1 + q^{-4}z^{-1}}+ \frac{\mu}{1 - q^{-1}z^{-1}} + \frac{q^{-2} - q^2}{(q^{\frac{1}{2}} - q^{-\frac{1}{2}})(1-q^{-3} z^{-1})} 
\\= \frac{1}{q^{\frac{1}{2}}-q^{-\frac{1}{2}}} \frac{z^2(-q^8 + q^6 + q^4 -q^3 -q^2) + z(q^2 - q^6) + q - q^{-1}}{(1-qz)(1-q^3z)(1+q^4z)}
= \frac{1}{q^{\frac{1}{2}}-q^{-\frac{1}{2}}}(q - q \frac{(1-q^{-1}z)(1+q^2z)(1-q^5z)}{(1-qz)(1-q^3z)(1+q^4z)}).$ 

$v_2$ : $v_3$ does not appear as $(-q^4)^{r+r'}\nu q^{-r'}  - q^{-r'}\nu(-q^4){r+r'} = 0$.
\\$\sum_{r\geq 0} ( q^{r+r'}\nu q^{-r'}  - (-q^4)^{r+r'}\nu (-q^4)^{-r'} )z^r 
= \nu (\frac{1}{1 - qz} - \frac{1}{1+q^4z})
\\= q^{-3/2} \frac{(1+q^5)q(1+q^3)}{1 - q + q^2}\frac{z}{(1-qz)(1+q^4z)}
= q^{-\frac{1}{2}} \frac{(1+q^5)(1+q)(1-q^{-1})}{(1-q^{-1})}\frac{z}{(1-qz)(1+q^4z)}
\\= \frac{1}{q^{\frac{1}{2}}-q^{-\frac{1}{2}}}(\frac{z(q^6 - q^{-1} + q - q^4)}{(1-qz)(1+q^4z)})
= \frac{1}{q^{\frac{1}{2}}-q^{-\frac{1}{2}}}(\frac{(1-q^{-1}z)(1+q^6z)}{(1-qz)(1+q^4z)} - 1).$ 
\\$\sum_{r\leq 0} ( q^{r+r'}\nu q^{-r'}  - (-q^4)^{r+r'}\nu (-q^4)^{-r'} )z^r 
= \nu (\frac{1}{1 - q^{-1}z^{-1}} - \frac{1}{1+q^{-4}z^{-1}})
\\= -\nu (\frac{(q^4 + q) z }{(1 - qz)(1+q^4z)})
= \frac{1}{q^{\frac{1}{2}}-q^{-\frac{1}{2}}}(-\frac{(1-q^{-1}z)(1+q^6z)}{(1-qz)(1+q^4z)} + 1).$ 

$v_3$ : $v_2$ does not appear as $q^{r+r'}\mu (-q^4)^{-r'}  - (-q^4)^{-r'}\mu q{r+r'} = 0$.
\\$\sum_{r\geq 0} ( (-q^4)^{r+r'}\mu (-q^4)^{-r'}  - q^{r+r'}\mu q^{-r'} )z^r 
\\= \mu (\frac{1}{1 + q^4 z} - \frac{1}{1 - qz})
= q^{-\frac{1}{2}}(-q(1+q^2)(1+q)^2)\frac{z}{(1-qz)(1+q^4z)}
\\= q^{-\frac{1}{2}} \frac{(1+q^2)(1 -q^2)(1+q))}{(1-q^{-1})}\frac{z}{(1-qz)(1+q^4z)}
\\= \frac{1}{q^{\frac{1}{2}}-q^{-\frac{1}{2}}}(\frac{z(1 - q^5 + q + q^4)}{(1-qz)(1+q^4z)})
= \frac{1}{q^{\frac{1}{2}}-q^{-\frac{1}{2}}}(\frac{(1 + z)(1-q^5z)}{(1+q^4z)(1-qz)} - 1)
\\= \frac{1}{q^{\frac{1}{2}}-q^{-\frac{1}{2}}}(\frac{(1 + z)(1+q^2z)(1-q^3z)(1-q^5z)}{(1+q^2z)(1+q^4z)(1-qz)(1-q^3z)} - 1).$ 
\\$\sum_{r\leq 0} ( (-q^4)^{r+r'}\mu (-q^4)^{-r'}  - q^{r+r'}\mu q^{-r'} )z^r 
= \mu (\frac{1}{1 + q^{-4}z^{-1}} - \frac{1}{1 - q^{-1}z^{-1}})
= -\mu (\frac{z(-q-q^4)}{(1 + q^4 z)(1 - qz)})
\\= \frac{1}{q^{\frac{1}{2}}-q^{-\frac{1}{2}}}(-\frac{(1 + z)(1+q^2z)(1-q^3z)(1-q^5z)}{(1+q^2z)(1+q^4z)(1-qz)(1-q^3z)} + 1).$ 
\qed

\subsection{The special property and the twisted $T$-system}

\begin{lem}\label{specsi} For $a\in\CC^*$, $k\geq 1$, the simple representation $W_{k,a}$ of $\U_q^\tau$ is special.\end{lem}
 
\demo Let us prove by induction on $k\geq 1$ that a monomial $m\in\mathcal{M}(W_{k,a})-\{m_{k,a}\}$ satisfies $m\leq m_{k,a}A_{aq^{2k}}^{-1}$. For $k = 1$ this follows from section \ref{exaide} and for $k = 2$ from Proposition \ref{kdeux}. Let $k\geq 3$ and $m\in\mathcal{M}(W_{k,a})$. Then $m$ is a monomial of $\chi_q^\tau(W_{1,a})\chi_q^\tau(W_{k-1,aq^2})$. From the induction hypothesis we can suppose that $m = m' m_{k-1,aq^2}$ where $m'\in\mathcal{M}(W_{1,a})$. In particular $v_{aq^5}(m(m_{k,a})^{-1}) = 0$. But $m$ occurs in $\chi_q^\tau(W_{2,a})\chi_q^\tau(W_{k-2,aq^4})$. So $m = m'' m_{k-2,aq^4}$ where $m''\in\mathcal{M}(W_{2,a})$. As $v_{aq^5}(m(m_{k,a})^{-1}) = 0$ it follows from the case $k = 2$ that $m'' = m_{2,a}$. So $m = m_{k,a}$.
\qed

In particular :

\begin{prop}\label{antau} For $k\geq 1,a\in\CC^*$ we have the explicit formula :
$$\chi_q^\tau(W_{k,a}) = m_{k,a} \sum_{0\leq R'\leq R\leq k}(\prod_{r = 1\cdots R}A_{aq^{2k + 1 - 2r}}^{-1})(\prod_{r' = 1\cdots R'}A_{-aq^{2k + 2 - 2r'}}^{-1}).$$
They satisfy the twisted $T$-system :
$$\chi_q^\tau(W_{k,a})\chi_q^\tau(W_{k,aq^2}) 
= \chi_q^\tau(W_{k+1,a})\chi_q^\tau(W_{k-1,aq^2}) + \chi_q^\tau(W_{k,-aq}).$$
\end{prop}

These formulas were considered in \cite{ks} as combinatorial solutions of the twisted $T$-systems. Here we prove that they correspond to characters of Kirillov-Reshetikhin modules.

\mk

\demo It is clear that the twisted $T$-system is combinatorially satisfied by the explicit formula (in \cite{ks} it is proved in a more general situation). It also clear that the formula contains a unique dominant monomials. 

From Section \ref{exfund} and Proposition \ref{kdeux}, $\chi_q^\tau(W_{k,a})$ is equal to the formula for $k\geq 2$.

Let $k\geq 2$. In the proof of Lemma \ref{specsi}, we get that the monomials of $\chi_q^\tau(W_{k,aq^2})-m_{k,aq^2}$ are lower than $m_{k,aq^2}A_{aq^{2k+2}}^{-1}$ in the sense of monomials. So the dominant monomials appearing in $\chi_q^\tau(W_{1,a})\chi_q^\tau(W_{k,aq^2})$ are $m_{k+1,a}$ and $Z_{-aq}m_{k-1,aq^4}$. The dominant monomials appearing in $\chi_q^\tau(W_{1,-aq})\chi_q^\tau(W_{k-1,aq^4})$ are $Z_{-aq}m_{k-1,aq^4}$ and $m_{k-2,aq^6}$. So :
$$\chi_q^\tau(W_{k+1,a}) = \chi_q^\tau(W_{1,a})\chi_q^\tau(W_{k,aq^2}) - \chi_q^\tau(W_{1,-aq})\chi_q^\tau(W_{k-1,aq^4}) + \chi_q^{\tau}(W_{k-2,aq^6}).$$
It is clear that this relation is also satisfied by the explicit relation. So by induction on $k$, the formula is equal to $\chi_q^\tau(W_{k,a})$.
\qed

Formulas for general types $A_n^{(2)}$ will be proved section \ref{expkr}.

\subsection{Weyl modules} For $m$ a dominant monomials, the Weyl module $W(m)$ corresponding to $m$ is by definition the maximal finite dimensional representations of $l$-highest weight $m$.

Kashiwara proved a cyclicity property of some tensor products of fundamental representations :

\begin{thm}\label{cyc}\cite{kas} 
For $m$ a dominant monomial, there is a tensor product of fundamental representations which is of $l$-highest weight $m$.
\end{thm}

As a consequence, $\text{dim}(W(m))$ is larger that the product $p(m)$ of the dimensions of the corresponding fundamental representations.

The following result is a direct consequence of \cite{bn} where important general results are obtained (see also \cite{a} for type $A_2^{(2)}$). Although we only need this result for $A_2^{(2)}$, we state the result in general :

\begin{thm}\label{weylt} For twisted quantum affine algebras, Weyl modules are tensor product of fundamental representations.\end{thm}

By Theorem \ref{cyc}, it suffices to prove that $\text{dim}(W(m))\leq p(m)$. All arguments are contained in \cite[Section 4]{bn}. As this is explained in \cite[Section 6.4]{cm} for untwisted case, we only sketch the proof :

Let $V(\lambda)$ be the extremal weight module of extremal weight $\lambda$ and extremal vector $u_\lambda$, and let $\tilde{V}(\lambda)$ be the corresponding tensor product of level $0$ extremal fundamental representations with $\tilde{u}(\lambda)$ the tensor product of the extremal vectors. There is an injective morphism of $\U_q(\Lo\Glie^\sigma)$-modules from $V(\lambda)$ to $\tilde{V}(\lambda)$ such that $u_\lambda\mapsto \tilde{u}_\lambda$. The action of certain elements of $\U_q(\Lo\Hlie^\sigma)$ in $V(\lambda)$ corresponds to symmetric functions of the Kashiwara automorphisms of $\tilde{V}(\lambda)$. As the quotient of $\tilde{V}(\lambda)$ for the Kashiwara automorphisms is a tensor product of finite dimensional fundamental representations of dimension $p(m)$, it is possible to conclude. All details are contained in \cite{bn}.

In particular :

\begin{cor}\label{Weyl} Let $P(u)\in\CC[u]$ such that $P(0) = 1$ and $m$ be the corresponding monomial. 
The Weyl module $W(m)$ corresponding to $m$ is of dimension $3^{\text{deg}(P)}$ and the monomial appearing in his $q$-character are lower than $m$ in the sense of monomials.
\end{cor}

\demo 
As this result is known for fundamental representations from section \ref{exfund}, the corollary is a consequence of Theorem \ref{weylt} and Theorem \ref{mult}.
\qed

\section{The Kirillov-Reshetikhin modules are special}\label{special} 

In this section we prove the special property of the Kirillov-Reshetikhin modules of twisted quantum affine algebras. This is the crucial point for the proof of the twisted $T$-systems and other results in this paper.

\subsection{Preliminary results} 

We can prove as \cite[Lemma 3.3]{her05} :

\begin{lem}\label{soush} Let $V$ be a finite dimensional $\U_q(\Lo\Glie^\sigma)$-module. For $W\subset V$ a $\U_q(\Lo\Hlie^\sigma)$-submodule of $V$ and $i\in I_\sigma$, $W_i'=\underset{r\in\ZZ}{\sum}x_{i,r}^-.W$ is a $\U_q(\Lo\Hlie^\sigma)$-submodule of $V$.\end{lem}

\begin{lem}\label{submod}\cite[Lemma 3.4]{her05} Let $V$ be a finite dimensional $\U_q(\Lo sl_2)$-module. For $p\in\ZZ$, let 
$$L_{\geq p}=\underset{q \geq p}{\sum}L_{q\Lambda}\text{ and }L'_{\geq p}=\underset{r\in\ZZ}{\sum}x_r^-.L_{\geq p}.$$ Then $L_{\geq p}, L_{\geq p}'$ are $\U_q(\Lo\Hlie)$-submodule of $L$ and we have
$$((L_{\geq p}')_{m}\neq 0)\Rightarrow (\exists m'\text{ , }(L_{\geq p})_{m'}\neq \{0\}\text{ and }m \leq m').$$\end{lem}

As a consequence of Corollary \ref{Weyl}, we have an analog result for $\U_q^\tau$ which can be proved exactly as \cite[Lemma 3.4]{her05} :

\begin{lem}\label{submodtau} Let $V$ be a finite dimensional $\U_q^\tau$-module. For $p\in\ZZ$ let 
$$L_{\geq p}=\underset{M|\beta(M) = z^q,q \geq p}{\sum}L_M\text{ and }L'_{\geq p}=\underset{r\in\ZZ}{\sum}x_r^-.L_{\geq p}.$$ Then $L_{\geq p}, L_{\geq p}'$ are $\U_q(\Lo\Hlie^\tau)$-submodule of $L$ and we have
$$((L_{\geq p}')_{m}\neq 0)\Rightarrow (\exists m'\text{ , }(L_{\geq p})_{m'}\neq \{0\}\text{ and }m \leq m').$$\end{lem}

\subsection{Proof of Theorem \ref{domkr}}\label{pf} 

First we can prove as \cite[Lemma 5.3]{hers} that :

\begin{prop}\label{debut} Let $m$ be a dominant monomial and $M_i$ be the sub $\hat{\U}_i$-module of $L(m)$ generated by an $l$-highest weight vector of $L(m)$. Then $M_i$ is simple. 
\end{prop}

The theorem \ref{domkr} is a direct consequence of Lemma \ref{rn} and the following result :

\begin{lem}\label{domkrlem} Let $m\in\mathcal{M}(W_{k,a}^{(i)})-\{m_{k,a}^{(i)}\}$. We have
\begin{equation*}
m\leq \begin{cases}
m_{k,a}^{(i)}A_{i,aq_i^{2k-1}}^{-1}\text{ if $\hat{\Glie}^\sigma$ is not of type $A_{2n}^{(2n)}$,}
\\m_{k,a}^{(i)}A_{i,aq^{2k-1}}^{-1}\text{ if $\hat{\Glie}^\sigma$ is of type $A_{2n}^{(2)}$.}
\end{cases}
\end{equation*}
 In particular $m$ is right-negative and not dominant.\end{lem}

The proof of the analog result for untwisted cases \cite[Lemma 4.4]{her06} is modified :

\mk

\demo For $m\leq m_{k,a}^{(i)}$ we denote $w(m)=v(m(m_{k,a}^{(i)})^{-1})$ (the extended definition of $v$ for general monomials was given in Section \ref{addef}). 

From Proposition \ref{debut}, the $\hat{\U}_i$-submodule $M_i$ of $W_{k,a}^{(i)}$ generated by $v$ is simple. As
$${\sum}_{\{m\leq m_{k,a}^{(i)}|w(m)=1\}}M_m\subset M_{k\Lambda_i^{\vee}-\alpha_i}\subset
\hat{\U}_i.v,$$
we have ${\sum}_{\{m\leq m_{k,a}^{(i)}|w(m)=1\}}(W_{k,a}^{(i)})_m\subset M_i$. So it follows from Proposition \ref{debut}, and from the (1) of Proposition \ref{aidesldeux} (resp. from Proposition \ref{antau}) if $C_{i,\sigma(i)}\geq 0$ (resp. if $C_{i,\sigma(i)} = -1$) that :
\begin{equation*}
\underset{m\leq m_{k,a}^{(i)}/w(m)=1}{\sum}(W_{k,a}^{(i)})_m =
\begin{cases}
(W_{k,a}^{(i)})_{m_{k,a}^{(i)}A_{i,aq_i^{2k-1}}^{-1}}\text{ if $\hat{\Glie}^\sigma$ is not of type $A_{2n}^{(2)}$,}
\\(W_{k,a}^{(i)})_{m_{k,a}^{(i)}A_{i,aq^{2k-1}}^{-1}}\text{ if $\hat{\Glie}^\sigma$ is of type $A_{2n}^{(2)}$,}
\end{cases}
\end{equation*}
and that this space is of dimension $1$.

\noindent Now consider $m\in\mathcal{M}(W_{k,a}^{(i)})$ such that $m\neq m_{k,a}^{(i)}$, and let us prove by induction on $w(m)\geq 1$ that 
\begin{equation*}
m\leq 
\begin{cases}m_{k,a}^{(i)}A_{i,aq_i^{2k-1}}^{-1} \text{ if $\hat{\Glie}^\sigma$ is not of type $A_{2n}^{(2)}$,} 
\\m_{k,a}^{(i)}A_{i,aq^{2k-1}}^{-1}\text{ if $\hat{\Glie}^\sigma$ is of type $A_{2n}^{(2)}$.}
\end{cases} 
\end{equation*}
For $w(m)=1$ we have proved that $m=m_{k,a}^{(i)}A_{i,aq_i^{2k-1}}^{-1}$. In general suppose
that $w(m)=p+1$ ($p\geq 1$). It follows from the structure of $M_i$ (which is also a $\U_q(\Lo\Hlie^\sigma)$-module) that we can suppose that $(M_i)_m=\{0\}$. Consider :
$$W={\bigoplus}_{\{m'\leq m_{k,a}^{(i)}|w(m')\leq p\}}(W_{k,a}^{(i)})_{m'}.$$ 
Here $\lambda - k\Lambda_i$ is 
Note that $W$ is a $\U_q(\Lo\Hlie^\sigma)$-submodule of $W_{k,a}^{(i)}$. As $W_{k,a}^{(i)}$ is a $l$-highest weight module, we have :
$${\bigoplus}_{\{m'\leq m_{k,a}^{(i)}|w(m')=p+1\}}(W_{k,a}^{(i)})_{m'} \subset \underset{j\in I}{\sum}W_j$$
$$\text{ where }W_j=\underset{r\in\ZZ}{\sum}x_{j,r}^-.W.$$ 
For $j\in I$, $W_j$ is
a $\U_q(\Lo\Hlie^\sigma)$-submodule of $W_{k,a}^{(i)}$ (Lemma \ref{soush}). So $\exists j\in I$, $(W_j)_m\neq 
\{0\}$.

\noindent Consider the decomposition $\tau_j(\chi_q^\sigma(W_{k,a}^{(i)}))={\sum}_rP_r\otimes Q_r$ of Lemma 
\ref{aidedeux} and the corresponding decomposition of $W_{k,a}^{(i)}$ as a $\hat{\U}_j$-module:  
$W_{k,a}^{(i)}=\underset{r}{\bigoplus}V_r$. For a given $r$, consider $M_r\in\mathcal{M}(W_{k,a}^{(i)})$ such that 
$\tau_j(M_r)$ appears in $P_r\otimes Q_r$. For another such $M$, we have 
\begin{equation*}
\begin{split}
v(MM_r^{-1}) = 
\begin{cases}
u_j(\tau_j(MM_r^{-1}))/2&\text{ if $(j,\hat{\Glie}^\sigma)\neq (0,A_{2n}^{(2)})$,}
\\u_j(\tau_j(MM_r^{-1}))&\text{ if $(j,\hat{\Glie}^\sigma) = (0,A_{2n}^{(2)})$,}
\end{cases}
\end{split}
\end{equation*} 
and so : 
\begin{equation*}
\begin{split}
u_j(\tau_j(M)) = 
\begin{cases}
u_j(\tau_j(M_r))-2w(M)+2w(M_r)&\text{ if $(j,\hat{\Glie}^\sigma)\neq (0,A_{2n}^{(2)})$,}
\\ u_j(\tau_j(M_r))-w(M)+w(M_r) &\text{ if $(j,\hat{\Glie}^\sigma) = (0,A_{2n}^{(2)})$.}
\end{cases}
\end{split}
\end{equation*}
This can be rewritten :
\begin{equation*}
\begin{split}
u_j(\tau_j(M)) = 
\begin{cases}
2(p-w(M))+p_r&\text{ if $(j,\hat{\Glie}^\sigma)\neq (0,A_{2n}^{(2)})$,}
\\ (p-w(M))+p_r &\text{ if $(j,\hat{\Glie}^\sigma) = (0,A_{2n}^{(2)})$,}
\end{cases}
\end{split}
\end{equation*}
where 
\begin{equation*}
\begin{split}
p_r=
\begin{cases}
-2p+2w(M_r)+u_j(\tau_j(M_r))&\text{ if $(j,\hat{\Glie}^\sigma)\neq (0,A_{2n}^{(2)})$,}
\\-p+w(M_r)+u_j(\tau_j(M_r)) &\text{ if $(j,\hat{\Glie}^\sigma) = (0,A_{2n}^{(2)})$.}
\end{cases}
\end{split}
\end{equation*}
Note that $p_r$ does not 
depend of $M$. So we have $w(M)\leq p\Leftrightarrow u_j(\tau_j(M))\geq p_r$. So 
$$W = {\bigoplus}_r((V_r)_{\geq p_r}) = {\bigoplus}_r(V_r\cap W).$$ 
As the $V_r$ are sub 
$\hat{\U}_j$-modules of $W_{k,a}^{(i)}$, we have $W_j=\underset{r}{\bigoplus}(V_r\cap W_j)$. Let $R$ such that 
$\tau_j(m)$ is a monomial of $P_R\otimes Q_R$. We can apply Lemma \ref{submod} if $C_{j,\sigma(j)}\geq 0$ and Lemma \ref{submodtau} if $C_{j,\sigma(j)} = -1$ to the $\hat{\U}_j$-module 
$V_R$ with $p_R$ for $Q_R^{-1}\tau_j(m)$ : we get a monomial $M'$ of $(\chi_q^\sigma)^j(W)$ (twisted $q$-character as $\hat{\U}_j$-module) such that 
$$Q_R^{-1}\tau_i(m)\in M' \ZZ[\tau_j(A_{j,a}^{-1})]_{a\in\CC^*}.$$ 
Consider $m'=\tau_j^{-1}(Q_R\otimes M')$ (it is a 
monomial of $\chi_q^\sigma(W)$). 

First suppose that $m'\neq m_{k,a}^{(i)}$. So $m < m'$ and $m\in m'\ZZ[A_{j,b}^{-1}]_{b\in\CC^*}$. As $m'\neq m_{k,a}^{(i)}$ we have $w(m')\geq 1$ and the induction hypothesis gives the result.

If $m'= m_{k,a}^{(i)}$, consider 
$Z=\hat{\U}_j.(W_{k,a}^{(i)})_{m_{k,a}^{(i)}}$. Then $Z$ is equal to $(W_{k,a}^{(i)})_{m_{k,a}^{(i)}}$ or $M_i$. As $Z$ is a sub $\hat{\U}_j$-module of $V_R$ and $Q_R^{-1}\tau_j(m)$ is not a monomial of $(\chi_q^\sigma)^j(Z)$, we can use the same arguments as above with $V_R/Z$ instead of $V_R$.
\qed

\section{Proof of Theorem \ref{tsyst}}\label{proofth}

The following proof of Theorem \ref{tsyst} relies on Theorem \ref{domkr}.

\subsection{Preliminary results}\label{prel} First we get exactly as \cite[Lemma 5.1]{her06} :

\begin{lem}\label{commute} Let $V$ be a special module. Suppose that $V\simeq V_1\otimes
\cdots\otimes V_r$ where $V_1,\cdots,V_r$ are $l$-highest weight modules. Then $V_1,\cdots,V_r$ are special and for all
$\sigma$ permutation of $\{1,\cdots,r\}$, we have $V\simeq V_{\sigma(1)}\otimes \cdots\otimes V_{\sigma(r)}$.\end{lem}

Indeed it is easy to produce a dominant non highest weight monomial of $V$ from such a monomial in one the $V_i$. Besides as for any $\sigma$ the tensor product is special, it is simple.

Let $i\in I, k\geq 1, a\in\CC^*$. Let 
\begin{equation*}
\begin{split}
M = 
\begin{cases}
m_{k,a}^{(i)}m_{k,aq_i^2}^{(i)}=m_{k+1,a}^{(i)}m_{k-1,aq_i^2}^{(i)}&\text{ if $\hat{\Glie}^\sigma$ is not of type $A_{2n}^{(2)}$,} 
\\m_{k,a}^{(i)}m_{k,aq^2}^{(i)}=m_{k+1,a}^{(i)}m_{k-1,aq^2}^{(i)}&\text{ if $\hat{\Glie}^\sigma$ is of type $A_{2n}^{(2)}$,}  
\end{cases}
\end{split}
\end{equation*}
and 
\begin{equation*}
\begin{split}
M' = 
\begin{cases}
MA_{i,aq_i^{2k-1}}^{-1}...A_{i,aq_i}^{-1}&\text{ if $\hat{\Glie}^\sigma$ is not of type $A_{2n}^{(2)}$,}
\\MA_{i,aq^{2k-1}}^{-1}...A_{i,aq}^{-1}&\text{ if $\hat{\Glie}^\sigma$ is of type $A_{2n}^{(2)}$.}
\end{cases}
\end{split}
\end{equation*}
Note that $M'$ is the highest weight monomial of $S_{k,a}^{(i)}$. Let us write the dominant monomial $M'$ in a normal way : 
\begin{equation}\label{normd}M'=\underset{l=1...L}{\prod}m_{k_l,a_l}^{(i_l)}.\end{equation} 
Consider the sets of monomials for $\hat{\Glie}^\sigma$ not of type $A_{2n}^{(2)}$ :
$$\mathcal{B}=\{m_{k,a}^{(i)}A_{i,aq_i^{2k-1}}^{-1}...A_{i,aq_i^{2(k-k')-1}}^{-1}A_{i_l,a_lq_{i_l}^{2k_l-1}}^{-1}|0\leq 
k'\leq k-1\text{ , }1\leq l\leq L\},$$
$$\mathcal{B'}=\{m_{k,a}^{(i)}, m_{k,a}^{(i)}A_{i,aq_i^{2k-1}}^{-1}, 
m_{k,a}^{(i)}A_{i,aq_i^{2k-1}}^{-1}A_{i,aq_i^{2k-3}}^{-1},...,m_{k,a}^{(i)}A_{i,aq_i^{2k-1}}^{-1}...A_{i,aq_i}^{-1}\},$$
and for $\hat{\Glie}^{\sigma}$ of type $A_{2n}^{(2)}$ :
$$\mathcal{B}=\{m_{k,a}^{(i)}A_{i,aq^{2k-1}}^{-1}...A_{i,aq^{2(k-k')-1}}^{-1}A_{i_l,a_lq^{2k_l-1}}^{-1}|0\leq 
k'\leq k-1\text{ , }1\leq l\leq L\},$$
$$\mathcal{B'}=\{m_{k,a}^{(i)}, m_{k,a}^{(i)}A_{i,aq^{2k-1}}^{-1}, 
m_{k,a}^{(i)}A_{i,aq^{2k-1}}^{-1}A_{i,aq^{2k-3}}^{-1},...,m_{k,a}^{(i)}A_{i,aq^{2k-1}}^{-1}...A_{i,aq}^{-1}\}.$$

\begin{lem}\label{check} The monomials of $m_{k,aq_i^2}^{(i)}.\mathcal{B}$ are right negative.\end{lem}

\demo First we suppose that $\hat{\Glie}^\sigma$ is not of type $A_{2n}^{(2)}$.

For $b\in\CC^*$ and $m$ a monomial, let us define 
$$\mu_b(m)=\text{Max}\{l\in\ZZ|\exists i\in I_\sigma,r\in\ZZ, z_{i,(bq^l\omega^r)^{d_i}}(m)\neq 0\}.$$ 
Let $\alpha=MA_{i,aq_i^{2k-1}}^{-1}...A_{i,aq_i^{2(k-k')-1}}^{-1}$, and $\alpha A_{i_l,aq_{i_l}^{2k_l-1}}^{-1}\in m_{k,aq_i^2}^{(i)}\mathcal{B}$. It suffices to check that 
$\mu_{a'}(\alpha)<\mu_{a'}(A_{i_l,a_lq_{i_l}^{2k_l-1}}^{-1})$ where $a'$ satisfies $(a')^{d_i} = a$.

Case 1: $d_i = 1$ : $\mu_a(\alpha)\leq 2k-1$.

\noindent If $d_{i_l} = 1$ : $a_lq_{i_l}^{2k_l - 2} = aq^{2k-1}$, $\mu_a(A_{i_l,a_lq_{i_l}^{2k_l - 1}}^{-1}) = 2k + 1$.

\noindent If $d_{i_l}\geq 2$ : $a_lq_{i_l}^{2k_l - 2} = (aq^{2k-1})^M$, $\mu_a(A_{i_l,a_lq_{i_l}^{2k_l - 1}}^{-1})=2k + 1$.

Case 2 : $d_i = 2$ : $\mu_a(\alpha)\leq 2k-1$.

\noindent If $d_{i_l}=1$ : $a_lq_{i_L}^{2k_l-2} = a'q^{2k-1}$ or $-a'q^{2k-1}$, $\mu_a(A_{i_l,a_lq_{i_l}^{2k_l - 1}}^{-1}) = 2k + 1$.

\noindent If $d_{i_l}=2$ : $a_lq_{i_L}^{2k_l-2} = aq_i^{4k - 2}$, $\mu_a(A_{i_l,a_lq_{i_l}^{2k_l - 1}}^{-1}) = 2k + 1$.

Case 3 : $d_i=3$ : $\mu_a(\alpha)\leq 2k-1$.

\noindent If $d_{i_l}=1$ : $a_lq_{i_L}^{2k_l-2} = a'q^{2k-1}$ or $ja'q^{2k-1}$ or $j^2a'q^{2k-1}$, $\mu_a(A_{i_l,a_lq_{i_l}^{2k_l - 1}}^{-1}) = 2k + 1$.

Now we suppose that $\hat{\Glie}^\sigma$ is of type $A_{2n}^{(2)}$. 
For $b\in\CC^*$ and $m$ a monomial, let us define 
$$\mu_b(m)=\text{Max}\{l\in\ZZ|\exists i\in I_\sigma, z_{i,bq^l}(m)\neq 0\text{ or }z_{i,-bq^l}(m)\neq 0\}.$$ 
Let $\alpha=MA_{i,aq^{2k-1}}^{-1}...A_{i,aq^{2(k-k')-1}}^{-1}$, and $\alpha A_{i_l,aq^{2k_l-1}}^{-1}\in m_{k,aq^2}^{(i)}\mathcal{B}$. It suffices to check that $\mu_a(\alpha)<\mu_a(A_{i_l,a_lq^{2k_l-1}}^{-1})$.

Case 1 : $d_i = 1$ : $\mu_a(\alpha)\leq 2k-1$.

\noindent If $d_{i_l} = 1$ : $a_lq^{2k_l-2} = aq^{2k-1}$, $\mu_a(A_{i_l,a_lq_{i_l}^{2k_l - 1}}^{-1}) = 2k + 1$.

\noindent If $d_{i_l} = \frac{1}{2}$ : $a_lq_{i_L}^{2k_l-2} = aq^{2k-1}$, $\mu_a(A_{i_l,a_lq_{i_l}^{2k_l - 1}}^{-1}) = 2k + 1$.

Case 2 : $d_i = \frac{1}{2}$ : $\mu_a(\alpha)\leq 2k-1$.

\noindent If $d_{i_l} = 1$ : $a_lq_{i_L}^{2k_l-2} = aq^{2k-1}$, $\mu_a(A_{i_l,a_lq_{i_l}^{2k_l - 1}}^{-1}) = 2k + 1$.

\noindent If $d_{i_l} = \frac{1}{2}$ : $a_lq_{i_L}^{2k_l-2} = -aq^{2k-1}$, $\mu_a(A_{i_l,a_lq_{i_l}^{2k_l - 1}}^{-1}) = 2k + 1$.
\qed

As a consequence of Lemma \ref{check}, we get as in \cite[Proposition 5.3]{her06} :

\begin{prop}\label{xspecial} For $i\in I, k\geq 1, a\in\CC^*$, the module $S_{k,a}^{(i)}$ is special. In particular $M'$ is the unique dominant monomial of $\chi_q^\sigma(S_{k,a}^{(i)})$.\end{prop}

\subsection{Proof of the theorem \ref{tsyst}}\label{partun} The two terms of the equality of the theorem \ref{tsyst} are in $\text{Im}(\chi_q^\sigma)$ and so are
characterized by the coefficient of
their dominant monomials. So it suffices to determine the dominant
monomials of each product.

\noindent First let us prove the following lemma about the monomials of $\chi_q^\sigma(W_{k,a}^{(i)})$ : 

\begin{lem}\label{moreinfo} The monomials $\chi_q^\sigma(W_{k,a}^{(i)})$ are lower than a monomial of 
$\mathcal{B}$ or are in $\mathcal{B'}$.\end{lem}

\noindent An analog result is proved in \cite{Nad} for the untwisted simply-laced cases. The general untwisted case is proved in \cite[Lemma 5.5]{her06}. We use this proof with some modification :

\demo For $m\in\mathcal{M}(W_{k,a}^{(i)})$ we prove the statement by induction on 
$$w(m)=v(m(m_{k,a}^{(i)})^{-1})\geq 0.$$ 
For $w(m)=0$ we have $m=m_{k,a}^{(i)}\in\mathcal{B'}$. For $w(m)\geq 1$ it follows from Theorem \ref{domkr} that there is $j\in I$ such that $m\notin B_j$. So we get from Proposition \ref{jdecomp} a monomial 
$m'\in\mathcal{M}(\chi_q^\sigma(W_{k,a}^{(i)}))$ such that $w(m')<w(m)$, $m$ is a monomial of $L_j(m')$, and $L_j(m')$ appears in the decomposition of $\chi_q^\sigma(W_{k,a}^{(i)})$. In particular $m\leq m'$, and if $m'$ is lower than a monomial in $\mathcal{B}$, so is $m$. 

So we can suppose that $m'\in\mathcal{B'}$. 

If $m'=m_{k,a}^{(i)}$, we have $j=i$. If moreover $(\hat{\Glie}^\sigma,i)\neq (A_{2n}^{(2)},0)$, the monomials of $L_i(m)$ are the monomials of $\mathcal{B'}$ from Proposition \ref{debut} and Proposition \ref{aidesldeux} (1). Otherwise it follows from Proposition \ref{debut} and Proposition \ref{antau} that the monomials of $L_n(m)$ are in $\mathcal{B'}$ or are lower than a monomial in $\mathcal{B}$ (with $i_l = n$ and $a_l\in -aq^{\ZZ}$).

Suppose that $m'\neq m_{k,a}^{(i)}$. 

If $(\hat{\Glie}^\sigma,j)\neq (A_{2n}^{(2)},0)$, $L_j(m)$ corresponds to the $q$-character of a tensor product of Kirillov-Reshetikhin modules of type $sl_2$ ((2) of Proposition \ref{aidesldeux}). Let $m_{k_{l'},a_{l'}}$ be the corresponding monomials (that it to say a normal form). For each $l'$, the complex number $a_{l'}q_{j}^{2(k_{l'}-1)}$ is equal to one $a_{l}q_{i_l}^{2(k_l-1)}$ with $i_l=j$ (see the decomposition (\ref{normd}) of the section \ref{prel}). We can conclude with (1) of Proposition \ref{aidesldeux}. 

If $(\hat{\Glie}^\sigma,j) = (A_{2n}^{(2)},0)$, $L_n(m)$ corresponds to a quotient of a tensor product of Kirillov-Reshetikhin modules of type $A_2^{(2)}$. Let $m_{k_{l'},a_{l'}}$ be the corresponding monomials. For each $l'$, $a_{l'}q^{2(k_{l'}-1)}$ is equal to one $a_{l}q^{2(k_l-1)}$ with $i_l = n$ (see the decomposition (\ref{normd}) of the section \ref{prel}). We can conclude with Proposition \ref{antau}.\qed

Consider 
\begin{equation*}
\chi_1 = 
\begin{cases}
\chi_q^\sigma(W_{k,a}^{(i)})\chi_q^\sigma(W_{k,aq_i^2}^{(i)})\text{ if $\hat{\Glie}^\sigma$ is not of type $A_{2n}^{(2)}$,}
\\\chi_q^\sigma(W_{k,a}^{(i)})\chi_q^\sigma(W_{k,aq^2}^{(i)})\text{ if $\hat{\Glie}^\sigma$ is of type $A_{2n}^{(2)}$,}
\end{cases}
\end{equation*}
and
\begin{equation*}
\chi_2 =
\begin{cases}
\chi_q^\sigma(W_{k+1,a}^{(i)})\chi_q^\sigma(W_{k-1,aq_i^2}^{(i)})\text{ if $\hat{\Glie}^\sigma$ is not of type $A_{2n}^{(2)}$,}
\\\chi_q^\sigma(W_{k+1,a}^{(i)})\chi_q^\sigma(W_{k-1,aq^2}^{(i)})\text{ if $\hat{\Glie}^\sigma$ is of type $A_{2n}^{(2)}$.}
\end{cases}
\end{equation*}

As a consequence of Theorem \ref{domkr} (more precisely of Lemma \ref{domkrlem}) :

\begin{lem}\label{premier} 1) The dominant monomials of $\chi_1$ are 
\begin{equation*}
M, MA_{i,aq_i^{2k-1}}^{-1},
MA_{i,aq_i^{2k-1}}^{-1}A_{i,aq_i^{2k-3}}^{-1},...,MA_{i,aq_i^{2k-1}}^{-1}...A_{i,aq_i}^{-1}\text{ for $\hat{\Glie}^\sigma\neq A_{2n}^{(2)}$,}
\end{equation*}
\begin{equation*}
M, MA_{i,aq^{2k-1}}^{-1},
MA_{i,aq^{2k-1}}^{-1}A_{i,aq^{2k-3}}^{-1},...,MA_{i,aq^{2k-1}}^{-1}...A_{i,aq}^{-1}\text{ for $\hat{\Glie}^\sigma = A_{2n}^{(2)}$.}
\end{equation*}
2) The dominant monomials of $\chi_2$ are
\begin{equation*}
M, MA_{i,aq_i^{2k-1}}^{-1}, 
MA_{i,aq_i^{2k-1}}^{-1}A_{i,aq_i^{2k-3}}^{-1},...,MA_{i,aq_i^{2k-1}}^{-1}...A_{i,aq_i^3}^{-1}\text{ for $\hat{\Glie}^\sigma\neq A_{2n}^{(2)}$,}
\end{equation*}
\begin{equation*}
M, MA_{i,aq_i^{2k-1}}^{-1}, 
MA_{i,aq^{2k-1}}^{-1}A_{i,aq^{2k-3}}^{-1},...,MA_{i,aq^{2k-1}}^{-1}...A_{i,aq^3}^{-1}\text{ for $\hat{\Glie}^\sigma = A_{2n}^{(2)}$.}
\end{equation*}
In each case the dominant monomials appear with multiplicity $1$.\end{lem}

{\it End of the proof of the theorem \ref{tsyst} :}

\noindent The unique dominant monomial that appears in $\chi_1 - \chi_2$ is 
$M'$, and it has a multiplicity $1$. We can conclude with Theorem \ref{imchiq} because $M'$ is the unique dominant monomial of $\chi_q^\sigma(S_{k,a}^{(i)})$ (Proposition \ref{xspecial}).\qed

Remark : Theorem \ref{tsyst} also follows from Theorem \ref{domkr}, from results proved below (Theorem \ref{transi}) and from the analog result for untwisted cases (but Theorem \ref{domkr} has to be directly proved as in the present paper). However we gave the proof of (Theorem \ref{domkr} $\Rightarrow$ Theorem \ref{tsyst}) in this subsection \ref{partun} to get all intermediate results, and also for the uniformity of the proof with the untwisted cases.

\subsection{Complement : asymptotic property}\label{partdeux} It is certainly possible to get as for the untwisted case an asymptotic property of characters and twisted $q$-characters of Kirillov-Reshetikhin modules of twisted quantum affine algebras, that is to say :

(1) The normalized character $\mathcal{Q}_k^{(i)}=e^{-k\Lambda_i}\overline{\chi}^\sigma(Q_k^{(i)})$ considered as a polynomial in $e^{-\alpha_j}$ has a limit as a formal 
power series :
$$\exists \underset{k\rightarrow \infty}{\text{lim}}\mathcal{Q}_k^{(i)}\in\ZZ[[e^{-\alpha_j}]]_{j\in I_\sigma}.$$

(2) The normalized $q$-character of $W_{k,a}^{(i)}$ considered as a polynomial in $A_{j,b}^{-1}$ has a limit as a
formal power series : 
$$\exists \underset{k\rightarrow
\infty}{\text{lim}}\frac{\chi_q^\sigma(W_{k,aq_i^{-2k}}^{(i)})}{m_{k,aq_i^{-2k}}^{(i)}}\in\ZZ[[A_{j,aq^m}^{-1}]]_{j\in I_\sigma, m\in\ZZ}.$$

Note that (1) is a consequence of (2). Moreover it is certainly possible to get also the convergence in the analytic sense as in \cite{her06}.

\section{The twisted Kirillov-Reshetikhin conjecture and general results}\label{expform}

In this section by using the results of previous sections we prove a close relation between the twisted types and the untwisted types (we construct an isomorphism between the Grothendieck rings of finite dimensional representations preserving Kirillov-Reshetikhin modules). We get explicit formulas for the character of an arbitrary tensor product of Kirillov-Reshetikhin modules for all types, and we so we prove the Kirillov-Reshetikhin conjecture.
 
\subsection{Relation between twisted and untwisted types} 
 
Let us define 
 $$\pi : \ZZ[Y_{i,a}^{\pm 1}]_{i\in I, a\in\CC^*} \rightarrow \ZZ[Z_{i,a}^{\pm 1}]_{i\in I_{\sigma},a\in\CC^*},$$ 
 as the ring morphism such that for $i\in I_{\sigma}$, $p\in\ZZ$, $a\in\CC^*$ : 
\begin{equation*}
\begin{split}
\pi(Y_{\sigma^p(i),a}) = 
\begin{cases}
Z_{i,(\omega^p a)^{d_i}}&\text{ if $\hat{\Glie}^\sigma$ is not of type $A_{2n}^{(2)}$,}
\\Z_{i,a(-1)^p}&\text{ if $\hat{\Glie}^\sigma$ is of type $A_{2n}^{(2)}$.}
\end{cases}
\end{split}
\end{equation*}

To avoid confusion, the Kirillov-Reshetikhin modules of the untwisted quantum affine algebra $\U_q(\Lo\Glie)$ are denoted by $\tilde{W}_{k,a}^{(i)}$ for $a\in\CC^*$, $k\geq 0$, $i\in I$. 

We denote by $\pi(\tilde{W}_{k,a}^{(i)})$ the Kirillov-Reshetikhin module of $\U_q(\Lo\Glie^\sigma)$ of highest monomial $\pi(M)$ where $M$ is the highest monomial of $\tilde{W}_{k,a}^{(i)}$. For $i\in I_{\sigma}$, $p\in\ZZ$, $k\geq 0$ and $a\in\CC^*$, we have :
\begin{equation*}
\begin{split}
\pi(\tilde{W}_{k,a}^{(\sigma^p(i))}) = 
\begin{cases}
W_{k,(a\omega^p)^{d_i}}^{(i)} &\text{ if $\hat{\Glie}^\sigma$ is not of type $A_{2n}^{(2)}$,}
\\W_{k,a(-1)^p}^{(i)}&\text{ if $\hat{\Glie}^\sigma$ is of type $A_{2n}^{(2)}$.}
\end{cases}
\end{split}
\end{equation*}
Here there is an abuse of notation as we use $\pi$ twice, but this does not lead to confusion thanks to the following :
\begin{thm}\label{transi} $\pi$ can be uniquely extended to a well-defined ring isomorphism 
$$\pi :  \text{Rep}(\U_q(\Lo\Glie)) \rightarrow \text{Rep}(\U_q(\Lo\Glie^\sigma)),$$
and the following diagram is commutative :
$$\begin{CD}
\text{Rep}(\U_q(\Lo\Glie))  @>{\chi_q}>>  \ZZ[Y_{i,a}^{\pm 1}]_{i\in I,a\in\CC^*}\\
@VV{\pi}V     @VV{\pi}V\\
\text{Rep}(\U_q(\Lo\Glie^\sigma))  @>{\chi_q^\sigma}>>  \ZZ[Z_{i,a}^{\pm 1}]_{i\in I_\sigma,a\in\CC^*}
\end{CD}.$$
In particular for $i\in I_{\sigma}, k\geq 0, p\in\ZZ, a\in\CC^*$, we have :
\begin{equation*}
\begin{split}
 \pi(\chi_q (\tilde{W}_{k,a}^{(\sigma^p(i))})) = 
\begin{cases}
 \chi_q^\sigma(W_{k,(a\omega^p)^{d_i}}^{(i)}) &\text{ if $\hat{\Glie}^\sigma$ is not of type $A_{2n}^{(2)}$,}
\\ \chi_q^\sigma(W_{k,(-1)^p a}^{(i)})&\text{ if $\hat{\Glie}^\sigma$ is of type $A_{2n}^{(2)}$.}
\end{cases}
\end{split}
\end{equation*}
\end{thm}

Remark : although this works for Kirillov-Reshetikhin modules, the result does not mean that in general we can compute the twisted $q$-character of an arbitrary simple $\U_q(\Lo\Glie^\sigma)$-module from the $q$-character of a simple $\U_q(\Lo\Glie)$.

\mk

\demo First let us prove the last identity of the Theorem.

For $i\in I_{\sigma}$, $p\in\ZZ$, $a\in\CC^*$ we have from the explicit defining formulas
\begin{equation*}
\begin{split}
\pi(A_{\sigma^p(i),a}) = 
\begin{cases}
A_{i,(\omega^p a)^{d_i}}&\text{ if $\hat{\Glie}^\sigma$ is not of type $A_{2n}^{(2)}$,} 
\\A_{i,(-1)^p a}&\text{ if $\hat{\Glie}^\sigma$ is of type $A_{2n}^{(2)}$.}
\end{cases} 
\end{split}
\end{equation*}
From Theorem \ref{imchiqu}, we have $\pi(\text{Im}(\chi_q))\subset \mathcal{K}_i$ for any $i\in I^\sigma$ : it is obvious if $C_{i,\sigma(i)}\neq -1$. For $C_{i,\sigma(i)} = -1$ we have : 
\begin{equation*}
\begin{split}
\text{Im}(\chi_q)&\subset \mathcal{K}_i\cap\mathcal{K}_{\sigma(i)} 
\\&= \ZZ[Y_{i,b}(1+A_{i,bq}^{-1}+A_{i,bq}^{-1}A_{\sigma(i),bq^2}^{-1}),Y_{\sigma(i),b}(1+A_{\sigma(i),bq}^{-1}+A_{\sigma(i),bq}^{-1}A_{i,bq^2}^{-1})]_{b\in\CC^*}
\\&\times\ZZ[Y_{j,d}^{\pm 1}]_{j\notin\{i,\sigma(i)\},d\in\CC^*},
\end{split}
\end{equation*}
from Theorem \ref{imchiqu} for the untwisted quantum affine algebra $\U_q(\Lo sl_3)$. And so we also have the result as for $r\in\ZZ$ :
\begin{equation*}
\begin{split}
&Z_{i,aq^r}(1+A_{i,aq^{r+1}}^{-1}+A_{i,aq^{r+1}}^{-1}A_{i,-aq^{r+2}}^{-1}) 
\\= &\pi(Y_{i,aq^r}(1+A_{i,aq^{r+1}}^{-1}+A_{i,aq^{r+1}}^{-1}A_{\sigma(i),aq^{r+2}}^{-1})),
\\&Z_{i,-aq^r}(1+A_{i,-aq^{r+1}}^{-1}+A_{i,-aq^{r+1}}^{-1}A_{i,aq^{r+2}}^{-1}) 
\\= &\pi(Y_{\sigma(i),aq^r}(1+A_{\sigma(i),aq^{r+1}}^{-1}+A_{\sigma(i),aq^{r+1}}^{-1}A_{i,aq^{r+2}}^{-1})).
\end{split}
\end{equation*}

In particular $\pi(\chi_q (\tilde{W}_{k,a}^{(\sigma^p(i))}))$ is in $\text{Im}(\chi_q^\sigma)$. 
\\Moreover for a monomial $m\in\ZZ[Y_{j,b}]_{j\in I,b\in\CC^*}$ we have : 
$$\text{($m$ is dominant) }\Longleftrightarrow\text{ ($\pi(m)$ is dominant).}$$ 
In particular as the $\U_q(\Lo\Glie)$-module $\tilde{W}_{k,a}^{(\sigma^p(i))}$ is special, $\pi(\chi_q(\tilde{W}_{k,a}^{(\sigma^p(i))}))$ has a unique dominant monomial. For $\hat{\Glie}^\sigma$ not of type $A_{2n}^{(2)}$ (resp. of type $A_{2n}^{(2)}$), the $\U_q(\Lo\Glie^\sigma)$-module $W_{k,(a\omega^p)^{d_i}}^{(i)}$ (resp. $W_{k,(-1)^pa}^{(i)}$) is special, so $\chi_q^\sigma(W_{k,(\omega^p a)^{d_i}}^{(i)})$ (resp. $\chi_q^\sigma(W_{k,(-1)^p a}^{(i)})$) has a unique dominant monomial equal to $\pi(M)$ where $M$ is the unique dominant monomial of $\chi_q(\tilde{W}_{k,a}^{(\sigma^p(i))})$. This implies the result.

As $\chi_q$ and $\chi_q^\sigma$ are injective ring morphisms, this result implies that $\pi$ is well-defined. Then the diagram is clearly commutative. As the fundamental representations generate the ring $\text{Rep}(\U_q(\Lo\Glie^\sigma))$, $\pi$ is surjective. The injectivity follows from the injectivity of $\pi\circ\chi_q$.
\qed

Let us define $\pi : \ZZ[y_i^{\pm 1}]_{i\in I} \rightarrow \ZZ[z_i^{\pm 1}]_{1\leq i\leq n}$ as the ring morphism such that for $i\in I$, 
\begin{equation*}
\begin{split}
\pi(y_i) = \begin{cases}
z_{\overline{i}}\text{ if $\hat{\Glie}^\sigma \neq A_{2n}^{(2)}$,}
\\z_{n-\overline{i}}\text{ if $\hat{\Glie}^\sigma = A_{2n}^{(2)}$.}
\end{cases}
\end{split}
\end{equation*}
There is an abuse of notation as we used $\pi$ for different maps, but this does not lead to confusion as from Proposition \ref{dag}, we have :
\begin{cor}\label{squn} The following diagram is commutative :
$$\begin{CD}
\text{Rep}(\U_q(\Lo\Glie))   @>{\chi\circ\text{res}}>>  \ZZ[y_i^{\pm 1}]_{i\in I}\\
@VV{\pi}V     @VV{\pi}V\\
\text{Rep}(\U_q(\Lo\Glie^\sigma))   @>{\chi^\sigma\circ\text{res}^\sigma}>>  \ZZ[z_i^{\pm 1}]_{1\leq i\leq n}
\end{CD}.$$
In particular, for $i\in I_{\sigma}, k\geq 0, a\in\CC^*, p\in\ZZ$, we have :
\begin{equation*}
\begin{split}
 \pi(\chi(\text{res}^\sigma(\tilde{W}_{k,a}^{(\sigma^p(i))}))) = 
\begin{cases}
   \chi^\sigma(\text{res}^\sigma(W_{k,(a\omega^p)^{d_i}}^{(i)})) &\text{ if $\hat{\Glie}^\sigma$ is not of type $A_{2n}^{(2)}$,}
\\ \chi^\sigma(\text{res}^\sigma(W_{k,(-1)^pa}^{(i)}))&\text{ if $\hat{\Glie}^\sigma$ is of type $A_{2n}^{(2)}$.}
\end{cases}
\end{split}
\end{equation*}
\end{cor}

Let us define $\overline{\pi} : \ZZ[y_i^{\pm 1}]_{i\in I} \rightarrow \ZZ[z_i^{\pm 1}]_{i\in I_{\sigma}}$ as the ring morphism such that 
for $i\in I$, 
\begin{equation*}
\begin{split}
\overline{\pi}(y_i) = 
\begin{cases}z_{\overline{i}}&\text{ if $(\hat{\Glie}^\sigma,i)\neq (A_{2n}^{(2)},0)$,}
\\z_0^2&\text{ if $(\hat{\Glie}^\sigma,i) = (A_{2n}^{(2)},0)$.}
\end{cases}
\end{split}
\end{equation*} 
From Proposition \ref{dagdeux}, we have :
\begin{cor}\label{sqdeux}The following diagram is commutative :
$$\begin{CD}
\text{Rep}(\U_q(\Lo\Glie))   @>{\chi\circ\text{res}}>>  \ZZ[y_i^{\pm 1}]_{i\in I}\\
@VV{\pi_{R,a}}V     @VV{\overline{\pi}}V\\
\text{Rep}(\U_q(\Lo\Glie^\sigma))   @>{\overline{\chi}^\sigma\circ\overline{\text{res}}^\sigma}>>  \ZZ[z_i^{\pm 1}]_{i\in I_\sigma}
\end{CD}.$$
In particular, for $i\in I_{\sigma}, k\geq 0, a\in\CC^*,p\in\ZZ$, we have :
\begin{equation*}
\begin{split}
 \overline{\pi}(\chi(\text{res}^\sigma(\tilde{W}_{k,a}^{(\sigma^p(i))}))) = 
\begin{cases}
  \overline{\chi}^\sigma(\overline{\text{res}}^\sigma(W_{k,(a\omega^p)^{d_i}}^{(i)})) &\text{ if $\hat{\Glie}^\sigma$ is not of type $A_{2n}^{(2)}$,}
\\ \overline{\chi}^\sigma(\overline{\text{res}}^\sigma(W_{k,(-1)^pa}^{(i)}))&\text{ if $\hat{\Glie}^\sigma$ is of type $A_{2n}^{(2)}$.}
\end{cases}
\end{split}
\end{equation*}
\end{cor}

Note that $\pi$ (resp. $\overline{\pi}$) makes sense as a map from $P$ to the intergral weight lattice of $\U_q(\Glie^\sigma)$ (resp. $\overline{\U}_q(\Glie^\sigma)$). By abuse of notation, we also denote this map by $\pi$ (resp. $\overline{\pi}$).

\subsection{The twisted Kirillov-Reshetikhin conjecture}

As we have explicit formulas \cite{Nad, her06} for the character of tensor product of the Kirillov-Reshetikhin modules of $\U_q(\Lo\Glie)$, we also get formulas for the Kirillov-Reshetikhin modules of $\U_q(\Lo\Glie^\sigma)$ (in fact here we use the simply-laced cases so we only need the results of \cite{Nad}).
 
Let $\mathcal{Q}_k^{(i)}=e^{-k\Lambda_i}\chi^\sigma(Q_k^{(i)})$ and $\overline{\mathcal{Q}}_k^{(i)}=e^{-k\overline{\pi}(\Lambda_i)}\overline{\chi}^\sigma(\overline{res}^\sigma(W_{k,a}^{(i)}))$.

\begin{defi}\label{explicit} For a sequence $\nu = (\nu_k^{(i)})_{i\in I_\sigma, k> 0}$ such that for all but finitely many 
$\nu_k^{(i)}$ are non zero let us define :
$$\mathcal{F}(\nu)=\underset{N=(N_k^{(i)})}{\sum}\underset{i\in I, k>0
}{\prod}\begin{pmatrix}P_k^{(i)}(\nu,N)+N_k^{(i)}\\N_k^{(i)}\end{pmatrix}e^{-k N_k^{(i)}\pi(\alpha_i)},$$ 
$$\overline{\mathcal{F}}(\nu)=\underset{N=(N_k^{(i)})}{\sum}\underset{i\in I, k>0
}{\prod}\begin{pmatrix}P_k^{(i)}(\nu,N)+N_k^{(i)}\\N_k^{(i)}\end{pmatrix}e^{-k N_k^{(i)}\overline{\pi}(\alpha_i)},$$ 
where 
$$P_k^{(i)}(\nu,N)=\underset{l=1...\infty}{\sum}\nu_l^{(i)}\text{min}(k,l)-\underset{j\in
I,l>0}{\sum}N_l^{(j)}r_iC_{i,j}\text{min}(k/r_j,l/r_i),$$
$$\begin{pmatrix}a\\b\end{pmatrix}=\frac{\Gamma(a+1)}{\Gamma(a-b+1)\Gamma(b+1)}.$$\end{defi}

\noindent The above formulas are obtained via $\pi$ from the non-deformed fermionic formulas (we use here the version of \cite{ki1, ki2, hkoty, knt}; the version of \cite{kr} is slightly different because the definition of binomial coefficients is a little changed, see \cite{knt}).

\begin{thm}[The twisted Kirillov-Reshetikhin conjecture] For a sequence $\nu=(\nu_k^{(i)})_{i\in I_\sigma, k> 0}$ such that for all but finitely many 
$\nu_k^{(i)}$ are zero. Then we have :
$$\mathcal{Q}_{\nu}\underset{\alpha\in\Delta_+}{\prod}(1-e^{-\pi(\alpha)}) = \mathcal{F}(\nu)\text{ , }\overline{\mathcal{Q}}_{\nu}\underset{\alpha\in\Delta_+}{\prod}(1-e^{-\overline{\pi}(\alpha)}) = \overline{\mathcal{F}}(\nu)$$
where
$$\mathcal{Q}_{\nu}=\underset{i\in I_\sigma, k\geq 1}{\prod}(\mathcal{Q}_k^{(i)})^{\nu_k^{(i)}}\text{ , }\overline{\mathcal{Q}}_{\nu}=\underset{i\in I_\sigma, k\geq 1}{\prod}(\overline{\mathcal{Q}}_k^{(i)})^{\nu_k^{(i)}}.$$
\end{thm}

\demo This result follows from Corollary \ref{squn} and the analog result for simply-laced untwisted case \cite{Nad} (see \cite{her06} for non-simply laced untwisted case).\qed

\noindent In particular the formula obtained for $Q_{\nu}=\underset{i\in I_\sigma, k\geq 1}{\prod}(\chi^\sigma(Q_k^{(i)}))^{\nu_k^{(i)}}$ is the character of a $\U_q(\Glie^\sigma)$-module (which is a purely combinatorial statement and is not clear a priori).

\noindent Note that for type $A_{2n}^{(2)}$, some formulas were known \cite{oss} when $\nu_k^{(a)} = 0$ for $k\geq 2, a\in\CC^*$.

\section{Branching rules}\label{branch}

Now we enter the second part of the present paper : from the general results of the first part we get as an application explicit formulas (several of them had been conjectured by different authors). 

First we get the proof of the branching rules conjectured in \cite{hkott} for the subalgebras of finite type $\tilde{\U}_q(\Glie^\sigma),\overline{\U}_q(\Glie^\sigma)\subset\U_q(\Lo\Glie^\sigma)$. 

We use the notations of Section \ref{un} for $\tilde{\U}_q(\Glie^\sigma), \overline{\U}_q(\Glie^\sigma)$, and for $\lambda\in P^+$, we denote by $V(\lambda)$ (resp. $\overline{V}(\lambda)$) the simple $\tilde{\U}_q(\Glie^\sigma)$-module (resp. $\overline{\U}_q(\Glie^\sigma)$-module) of highest weight $\lambda$.
See Section \ref{subalg} for the type of $\tilde{U}_q(\Glie^\sigma)$ and $\overline{\U}_q(\Glie^\sigma)$. For type $\hat{\Glie}^\sigma \neq A_{2n}^{(2)}$, as $\U_q(\Glie^\sigma) \simeq \overline{\U}_q(\Glie^\sigma)$, we do not write the branching rules twice. The branching rules correspond to the decomposition obtained in \cite{ns} for conjectural crystals of Kirillov-Reshetikhin modules.

\subsection{Fundamental representations}

Let us start with fundamental representations. For these representations, the branching rules were obtained in \cite{hn} (in several cases it was already known, see the references in \cite{hn}), except for the node $3$ of $E_6^{(2)}$. These branching rules will also be obtained from explicit formulas of twisted $q$-characters in section \ref{expfund}. The two methods are very different (and the result match). 

Let us give these branching rules as this is the starting point for the proof of Theorem \ref{brkr} (we also give the dimension $D_i$ of the fundamental representation $V_i(a)$) :

\begin{thm}\label{br} We have the following branching rules for fundamental representations :

\noindent Type $A_{2n}^{(2)}$ and $0\leq i\leq n-1$ : 
$$\text{res}^\sigma(V_i(a)) = V(\Lambda_{n-i}) \oplus V(\Lambda_{n-i-1}) \oplus \cdots \oplus V(\Lambda_1) \oplus V(0)\text{ , }D_i = \begin{pmatrix}2n + 1\\n-i\end{pmatrix}.$$
$$\overline{\text{res}}^\sigma(V_i(a)) = \overline{V}(\Lambda_i)\text{ , }D_i = \begin{pmatrix}2n + 1\\n-i\end{pmatrix}.$$
Type $A_{2n-1}^{(2)}$ and $1\leq i\leq n$ :
$$\text{res}^\sigma(V_i(a)) = V(\Lambda_i) \oplus V(\Lambda_{i-2})\oplus \cdots \oplus V(\delta_{i,0[2]}\Lambda_0)\text{ , }D_i = \begin{pmatrix}2n \\i\end{pmatrix}.$$
Type $D_{n+1}^{(2)}$ and $1\leq i\leq n - 1$ :
$$\text{res}^\sigma(V_i(a)) = V(\Lambda_i) \oplus V(\Lambda_{i-1}) \oplus \cdots \oplus V(\Lambda_1) \oplus V(0),$$
$$D_i = \begin{pmatrix}2n+1\\i\end{pmatrix} + \sum_{k\geq 1} \begin{pmatrix}2n+2\\i-2k+1\end{pmatrix},$$
$$\text{res}^\sigma(V_n(a)) = V(\Lambda_n)\text{ , }D_n = 2^n.$$
Type $D_4^{(3)}$ :
$$\text{res}^\sigma(V_1(a)) = V(\Lambda_1) \oplus V(0)\text{ , }D_1 = 8,$$
$$\text{res}^\sigma(V_2(a)) = V(\lambda_2) \oplus V(\lambda_1)^{\oplus 2} \oplus V(0)\text{ , }D_2 = 29,$$
Type $E_6^{(2)}$ :
$$\text{res}^\sigma(V_1(a)) = V(\Lambda_1) \oplus V(0)\text{ , }D_1 = 27,$$
$$\text{res}^\sigma(V_2(a)) = V(\Lambda_2) \oplus V(\Lambda_4)\oplus V(\Lambda_1)^{\oplus 2}\oplus V(0)\text{ , }D_2 = 378,$$
$$\text{res}^\sigma(V_4(a)) = V(\Lambda_4) \oplus V(\Lambda_1)\oplus V(0)\text{ , }D_4 = 79.$$
\end{thm}

For the remaining fundamental representation in type $E_6^{(2)}$, we have $D_3 = 3732$ (see section \ref{expfund}) and it should be possible to check with a computer from the result of Section \ref{expfund} that the following conjectural formula of \cite{hkott} is satisfied :
$$\text{res}^\sigma(V_3(a)) = V(\Lambda_3) \oplus V(\Lambda_1 + \Lambda_4) \oplus V(2\Lambda_1) \oplus V(\Lambda_2)^{\oplus 3} \oplus V(\Lambda_4)^{\oplus 3} \oplus V(\Lambda_1)^{\oplus 4} \oplus V(0)^{\oplus 2}.$$

\subsection{Kirillov-Reshetikhin modules}

For Kirillov-Reshetikhin modules, conjectural branching rules are given in \cite{hkott}. It is proved \cite[Theorem 6.2]{hkott} that these formulas satisfy the twisted $Q$-system. As we have proved the twisted $Q$-system and that the formulas match for $k=0$ (trivial) and $k=1$ (Theorem \ref{br}), we get by induction on $k\geq 0$ the following :

\begin{thm}\label{brkr} For $a\in\CC^*$, $k\geq 0$, we have the following branching rules :

\noindent Type $A_{2n}^{(2)}$. For $0\leq i\leq n-1$ :
$$\text{res}^\sigma(W_{k,a}^{(i)}) = \bigoplus_{m_1\geq 0,\cdots, m_{n-i}\geq 0\text{ and }m_1+\cdots +m_{n-i} \leq k} V(m_1\Lambda_1 + m_2\Lambda_2 + \cdots + m_{n-i} \Lambda_{n-i}).$$
Type $A_{2n}^{(2)}$. For $1\leq i\leq n - 1$ :
$$\overline{\text{res}}^\sigma(W_{k,a}^{(i)}) = \bigoplus_{m_{n-1}\geq 0,\cdots, m_i\geq 0\text{ and }m_{n-1}+\cdots +m_i \leq k\text{ and }m_j = k\delta_{i,j}[2]} \overline{V}(m_{n-1}\Lambda_{n-1} + \cdots + m_i \Lambda_i),$$
and $\overline{\text{res}}^\sigma(W_{k,a}^{(0)})$ is equal to
$$\bigoplus_{m_0\geq 0,\cdots, m_{n-1}\geq 0\text{ and }m_0+\cdots +m_{n-1} \leq k\text{ and }m_j = k\delta_{i,0}[2]} V(m_1\Lambda_1 +  \cdots + m_{n-1} \Lambda_{n-1} + 2m_0 \Lambda_0).$$
Type $A_{2n-1}^{(2)}$. For $1\leq i\leq n$ and $j\in\{0,1\}$ such that $i = j\text{ mod }[2]$ :
$$\text{res}^\sigma(W_{k,a}^{(i)}) = \bigoplus_{m_j\geq 0,\cdots, m_i\geq 0\text{ and }m_j+\cdots +m_i = k} V(m_j\Lambda_j + m_{j+2}\Lambda_{j+2} + \cdots + m_i \Lambda_i).$$
Type $D_{n+1}^{(2)}$. For $1\leq i\leq n - 1$ :
$$\text{res}^\sigma(W_{k,a}^{(i)}) = \bigoplus_{m_1\geq 0,\cdots, m_i\geq 0\text{ and }m_1+\cdots +m_i \leq k} V(m_1\Lambda_1 + m_2\Lambda_2 + \cdots + m_i \Lambda_i).$$
$$\text{res}^\sigma(W_{k,a}^{(i)}) = V(k\Lambda_n).$$
Type $D_4^{(3)}$.
$$\text{res}^\sigma(W_{k,a}^{(1)}) = \bigoplus_{m = 0\cdots k} V(m\Lambda_1),$$
and $\text{res}^\sigma(W_{k,a}^{(2)})$ is equal to :
$$\bigoplus_{m_1 + m_2\leq k\text{ and }m_1, m_2\geq 0} (m_1 + 1)\text{min}(1 + m_2, 1+ k -m_1 - m_2)V(m_1\Lambda_1 + m_2\Lambda_2).$$
\end{thm}

\section{Explicit formulas for Kirillov-Reshetikhin modules}\label{expkr}

In this section we prove explicit formulas for the twisted $q$-characters of Kirillov-Reshetikhin in types $A_{2n}^{(2)}$, $A_{2n-1}^{(2)}$, $D_4^{(3)}$, $D_4^{(2)}$. Several of them had been conjectured in various papers. 

Remark : The formulas are given in terms of tableaux. It should be possible to directly compare them with the solutions of the twisted $T$-system given in \cite{t} in terms of determinant as this is done in \cite{nn1, nn2} for untwisted cases.

\subsection{Twisted $q$-characters of Kirillov-Reshetikhin modules in type $A_n^{(2)}$ ($n\geq 2$)}\label{krexp}

The following formulas appeared as combinatorial solutions of the twisted $T$-system in \cite{ks}. Here we prove that they are the twisted $q$-characters of Kirillov-Reshetikhin modules. In particular we get a new proof of the combinatorial statement that they satisfy the twisted $T$-system.

As for Kirillov-Reshetikhin modules in the untwisted case $A^{(1)}$ explicit formulas are known (see \cite{che1, Fre3, nt}), by using Theorem \ref{transi} we get the formulas for types $A_{2n-1}^{(2)}$ and types $A_{2n}^{(2)}$ :

\subsubsection{Type $A_{2n - 1}^{(2)}$ ($n\geq 2$)}

For $a\in\CC^*$ and $1\leq i\leq 2n$, let
\begin{equation*}
\begin{split}
\ffbox{i}_a = 
\begin{cases} Z_{1,a} &\text{ if $i = 1$,}
    \\Z_{i-1,aq^i}^{-1} Z_{i,aq^{i - 1}} &\text{ if $2 \leq i \leq n - 1$,}
    \\Z_{n-1,aq^n}^{-1} Z_{n,a^2q^{2n - 2}} &\text{ if $i = n$,}
    \\Z_{n,a^2q^{2n+2}}^{-1} Z_{n - 1,-aq^n} &\text{ if $i = n+1$,}
    \\Z_{2n - i + 1,-aq^i}^{-1} Z_{2n - i,-aq^{i - 1}} &\text{ if $n+2 \le i \le 2n - 1$,}
    \\Z_{1,-aq^{2n}}^{-1} &\text{ if $i=2n$.}
\end{cases}
\end{split}
\end{equation*}

For $1\leq i_0\leq n$, let $\text{Tab}(i_0,k)$ be the set of tableaux $(T_{i,j})_{1\leq i\leq i_0, 1\leq j\leq k}$ with coefficients in $\{1,\cdots,2n\}$ satisfying the two conditions :
\begin{itemize}
\item $T_{i,j}\leq T_{i,j+1}$ for any $1\leq i\leq i_0$ and $1\leq j\leq k-1$,

\item $T_{i,j} < T_{i+1,j}$ for any $1\leq i\leq i_0 - 1$ and $1\leq j\leq k$.
\end{itemize}

For such a tableaux $T\in\text{Tab}(i_0,k)$ and $a\in\CC^*$ we set 
$$m_{T,a} = \prod_{1\leq i\leq i_0,1\leq j\leq k} \ffbox{T_{i,j}}_{aq^{2(j-i)}}.$$

\begin{prop} For $a\in\CC^*$ and $1\leq i_0\leq n-1$, we have :
$$\chi_q^\sigma(W_{k,a}^{(i_0)}) = \sum_{T\in \text{Tab}(i_0,k)} m_{T,aq^{i_0-1}}.$$
Let $b\in\CC^*$ such that $b^2 = a$. Then we have :
$$\chi_q^\sigma(W_{k,a}^{(n)}) = \sum_{T\in \text{Tab}(i_0,k)} m_{T,bq^{n-1}}.$$
\end{prop}

\subsubsection{Type $A_{2n}^{(2)}$}

For $a\in\CC^*$ and $1\leq i\leq 2n+1$, let
\begin{equation*}
\begin{split}
 \ffbox{i}_a = \begin{cases} Z_{n-1,a}&\text{ if $i=1$,} 
    \\
     Z_{n-i+1,aq^i}^{-1} Z_{n-i,aq^{i - 1}} &\text{ if $2 \leq i \leq n$,}
    \\
    Z_{0,aq^{n+1}}^{-1} Z_{0,-aq^n} &\text{ if $i=n+1$,}
    \\
    Z_{i-n-2,-a q^i}^{-1} Z_{i-n-1,-aq^{i - 1}} &\text{ if $n+2 \leq i \leq 2n$,}
    \\
    Z_{n-1,-aq^{2n+1}}^{-1} &\text{ if $i=2n+1$.}
    \end{cases}
\end{split}
\end{equation*}

For $0\leq i_0\leq n-1$, let $\text{Tab}'(i_0,k)$ be the set of tableaux $(T_{i,j})_{1\leq i\leq n-i_0, 1\leq j\leq k}$ with coefficients in $\{1,\cdots,2n+1\}$ satisfying the conditions :
\begin{itemize}
\item $T_{i,j}\leq T_{i,j+1}$ for any $1\leq i\leq n-i_0$ and $1\leq j\leq k-1$,

\item $T_{i,j} < T_{i+1,j}$ for any $1\leq i\leq n-i_0 - 1$ and $1\leq j\leq k$.
\end{itemize}

For such a tableaux $T\in\text{Tab}'(i_0,k)$ and $a\in\CC^*$ we set 
$$m_{T,a} = \prod_{1\leq i\leq n-i_0,1\leq j\leq k} \ffbox{T_{i,j}}_{aq^{2(j-i)}}.$$

\begin{prop} For $a\in\CC^*$, $0\leq i_0\leq n-1$, we have :
$$\chi_q^\sigma(W_{k,a}^{(i_0)}) = \sum_{T\in \text{Tab}'(i_0,k)} m_{T,aq^{n-i_0-1}}.$$
\end{prop}

\subsection{Twisted $q$-characters of Kirillov-Reshetikhin modules in type $D_4^{(3)}$ and conjecture of \cite{ks}}\label{dtrois}

\subsubsection{The formulas}

The following formulas appeared in \cite{ks} as combinatorial conjectural solutions of the twisted $T$-system. In particular we prove the conjecture that they are solution of the twisted $T$-system.

For $a\in\CC^*$, let
{\allowdisplaybreaks
\begin{equation*}
\begin{aligned}[c]
    & \ffbox{1}_a = Z_{1,a}, 
    \\
    & \ffbox{2}_a = Z_{1,aq^2}^{-1} Z_{2,a^3q^3},
    \\
    & \ffbox{3}_a = Z_{2,a^3q^9}^{-1} Z_{1,aq^2j} Z_{1,aq^2j^2},
    \\
    & \ffbox{4}_a = Z_{1,aq^2j}Z_{1,aq^4j^2}^{-1},
    \\
    & \ffbox{\overline{4}}_a = Z_{1,aq^2j^2}Z_{1,aq^4j}^{-1},
    \\    
    & \ffbox{\overline{3}}_a = Z_{1,aq^4j}^{-1}Z_{1,aq^4j^2}^{-1}Z_{2,a^3q^9},
    \\
    & \ffbox{\overline{2}}_a = Z_{1,aq^4}Z_{2,a^3q^{15}}^{-1},
    \\
    & \ffbox{\overline{1}}_a = Z_{1,aq^6}^{-1}.
\end{aligned}
\end{equation*}}

Let $\mathbf B = \{ 1, 2, 3, 4, \overline{4}, \overline{3}, \overline{2}, \overline{1}\}$. 
We give the ordering $\prec$ on the set $\mathbf B$ by
\begin{equation*}
  1 \prec 2 \prec  3 \prec \begin{matrix}4\\\overline{4}\end{matrix}\prec \overline{3} \prec \overline{2}\prec \overline{1}.
\end{equation*}

Let $\text{Tab}(1,k)$ be the set of tableaux $(T_j)_{1\leq j\leq k}$ with coefficients in $\mathbf B$ satisfying $T_j\preceq T_{j+1}$ for any $1\leq j\leq k-1$.

For such a tableaux and $a\in\CC^*$ we set 
$$m_{T,a}^{(1)} = \prod_{1\leq j\leq k} \ffbox{T_j}_{aq^{2(j-1)}}.$$

\begin{thm}\label{explidun} For $a\in\CC^*$, $k\geq 1$ we have :
$$\chi_q^\sigma(W_{k,a}^{(1)}) = \sum_{T\in \text{Tab}(1,k)} m_{T,a}^{(1)}.$$
\end{thm}

This result will be proved in section \ref{pun}.

Let $\text{Tab}(2,k)$ be the set of tableaux $(T_{i,j})_{1\leq i\leq 2, 1\leq j\leq k}$ with coefficients in $\mathbf B$ satisfying the following conditions :
\begin{itemize}
\item($\theta$1) $(T_{i,1},\cdots,T_{i,k})\in \text{Tab}(1,k)$ for $i=1,2$,

\item($\theta$2) $T_{1,j}\nsucceq T_{2,j}$ for any $1\leq j\leq k$,

\item($\theta$3) $\begin{pmatrix}T_{1,j} & T_{1,j+1} & \cdots & T_{1,j'}\\T_{2,j}&T_{2,j+1} &\cdots &T_{2,j'}\end{pmatrix} \neq \begin{pmatrix}3 & 3& \cdots & 3& 4\\4&\overline{3}&\cdots & \overline{3} & \overline{3}\end{pmatrix}$ for any $1\leq j < j'\leq k$,

\item($\theta$4) $\begin{pmatrix}T_{1,j} & T_{1,j+1} & \cdots & T_{1,j'}\\T_{2,j}&T_{2,j+1} &\cdots &T_{2,j'}\end{pmatrix} \neq \begin{pmatrix}3 & 3& \cdots & 3& \overline{4}\\\overline{4}&\overline{3}&\cdots & \overline{3} & \overline{3}\end{pmatrix}$ for any $1\leq j < j'\leq k$.
\end{itemize}

The explanation for conditions ($\theta$3) and ($\theta$4) will be given in Section \ref{pdeux}.

\begin{lem} Under the conditions ($\theta$1) and ($\theta$2), the conditions ($\theta$3) and ($\theta$4) are respectively equivalent to the conditions :

($\theta$3') $\begin{pmatrix}T_{1,j} & T_{1,j'}\\T_{2,j} &T_{2,j'}\end{pmatrix}\neq \begin{pmatrix} 3&4\\4 &\overline{3}\end{pmatrix}$ for any $1\leq j<j'\leq k$,

($\theta$4') $\begin{pmatrix}T_{1,j} & T_{1,j'}\\T_{2,j} & T_{2,j'}\end{pmatrix}\neq \begin{pmatrix}3 & \overline{4}\\\overline{4}&\overline{3}\end{pmatrix}$ for any $1\leq j<j'\leq k$.
\end{lem}

The conditions ($\theta$3') and ($\theta$4') were originally used in \cite{ks} instead of ($\theta$3) and ($\theta$4).

\mk

\demo It is clear that ($\theta$3')$\Rightarrow$($\theta$3) and ($\theta$3')$\Rightarrow$($\theta$3). Suppose that ($\theta$3) is satisfied and suppose that there are $1\leq j < j'\leq k$ such that 
$$\begin{pmatrix}T_{1,j} & T_{1,j'}\\T_{2,j} & T_{2,j'}\end{pmatrix} = \begin{pmatrix} 3 & \overline{4}\\\overline{4}&\overline{3}\end{pmatrix}.$$ 
We can suppose that $j'-j$ is minimal for this property. By ($\theta$3) we have  $j'-j\geq 2$. We have $3\preceq T_{1,j+1}\preceq 4$. If $T_{1,j+1} = 4$, we have $T_{2,j+1}=\overline{3}$. If $T_{1,j+1} = 3$ we have $T_{2,j+1}\neq 4$ and so $T_{2,j+1} = \overline{3}$. So $T_{2,j+1} = \cdots = T_{2,j'} = \overline{3}$. In the same way $T_{1,j'-1} = \cdots = T_{1,j} = 3$. Contradiction by ($\theta$3).

We get ($\theta$4)$\Rightarrow$($\theta$4)' in an analog way.
\qed

For such a tableaux and $a\in\CC^*$ we set 
$$m_{T,a}^{(2)} = \prod_{1\leq i\leq 2, 1\leq j\leq k} \ffbox{T_{i,j}}_{aq^{2(j-i)}}.$$

\begin{thm}\label{expliddeux} For $a\in\CC^*$, $k\geq 1$ we have :
$$\chi_q^\sigma(W_{k,a}^{(2)}) = \sum_{T\in \text{Tab}(2,k)} m_{T,bq}^{(2)},$$
where $b\in\CC^*$ satisfies $b^3 = a$. 
\end{thm}

This result will be proved in section \ref{pdeux}.

\subsubsection{Notations} For $i\in\mathbf B$, we denote by $\text{succ}(i)$ (resp. $\text{prec}(i)$) the set of minimal (resp. maximal) elements of $\{j\in \mathbf B| j\succ i\}$ (resp. $\{j\in \mathbf B| j\prec i\}$).

\begin{defi} 
The affine degree of a monomial $m$ is 
$$d(m) = \text{max}\{z_{i,a}(m)|i\in I_\sigma,a\in\CC^*,u_{i,a}(m)\geq 0\}.$$ 
A monomial is said to be thin if it has affine degree $1$.

The affine degree of a $\U_q^\sigma(\hat{\Glie})$-module $V$ is the maximal affine degree of the monomials occurring in $\chi_q^\sigma(V)$. A $\U_q^\sigma(\hat{\Glie})$-module $V$ is said to be thin if it has affine degree $1$.  
\end{defi}

The notion of affine degree will be used in the following proofs and will be investigated more systematically in another paper in relation to some other problems for representations of quantum affine algebras in the continuation of \cite{herma}.

For $m$ a monomial we denote 
$$(m)^{\pm} = \prod_{\{(i,a)\in I\times\CC^|\pm z_{i,a}(m)\geq 0\}}Z_{i,a}^{z_{i,a}(m)}.$$ 
For $m,m'$ monomials, we say that $(m)^+$ is partly canceled by $(m')^-$ if there is $(i,a)$ such that $z_{i,a}(m) > 0$ and $z_{i,a}(m') < 0$.

\subsubsection{Proof of Theorem \ref{explidun}}\label{pun} 

\begin{lem}\label{superun} 
Let $T\in \text{Tab}(1,k)$ and $a\in\CC^*$. Let $1\leq j \neq j'\leq k$, $\alpha = T_j$ and $\beta = T_{j'}$. Then $(\ffbox{\alpha}_{aq^{2(j-1)}})^-$ is not partly canceled by $(\ffbox{\beta}_{aq^{2(j'-1)}})^+$.
\end{lem}

\demo Indeed we would have $\beta \prec \alpha$ or $(\alpha,\beta)\in\{(3,\overline{3}),(2,\overline{2})\}$. In the first case $j' < j $, contradiction. If $(\alpha,\beta) = (2,\overline{2})$, then $j' > j$ and 
$$(aq^{2(j'-1)})q^4 = (aq^{2(j-1)})q^2 \Rightarrow j ' + 1 = j,$$ 
contradiction. If $(\alpha,\beta) = (3,\overline{3})$, then $j' > j$ and 
$$(aq^{2(j'-1)})^3q^9 = (aq^{2(j-1)})^3q^9\Rightarrow j ' = j,$$ 
contradiction.\qed

\begin{lem}\label{thinun} For $T\in \text{Tab}(1,k)$ and $a\in\CC^*$, the monomial $m_{T,a}^{(1)}$ is thin.
\end{lem}

\demo  Indeed let $1\leq j < j'\leq k$, $\alpha = T_{j}$ and $\beta = T_{j'}$. We suppose that $$(\ffbox{\alpha}_{aq^{2(j-1)}})^+ = (\ffbox{\beta}_{aq^{2(j'-1)}})^+ \neq 1.$$ 
We have $\alpha = \beta$ or $(\alpha,\beta)\in\{(1,\overline{2}),(2,\overline{3}),(3,4),(3,\overline{4})\}$. If $\alpha = \beta$, we have $j = j'$, contradiction. Otherwise as $j < j'$ and $\alpha \prec \beta$, the power $r$ of $q$ for the term with positive exponent $Z_{l,aq^r}$ occurring for $\beta$ is strictly larger than the one for $\alpha$, contradiction.

 $(\ffbox{\alpha}_{aq^{2(j'-1)}})^+$ is necessarily strictly larger than in $(\ffbox{\beta}_{aq^{2(j-1)}})^+$, contradiction.\qed

Let us complete the proof of Theorem \ref{explidun} :

\mk

\demo From Lemma \ref{superun} there is a unique dominant monomial in the formula and it corresponds to $(1,1,\cdots,1)\in\text{Tab}(1,k)$. So it suffices to prove that the formula is in $\text{Im}(\chi_q^\sigma)$. 

Note that is also follows from Lemma \ref{superun} that the map $T\mapsto m_{T,a}^{(i)}$ is injective.

Let $1\leq i\leq 2$. We want to give a decomposition as in Proposition \ref{jdecomp} for $J = \{i\}$. From Proposition \ref{aidesldeux} and Lemma \ref{thinun}, the $L_i(M)$ that should appear in this decomposition are thin. It suffices to prove that the set $\text{Tab}(1,k)$ is in bijection with a disjoint union of sets $\mathcal{M}(L_i(M))$ via $T\mapsto m_{T,a}^{(1)}$.
 
We define a partial ordering on $\text{Tab}(1,k)$ : for $T,T'\in\text{Tab}(1,k)$, we denote 
$$T\preceq T'\text{ if and only if }(T_j\preceq T_j'\text{ for any $1\leq j\leq k$}).$$
Consider the set 
$$\mathcal{M}_1 = \{T\in\text{Tab}(1,k)|\forall 1\leq j\leq k, T_j\in\{1, 3, \overline{2}\}\}.$$ 
Then by Lemma \ref{superun}, $\mathcal{M}_1$ is in bijection with the set of $1$-dominant monomials occurring in the formula. 

Let $T\in\mathcal{M}_1$ and $\tilde{T}\in\text{Tab}(1,k)$ obtained from $T$ by replacing $1$ (resp. $3$, $\overline{2}$) by $2$ (resp. $\overline{3}$, $\overline{1}$). Let 
$$\mathcal{M}_1(T) = \{T'\in\text{Tab}(1,k)|T\preceq T'\preceq \tilde{T}\}.$$ 
Consider the decomposition $m_T = m_1m_3m_{\overline{2}}$ where $m_1$ (resp. $m_3$, $m_{\overline{2}}$) is obtained from $m_T$ by using only the boxes of values $1$ (resp. $3$, $\overline{2}$). From Proposition \ref{aidesldeux} we have 
$$L_1(m_T) = L_1(m_1)L_1(m_3)L_1(m_{\overline{2}}),$$ 
and $L_1(m_T)$ is thin. So $\mathcal{M}_1(T)$ is in bijection with the set of monomials of $L_1(m_T)$. In particular if $T\neq T'\in\mathcal{M}_1$ then $\mathcal{M}_1(T)$ and $\mathcal{M}_1(T')$ are disjoint. 

Moreover $(\mathcal{M}_1(T))_{T\in\mathcal{M}_1}$ is a partition of $\text{Tab}(1,k)$. Indeed for $T\in\text{Tab}(1,k)$ let $j_1$ minimal such that $T_{j_1}\succeq 3$ and $j_2$ minimal such that $T_{j_2}\succeq \overline{2}$. Consider $T'$ defined by 
\begin{equation*}
T_j' = 
\begin{cases}1 \text{ if $j < j_1 $,} 
\\3 \text{ if $j_1\leq j < j_2$,}
\\\overline{2}\text{ if $j\geq j_2$.}
\end{cases} 
\end{equation*}
Then $T'\in\mathcal{M}_1$ and $T$ is in $\mathcal{M}_1(T')$.

Consider the set 
$$\mathcal{M}_2 = \{T\in\text{Tab}(1,k)|\forall 1\leq j\leq k, T_j\in\{1, 2, 4,\overline{4}, \overline{3}, \overline{1}\}\}.$$ 
Then by Lemma \ref{superun} $\mathcal{M}_2$ is in bijection with the set of $2$-dominant monomials occurring in the formula.

Let $T\in\mathcal{M}_2$ and $\tilde{T}\in\text{Tab}(1,k)$ obtained from $T$ by replacing $2$ (resp. $\overline{3}$) by $3$ (resp. $\overline{2}$). Let 
$$\mathcal{M}_2(T) = \{T'\in\text{Tab}(1,k) | T\preceq T'\preceq \tilde{T}\}.$$ 
Consider the decomposition 
$$m_T = m_1m_2m_4m_{\overline{4}}m_{\overline{3}}m_{\overline{1}}$$ 
where $m_1$ (resp. $m_2$, $m_4$, $m_{\overline{4}}$, $m_{\overline{3}}$, $m_{\overline{1}}$) is obtained from $m_T$ by using only the boxes of values $1$ (resp. $2$, $4$, $\overline{4}$, $\overline{3}$, $\overline{1}$). From Proposition \ref{aidesldeux} we have $$L_2(m_T) = m_1 m_4 m_{\overline{4}}m_{\overline{1}}L_2(m_2)L_2(m_{\overline{3}})$$ 
and $L_2(m_T)$ is thin. So $\mathcal{M}_2(T)$ is in bijection with the set of monomials of $L_2(m_T)$. In particular if $T\neq T'\in\mathcal{M}_2$ then $\mathcal{M}_2(T)$ and $\mathcal{M}_2(T')$ are disjoint. 

Moreover $(\mathcal{M}_2(T))_{T\in\mathcal{M}_2}$ is a partition of $\text{Tab}(1,k)$. Indeed let $T\in\text{Tab}(1,k)$. We suppose that $\overline{4}$ does not occur in $T$ (we can treat in a similar way the case when $4$ does not occur in $T$). Then we wan construct $T'\in\mathcal{M}_2$ such that $T\in\mathcal{M}_2(T')$, by analogy to the case of $\mathcal{M}_1$, by using $(1,2,4,\overline{3},\overline{1})$ instead of $(1,3,\overline{2})$.
\qed

\subsubsection{Proof of Theorem \ref{expliddeux}}\label{pdeux}

\begin{lem}\label{superdeux} Let $T\in \text{Tab}(2,k)$ and $a\in\CC^*$. Let $1\leq i,i'\leq 2$, $1\leq j, j'\leq k$, $\alpha = T_{i,j}$ and $\beta = T_{i',j'}$. If $(\ffbox{\alpha}_{aq^{2(j-i)}})^-$ is partly canceled by $(\ffbox{\beta}_{aq^{2(j'-i')}})^+$, then one of the following condition is satisfied : 
\begin{itemize}
\item $i=2$, $i'=1$, $j' = j$ and ($\alpha\in\text{Succ}(\beta)$ or $(\alpha,\beta) = (3,\overline{3})$).

\item $i=2$, $i'=1$, $j' = j+1$ and $(\alpha,\beta)=(\overline{2},2)$.

\item $i=2$, $i'=1$, $j' = j+2$ and $(\alpha,\beta)=(\overline{1},1)$.

\item $i=1$, $i'=2$, $j = j'$ and $(\alpha,\beta) = (2,\overline{2})$.

\item $i=1$, $j=2$, $j = j' - 1$ and $(\alpha,\beta) = (3,\overline{3})$.
\end{itemize}
\end{lem}

\demo First from Lemma \ref{superun}, we have necessarily $i\neq i'$. 

We have $\alpha\neq\beta$. Moreover 
$$(\alpha\prec\beta\Leftrightarrow (j-i)\geq (j'-i'))\text{ and }(\beta\prec\alpha\Leftrightarrow (j-i) < (j'-i')).$$
Suppose that $i=2$ and $i'=1$.

If $j'< j$, we have $\beta\prec \alpha$. As $(j-i) \geq (j'-i')$, contradiction.

If $j' \geq j$, $(j'-i')\geq (j-i) + 1$, so $\beta\prec \alpha$. 

Suppose that $i=1$ and $i'=2$.

If $j < j'$, we have $\alpha\prec\beta$ and $(j-i)\leq (j'-i')$. So $j-i=j'-i'$ and $(\alpha,\beta) = (3,\overline{3})$. In particular $j=j'-1$.

If $j\geq j'$, we have $(j-i) \geq (j'-i') + 1$. So $\alpha\prec\beta$.\qed

\begin{lem}\label{degdeux} Let $T\in\text{Tab}(2,k)$, $\alpha = T_{i,j}$ and $\beta = T_{i',j'}$ where $1\leq i\leq i'\leq 2$ and $1\leq j,j'\leq k$. Suppose that 
$$(\ffbox{\alpha}_{aq^{2(j-i)}})^+ = (\ffbox{\beta}_{aq^{2(j'-i')}})^+\neq 1,$$ 
and denote by $M$ this monomial. Then $i=1$ and $i'=2$ and one of the following conditions is satisfied :
\begin{itemize}
\item($j = j'$ and $(\alpha,\beta)=(2,\overline{3})$).

\item($j-j'=1$ and $(\alpha,\beta) = (1,\overline{2})$). 

\item($j'-j=1$ and $\alpha =3$ and $\beta\in\{4,\overline{4}\}$). 
\end{itemize}
In the first case $M$ is $2$-dominant, and in the last two cases $M$ is $1$-dominant.
\end{lem}

\demo By Lemma \ref{thinun}, we have $i\neq i'$. So $i=1$ and $i'=2$. We have $\alpha = \beta$ or $$\{\alpha,\beta\}\in\{\{1,\overline{2}\},\{2,\overline{3}\},\{3,4\},\{3,\overline{4}\}\}.$$ 

If $\alpha = \beta$, we have $j+1 = j'$, contradiction as $T_{1,j}\prec T_{2,j+1}$.

If $(\alpha,\beta) = (1,\overline{2})$, we have $2(j-1) =2(j'-2)+4$ and $j' = j - 1$.

If $(\alpha,\beta) = (\overline{2},1)$, we have $2(j-1)+4 =2(j'-2)$ and $j' = j + 3$, contradiction, as $T_{1,j}\prec T_{2,j'}$. 

If $(\alpha,\beta) = (2,\overline{3})$, we have $6(j-1) +3= 6(j'-2)+9$ and $j=j'$.

If $(\alpha,\beta) = (\overline{3},2)$, we have $6(j-1) +9= 6(j'-2)+j'$ and $j'=j+3$, contradiction as $T_{1,j}\preceq T_{2,j'}$.

If $(\alpha,\beta) = (3,4)$, we have $j-1=j'-2$ and $j'=j+1$.

If $(\alpha,\beta) = (4,3)$, we have $j-1=j'-2$ and $j'=j+1$, contradiction as $T_{1,j}\prec T_{2,j+1}$.

The cases $(\alpha,\beta) = (3,\overline{4})$ and $(\alpha,\beta) = (\overline{4},3)$ are treated in an analog way.
\qed

We define a partial ordering on $\text{Tab}(2,k)$ : 

\begin{defi} For $T,T'\in\text{Tab}(2,k)$, we denote 
$$T\preceq T'\text{ if and only if }(\forall 1\leq j\leq k, T_{1,j}\preceq T_{1,j}'\text{ and }T_{2,j}\preceq T_{2,j}').$$
\end{defi}

\begin{lem}\label{prevun} Consider the elements of $\text{Tab}(2,k)$: 
$$T = \begin{pmatrix} 3 & 3 &\cdots & 3\\ 4 & 4 &\cdots & 4\end{pmatrix}\text{ , }\tilde{T} = \begin{pmatrix} \overline{4} & \overline{4} &\cdots & \overline{4}\\ \overline{3} & \overline{3} &\cdots & \overline{3}\end{pmatrix}\text{ , }T_r = \begin{pmatrix} 3 & 3 &\cdots & 3\\ \overline{3} & \overline{3} &\cdots & \overline{3}\end{pmatrix}.$$
Then we have :
$$\sum_{\{T'\in\text{Tab}(2,k)|T\preceq T'\preceq \tilde{T}\}} m_{T'} = L_1(m_T) + m_{T_r} = L_1(m_T) + L_1(m_{T_r}).$$
We have the same result with 
$$T = \begin{pmatrix} 3 & 3 &\cdots & 3\\ \overline{4} & \overline{4} &\cdots & \overline{4}\end{pmatrix}\text{ , }\tilde{T} = \begin{pmatrix} 4 & 4 &\cdots & 4\\ \overline{3} & \overline{3} &\cdots & \overline{3}\end{pmatrix}\text{ , }T_r = \begin{pmatrix} 3 & 3 &\cdots & 3\\ \overline{3} & \overline{3} &\cdots & \overline{3}\end{pmatrix},$$
and
$$T = \begin{pmatrix} 1 & 1 &\cdots & 1\\ \overline{2} & \overline{2} &\cdots & \overline{2}\end{pmatrix}\text{ , }\tilde{T} = \begin{pmatrix} 2 & 2 &\cdots & 2\\ \overline{1} & \overline{1} &\cdots & \overline{1}\end{pmatrix}\text{ , }T_r = \begin{pmatrix} 2 & 2 &\cdots & 2\\ \overline{2} & \overline{2} &\cdots & \overline{2}\end{pmatrix}.$$
\end{lem}

\demo We study the first case (the two other cases are analog). The monomials $(m_T)^{(1)}$ is of the form $m_{b,k+1}^{(1)}m_{bq^2,k-1}^{(1)}$ where $b\in\CC^*$. From Proposition \ref{aidesldeux} (2), the unique monomial of the form $m_{T'}$ where $T\preceq T'\preceq \tilde{T}$ that does not occur in $L_1(m_T)$ is the dominant monomial 
$$m_TA_{1,bq^{2k-1}}^{-1}A_{1,bq^{2k-3}}^{-1}\cdots A_{1,bq^3}^{-1} = m_{T'}.$$ 
\qed

Remark : Consider the tableaux $t = \begin{pmatrix} 3 & 3 &\cdots & 3&4\\ 4 & \overline{3} &\cdots & \overline{3}&\overline{3}\end{pmatrix}$ which is not in $\text{Tab}(2,k)$ by condition $(\theta_3)$. We have 
$$\begin{pmatrix} 3 & 3 &\cdots & 3\\ 4 & 4 &\cdots & 4\end{pmatrix}\preceq t \preceq \begin{pmatrix} \overline{4} & \overline{4} &\cdots & \overline{4}\\ \overline{3} & \overline{3} &\cdots & \overline{3}\end{pmatrix}.$$
Lemma \ref{prevun} gives an explanation for the condition $(\theta_3)$ : the tableaux $t$ does not correspond to a monomial in $L_1(m_T)$ or $L_1(m_{T_r})$. We have an analog remark for the condition $(\theta 4)$. As $3$ is the unique element of $\mathbf B$ satisfying $|\text{succ}(3)| \geq 2$, this kind of situation can only appear for the first two tableaux $T$ considered in Lemma \ref{prevun}. That is why no other condition than $(\theta_3)$ and $(\theta_4)$ are required in the definition of $\text{Tab}(2,k)$.

As for Lemma \ref{prevun}, we have (but here $2$-dominant monomials are involved instead of $1$-dominant monomials) :

\begin{lem}\label{prevdeux} Consider the elements of $\text{Tab}(2,k)$: 
$$T = \begin{pmatrix} 1 & 2 &\cdots &2&2\\\overline{3}&\overline{3}&\cdots &\overline{3}&\overline{1}\end{pmatrix}
\text{ , }
\tilde{T} = \begin{pmatrix} 1 & 3 &\cdots &3&3\\\overline{2}&\overline{2}&\cdots &\overline{2}&\overline{1}\end{pmatrix}
\text{ , }
T_r = \begin{pmatrix} 1 & 2 &\cdots &2&2\\\overline{2}&\overline{2}&\cdots &\overline{2}&\overline{1}\end{pmatrix}.$$
Then we have :
$$\sum_{\{T'\in\text{Tab}(2,k)|T\preceq T'\preceq \tilde{T}\}} m_{T'} = L_2(m_T) + m_{T_r} = L_2(m_T) + L_2(m_{T_r}).$$
We have the same result with 
$$T = \begin{pmatrix} 2 & 2 &\cdots &2&4\\ \overline{4}&\overline{3}&\cdots &\overline{3}&\overline{3}\end{pmatrix}
\text{ , }
\tilde{T} = \begin{pmatrix} 3 & 3 &\cdots &3&4\\ \overline{4}&\overline{2}&\cdots &\overline{2}&\overline{2}\end{pmatrix}
\text{ , }
T_r = \begin{pmatrix} 3 & 3 &\cdots & 3 & 4\\ \overline{4}&\overline{3}&\cdots &\overline{3}&\overline{3}\end{pmatrix},$$
and
$$T = \begin{pmatrix} 2 & 2 &\cdots &2&\overline{4}\\ 4&\overline{3}&\cdots &\overline{3}&\overline{3}\end{pmatrix}
\text{ , }
\tilde{T} = \begin{pmatrix} 3 & 3 &\cdots &3&\overline{4}\\ 4&\overline{2}&\cdots &\overline{2}&\overline{2}\end{pmatrix}
\text{ , }
T_r = \begin{pmatrix} 3 & 3 &\cdots & 3 & \overline{4}\\ 4&\overline{3}&\cdots &\overline{3}&\overline{3}\end{pmatrix}.$$
\end{lem}

We also have the following result which is a little bit different than Lemma \ref{prevdeux} because of conditions $(\theta 3)$ and $(\theta 4)$ :

\begin{lem}\label{prevtrois} Consider the elements of $\text{Tab}(2,k)$: 
$$T = \begin{pmatrix} 2 & 2 &\cdots &2&4\\ 4&\overline{3}&\cdots &\overline{3}&\overline{3}\end{pmatrix}
\text{ , }
\tilde{T} = \begin{pmatrix} 3 & 3 &\cdots &3&4\\ 4&\overline{2}&\cdots &\overline{2}&\overline{2}\end{pmatrix}
.$$
Then we have :
$$\sum_{\{T'\in\text{Tab}(2,k)|T\preceq T'\preceq \tilde{T}\}} m_{T'} = L_2(m_T).$$
We have the same result with 
$$T = \begin{pmatrix} 2 & 2 &\cdots &2&\overline{4}\\ \overline{4}&\overline{3}&\cdots &\overline{3}&\overline{3}\end{pmatrix}
\text{ , }
\tilde{T} = \begin{pmatrix} 3 & 3 &\cdots &3&\overline{4}\\ \overline{4}&\overline{2}&\cdots &\overline{2}&\overline{2}\end{pmatrix}
.$$
\end{lem}

Let us complete the proof of Theorem \ref{expliddeux} :

\mk

\demo Suppose that there is a dominant monomial $m_T^{(2)}$ where 
$$T\neq \begin{pmatrix}1&1&\cdots&1\\1&1&\cdots&1\end{pmatrix}\in\text{Tab}(2,k).$$ 
We have necessarily $\alpha = T_{2,k}\neq 1$. Suppose that $(\ffbox{\alpha}_{aq^{2(k-2)}})^-$ is partly canceled by $(\ffbox{\beta}_{aq^{2(j'-i')}})^+$ (let $n_1$ such that the canceled variable is of the form $Z_{l,(a\omega^s q^{n_1})^{d_l}}$). By Lemma \ref{superdeux}, $j'=k$, $i'=2$. Suppose that  $(\ffbox{\beta}_{aq^{2(k-1)}})^-$ is partly canceled by $(\ffbox{\gamma}_{aq^{2(j''-2)}})^+$ (let $n_1$ such that the canceled variable is of the form $Z_{l',(a\omega^{s'}q^{n_2})^{d_{l'}}}$). By Lemma \ref{superdeux}, $j''\leq j$ and $\beta\prec\gamma$. But we have $n_2 > n_1$. So $Z_{l',(a\omega^{s'}q^{n_2})^{d_{l'}}}$ appears in $(\ffbox{\gamma}_{aq^{2(j''-2)}})$. As $\gamma\preceq \alpha$, we have an element $Z_{l'',(a\omega^{s''}q^{n_3})^{d_{l''}}}$ which appears in $(\ffbox{\alpha}_{aq^{2(k-2)}})^-$ with $n_3 > n_1$, contradiction.

So there is a unique dominant monomial in the formula.

From Lemma \ref{degdeux}, the monomials $m_T^{(2)}$ have the maximal affine degree equal to $2$ (as one $(i,j)$ can only be involved in one of the relations described in Lemma \ref{degdeux}). 

Let $1\leq i\leq 2$. We want to give a decomposition as in Proposition \ref{jdecomp} for $J = \{i\}$. For $k=1$ the result is clear from the study of fundamental representations (see section \ref{fundhe}). We suppose that $k\geq 2$. From Proposition \ref{aidesldeux}, the $L_i(M)$ that should appear in this decomposition have affine degree at most equal to $2$ and several dominant $i$-dominant monomials may appear inside.

Consider the set $\mathcal{M}_1$ of tableaux $T\in\text{Tab}(2,k)$ satisfying
\begin{itemize}
\item for any $1\leq j\leq k$, 
$$(T_{1,j},T_{2,j})\in\{(1,2), (1,3), (1,\overline{2}), (3,\overline{2}),(3,4), (3,\overline{4}), (3,\overline{3}),(\overline{2},\overline{1}),(2,\overline{2})\},$$ 

\item $(3,\overline{3})$ and $(3,4)$ do not appear simultaneously,

\item $(3,\overline{3})$ and $(3,\overline{4})$ do not appear simultaneously,

\item $(2,\overline{2})$ and $(1,\overline{2})$ do not appear simultaneously.
\end{itemize}

\begin{rem}\label{aiderem} Note that the $T\in\mathcal{M}_1$ define $1$-dominant monomials $m_T$, but there are other elements of $\text{Tab}(2,k)$ with the same property. In fact we can obtain all the others in the following way. Let us start with $T\in\mathcal{M}_1$ such that the columns $(3,4)$ or $(3,\overline{4})$ or $(1,\overline{2})$ appear several times. Then we replace some $(3,4)$, some $(3,\overline{4})$, some $(1,\overline{2})$ respectively by $(3,\overline{3})$, $(3,\overline{3})$, $(2,\overline{2})$ but in each case not all of them. The $T'$ that we get by this process are all elements of $\text{Tab}(2,k)-\mathcal{M}_1$ corresponding to $1$-dominant monomials.\end{rem}

Let $T\in\mathcal{M}_1$ and $\tilde{T}\in\text{Tab}(2,k)$ obtained from $T$ by replacing $(1,3)$ (resp. $(1,\overline{2})$, $(3,\overline{2})$, $(3,4)$, $(3,\overline{4})$, $(2,\overline{2})$) by $(2,\overline{3})$ (resp. $(2,\overline{1})$, $(\overline{3},\overline{1})$, $(\overline{4},\overline{3})$, $(4,\overline{3})$, $(2,\overline{2})$). 

Let $\mathcal{M}_1(T)$ be the set of tableaux $T'\in\text{Tab}(1,k)$ satisfying
\begin{itemize}
\item $T\preceq T'\preceq \tilde{T}$,

\item if $(3,4)$ appears in $T$, then
$$(T_{1,l}',T_{2,l}')\neq (3,\overline{3})\text{ where }j = \text{Min}\{1\leq l\leq k|(T_{1,l},T_{2,l}) = (3,4)\},$$ 

\item if $(3,\overline{4})$ appears in $T$, then
$$(T_{1,j}',T_{2,j}')\neq (3,\overline{3})\text{ where }j = \text{Min}\{1\leq l\leq k| (T_{1,l},T_{2,l})=(3,\overline{4})\},$$

\item if $(1,\overline{2})$ appears in $T$, then 
$$(T_{1,j}',T_{2,j}')\neq (2,\overline{2})\text{ where }j = \text{Min}\{1\leq l\leq k| (T_{1,l},T_{2,l})=(1,\overline{2})\}.$$ 
\end{itemize}

In particular if $T\neq T'\in\mathcal{M}_1$ then $\mathcal{M}_1(T)$ and $\mathcal{M}_1(T')$ are disjoint. 

Let $T\in\mathcal{M}_1$ and let us prove that 
$$L_1(m_T)=\sum_{T'\in\mathcal{M}_1(T)}m_{T'}.$$
This can be done as for the proof of Theorem \ref{explidun}, except that by Lemma \ref{degdeux} $m_T$ might be not thin. By Lemma \ref{degdeux} this can happen when one of the columns $(3,4)$, $(3,\overline{4})$ or $(1,\overline{2})$ appears several times. These situations can be reduced to the cases studied in Lemma \ref{prevun}. In each case the monomial $m_{T_r}$ is precisely excluded of $\mathcal{M}_1(T)$ by the definition of $\mathcal{M}_1(T)$, but is in $\mathcal{M}_1$ (the other monomials excluded from $\mathcal{M}_1$ are in $\mathcal{M}_1(T)$, see Remark \ref{aiderem}).

Moreover $(\mathcal{M}_1(T))_{T\in\mathcal{M}_1}$ is a partition of $\text{Tab}(1,k)$. Indeed for $T\in\text{Tab}(1,k)$ by analogy to the proof of Theorem \ref{explidun}, it is easy to construct $T'\in\mathcal{M}_1$ such that $T\in\mathcal{M}_1(T')$. Here we use the partial ordering on ${\bf B}\times {\bf B}$ instead of the partial ordering on ${\bf B}$ used in the proof of Theorem \ref{explidun} : 
$$((i,j)\preceq (i',j'))\Leftrightarrow (i\preceq i'\text{ and }j\preceq j').$$
Consider the set $\mathcal{M}_2$ of tableaux $T\in\text{Tab}(2,k)$ satisfying for any $1\leq j\leq k$ :
\begin{itemize}
\item $T_{1,j}\neq \overline{2}$,

\item if $T_{2,j} = 3$ then $T_{1,j} = 2$,

\item if $T_{2,j} = \overline{2}$ and $T_{1,j} \neq \overline{3}$, then there are $1\leq j_1\leq j < j_2\leq k$ such that 
$$\begin{pmatrix}T_{1,j_1} & T_{1,j_1+1}&\cdots&T_{1,j_2}\\T_{2,j_1}&T_{2,j_1+1}&\cdots&T_{2,j_2} \end{pmatrix} = \begin{pmatrix}1 & 2 &\cdots & 2 & 2\\\overline{2}&\overline{2}& \cdots   &\overline{2}&\overline{1} \end{pmatrix},$$

\item if $T_{1,j} = 3$, then there are $1\leq j_1 \leq j < j_2 \leq k$ such that 
$$\begin{pmatrix}T_{1,j_1} & T_{1,j_1+1}&\cdots&T_{1,j_2}\\T_{2,j_1}&T_{2,j_1+1}&\cdots&T_{2,j_2} \end{pmatrix} = \begin{pmatrix}3 & 3 &\cdots & 3 & \alpha\\\beta&\overline{3}& \cdots   &\overline{3}&\overline{3} \end{pmatrix},$$
where $(\alpha,\beta) = (4,\overline{4})$ or $(\alpha,\beta) = (\overline{4}, 4)$.
\end{itemize}

\begin{rem}\label{aideremdeux}
Note that the $T\in\mathcal{M}_2$ define $2$-dominant monomials $m_T$, but there are other elements of $\text{Tab}(2,k)$ with the same property. In fact we can obtain all the others in the following way. Let us start with $T\in\mathcal{M}_2$ such that the columns $(2,\overline{3})$ and $(2,\overline{1})$ (resp. and $(\alpha,\overline{3})$ where $\alpha\in\{4,\overline{4}\}$) appear simultaneously. Then we replace some $(2,\overline{3})$ by $(2,\overline{2})$ (resp. by $(3,\overline{3})$) such that the $T'$ that we get is in $\text{Tab}(2,k)$. By this process we get all elements $T'\in\text{Tab}(2,k)-\mathcal{M}_2$ corresponding to $2$-dominant monomials.
\end{rem}

\noindent Let $T\in\mathcal{M}_2$ and $\tilde{T}\in\text{Tab}(2,k)$ obtained from $T$ by the following process :
\begin{itemize}
\item if $T_{2,j} = \overline{3}$ and $T_{1,j-1}\neq 3$, we replace this $\overline{3}$ by $\overline{2}$,

\item if $T_{2,j} = 2$, we replace this $2$ by $3$,

\item if $T_{1,j} = \overline{3}$ and $T_{2,j}\neq \overline{2}$, we replace this $\overline{3}$ by $\overline{2}$,

\item if $T_{1,j} = 2$ and $T_{2,j}\neq 3$ and $T_{2,j-1}\neq \overline{2}$, we replace this $2$ by $3$.
\end{itemize}

Let $T\in\mathcal{M}_2$. Then we define $\mathcal{M}_2(T)$ as the set of $T'\in\text{Tab}(2,k)$ satisfying :
\begin{itemize}
\item $T\preceq T'\preceq \tilde{T}$,

\item if $(1,\overline{3})$ and $(2,\overline{3})$ and $(2,\overline{1})$ appear simultaneously in $T$, for $j$ such that \\$\begin{pmatrix}T_{1,j} & T_{1,j+1}\\ T_{2,j} & T_{2,j+1}\end{pmatrix} = \begin{pmatrix} 1 & 2\\ \overline{3} & \overline{3}\end{pmatrix}$, we have ($T_{1,j+1}'=2\Rightarrow T_{2,j}'\neq \overline{2}$),

\item if $(2,\alpha)$ and $(2,\overline{3})$ and $(\beta,\overline{3})$ appear simultaneously in $T$ where $\alpha,\beta\in\{4,\overline{4}\}$, for $j$ such that $\begin{pmatrix}T_{1,j} & T_{1,j+1} \\T_{2,j} & T_{2,j+1}\end{pmatrix}=\begin{pmatrix} 2 & 2 \\\alpha & \overline{3}\end{pmatrix}$, we have ($T_{2,j+1}'=\overline{3}\Rightarrow T_{1,j}'\neq 3$).
\end{itemize}

Note that by $(\theta 3)$ and $(\theta 4)$, the last condition is automatic in the cases $\alpha = \beta$.

Then we can prove as above by using Lemma \ref{degdeux}, Lemma \ref{prevdeux}, Lemma \ref{prevtrois} and Remark \ref{aideremdeux} that for $T\in\mathcal{M}_2$ we have 
$$L_2(m_T)=\sum_{T'\in\mathcal{M}_2(T)}m_{T'},$$ 
and that $(\mathcal{M}_2(T))_{T\in\mathcal{M}_2}$ is a partition of $\text{Tab}(2,k)$.
\qed

For recent results on conjectural crystals (pseudo basis) for $i=1$ in type $D_4^{(3)}$ see \cite{kmoy}.

\begin{cor} 
$W_{k,a}^{(1)}$ is thin and $W_{k,a}^{(2)}$ has affine degree $2$.
\end{cor}

\demo
From the study of fundamental representations, the affine degree of $W_{k,a}^{(2)}$ is at least equal to $2$. As the monomials $m_T^{(2)}$ have the maximal affine degree equal to $2$ (see the proof of Theorem \ref{expliddeux}), the affine degree of $W_{k,a}^{(2)}$ is $2$.
\qed

\subsection{Relation to untwisted type $D_4^{(1)}$}

Our result for twisted type $D_4^{(3)}$ is coherent with explicit formulas conjectured \cite{nn1, nn2} for $q$-characters of Kirillov-Reshetikhin modules in untwisted types $D_4^{(1)}$. As it is proved in \cite{kns2} that these explicit formulas solve the $T$-system, these formulas were first proved in \cite{Nad} as a direct consequence of the proof of the $T$-system for simply-laced cases.

Here we explain the connection to the results of the present paper. The strategy of our proof for type $D_4^{(3)}$ is different as we do not use the $T$-system but we use directly the special property. Additional formulas conjectured in \cite{nn1, nn2} for type $D$ will be proved in a separate publication by using the results of \cite{herma} and this strategy. For type $D_4^{(1)}$ the special property of Kirillov-Reshetikhin modules was proved by Nakajima \cite{Nab, Nad}, and so to prove the explicit formulas for $D_4^{(1)}$ we could also have directly rewritten the proof of Section \ref{dtrois} in this context.

Let us consider the untwisted quantum affine algebra of type $D_4^{(1)}$.

For $a\in\CC^*$, let
{\allowdisplaybreaks
\begin{equation*}
\begin{aligned}[c]
    & \ffbox{1}_{u,a} = Y_{1,a}, 
    \\
    & \ffbox{2}_{u,a} = Y_{1,aq^2}^{-1} Y_{2,aq},
    \\
    & \ffbox{3}_{u,a} = Y_{2,aq^3}^{-1} Y_{4,aq^2} Y_{3,aq^2},
    \\
    & \ffbox{4}_{u,a} = Y_{4,aq^2}Y_{3,aq^4}^{-1},
    \\
    & \ffbox{\overline{4}}_{u,a} = Y_{3,aq^2}Y_{4,aq^4}^{-1},
    \\    
    & \ffbox{\overline{3}}_{u,a} = Y_{4,aq^4}^{-1}Y_{3,aq^4}^{-1}Y_{2,aq^3},
    \\
    & \ffbox{\overline{2}}_{u,a} = Y_{1,aq^4}Y_{2,aq^5}^{-1},
    \\
    & \ffbox{\overline{1}}_{u,a} = Y_{1,aq^6}^{-1}.
\end{aligned}
\end{equation*}}

We define $\text{Tab}(1,k)$ and $\text{Tab}(2,k)$ as in the previous section.

For $T\in\text{Tab}(1,k)$ and $a\in\CC^*$ we set 
$$m_{u,T,a}^{(1)} = \prod_{1\leq j\leq k} \ffbox{T_j}_{u,aq^{2(j-1)}}.$$
For $T\in\text{Tab}(2,k)$ and $a\in\CC^*$ we set 
$$m_{u,T,a}^{(2)} = \prod_{1\leq i\leq 2, 1\leq j\leq k} \ffbox{T_{i,j}}_{u,aq^{2(j-i)}}.$$
As a consequence of the results in \cite{kns2, nn1, Nad} :

\begin{thm}\label{expu} For $a\in\CC^*$, $k\geq 1$ we have :
$$\chi_q(W_{k,a}^{(1)}) = \sum_{T\in \text{Tab}(1,k)} m_{u,T,a}^{(1)}\text{ and }\chi_q(W_{k,a}^{(2)}) = \sum_{T\in \text{Tab}(2,k)} m_{u,T,aq}^{(2)}.$$ 
\end{thm}

Theorem \ref{transi}, Theorem \ref{explidun} and Theorem \ref{expliddeux} also imply this result.

Note that $\chi_q(W_{k,a}^{(3)})$ and $\chi_q(W_{k,a}^{(4)})$ are obtained in an analog way by permuting the numbering of the nodes $(2,3,4)$.

\begin{cor} 
$W_{k,a}^{(1)}$, $W_{k,a}^{(3)}$, $W_{k,a}^{(4)}$ are thin and $W_{k,a}^{(2)}$ has affine degree $2$.
\end{cor}

To check that we get the explicit formulas of \cite{nn2}, we see we have the same condition to define the set of tableaux. Indeed the condition of \cite{nn2} are :
\begin{itemize}
\item $T_{i,j}\preceq T_{i,j+1}$ or $\{T_{i,j},T_{i,j+1}\} = \{4,\overline{4}\}$ for $1\leq i\leq 2$, $1\leq j\leq k-1$,

\item $T_{1,j}\nsucceq T_{2,j}$ for $1\leq j\leq k$,

\item $T$ does not have any odd $II$-region.
\end{itemize}

\noindent The last condition is interpreted in \cite{nn2} in terms of tableaux :

The one row Example 5.10 of \cite{nn2} gives the result for $\chi_q(W_{k,a}^{(1)})$ : the extra condition is equivalent to $(T_j,T_{j+1})\notin\{(4,\overline{4}),(\overline{4},4)\}$ for $1\leq j\leq k-1$.

The two row Example 5.12 of \cite{nn2} (see also \cite{ss}) gives the result for $\chi_q(W_{k,a}^{(2)})$ : the extra condition is equivalent to conditions ($\theta$3'), ($\theta$4').

\subsection{Type $D_4^{(2)}$} Let us consider the twisted quantum affine algebra of type $D_4^{(2)}$.

For $a\in\CC^*$, let $b\in\CC^*$ such that $b^2 = 1$ and :
{\allowdisplaybreaks
\begin{equation*}
\begin{aligned}[c]
    & \ffbox{1}_{t,a} = Z_{1,a}, 
    \\
    & \ffbox{2}_{t,a} = Z_{1,aq^4}^{-1} Z_{2,aq^2},
    \\
    & \ffbox{3}_{t,a} = Z_{2,aq^6}^{-1} Z_{3,bq^2} Z_{3,-bq^2},
    \\
    & \ffbox{4}_{t,a} = Z_{3,bq^2}Z_{3,-bq^4}^{-1},
    \\
    & \ffbox{\overline{4}}_{t,a} = Z_{3,-bq^2}Z_{3,bq^4}^{-1},
    \\    
    & \ffbox{\overline{3}}_{t,a} = Z_{3,bq^4}^{-1}Z_{3,-bq^4}^{-1}Z_{2,aq^6},
    \\
    & \ffbox{\overline{2}}_{t,a} = Z_{1,aq^8}Z_{2,aq^{10}}^{-1},
    \\
    & \ffbox{\overline{1}}_{t,a} = Z_{1,aq^{12}}^{-1}.
\end{aligned}
\end{equation*}}

We define $\text{Tab}(1,k)$ and $\text{Tab}(2,k)$ as in the previous section.

For $T\in\text{Tab}(1,k)$ and $a\in\CC^*$ we set 
$$m_{t,T,a}^{(1)} = \prod_{1\leq j\leq k} \ffbox{T_j}_{t,aq^{2(j-1)}}.$$
For $T\in\text{Tab}(2,k)$ and $a\in\CC^*$ we set 
$$m_{t,T,a}^{(2)} = \prod_{1\leq i\leq 2, 1\leq j\leq k} \ffbox{T_{i,j}}_{t,aq^{2(j-i)}}.$$
As a consequence of Theorem \ref{transi} and Theorem \ref{expu} we get :
\begin{thm} For $a\in\CC^*$, $k\geq 1$ we have :
$$\chi_q^\sigma(W_{k,a}^{(1)}) = \sum_{T\in \text{Tab}(1,k)} m_{t,T,a}^{(1)}\text{ and }\chi_q^\sigma(W_{k,a}^{(2)}) = \sum_{T\in \text{Tab}(2,k)} m_{t,T,a}^{(2)}.$$ 
\end{thm}
For $a\in\CC^*$, let :
{\allowdisplaybreaks
\begin{equation*}
\begin{aligned}[c]
    & \ffbox{1}_{s,a} = Z_{3,a}, 
    \\
    & \ffbox{2}_{s,a} = Z_{3,aq^2}^{-1} Z_{2,a^2q^2},
    \\
    & \ffbox{3}_{s,a} = Z_{2,a^2q^6}^{-1} Z_{1,a^2q^4} Z_{3,-aq^2},
    \\
    & \ffbox{4}_{s,a} = Z_{1,a^2q^4}Z_{3,-aq^4}^{-1},
    \\
    & \ffbox{\overline{4}}_{s,a} = Z_{3,-aq^2}Z_{1,a^2q^8}^{-1},
    \\    
    & \ffbox{\overline{3}}_{s,a} = Z_{3,-aq^4}^{-1}Z_{1,a^2q^8}^{-1}Z_{2,a^2q^6},
    \\
    & \ffbox{\overline{2}}_{s,a} = Z_{3,aq^4}Z_{2,a^2q^{10}}^{-1},
    \\
    & \ffbox{\overline{1}}_{s,a} = Z_{3,aq^6}^{-1}.
\end{aligned}
\end{equation*}}
For $T\in\text{Tab}(1,k)$ and $a\in\CC^*$ we set 
$$m_{s,T,a}^{(3)} = \prod_{1\leq j\leq k} \ffbox{T_j}_{s,aq^{2(j-1)}}.$$
As a consequence of Theorem \ref{transi} and Theorem \ref{expu} we get :
\begin{thm} For $a\in\CC^*$, $k\geq 1$ we have :
$$\chi_q^\sigma(W_{k,a}^{(3)}) = \sum_{T\in \text{Tab}(1,k)} m_{t,T,a}^{(1)}.$$ 
\end{thm}

\begin{cor} 
$W_{k,a}^{(1)}$, $W_{k,a}^{(3)}$ are thin and $W_{k,a}^{(2)}$ has affine degree $2$.
\end{cor}

\section{Twisted $q$-characters of fundamental representations}\label{expfund}
 
In this section we get explicit formulas for the twisted $q$-character of fundamental representations of twisted quantum affine algebras (that it to say the first term in the inductive twisted $T$-system).

For type $A_{2n}^{(2)}$ ($n\geq 1$) and type $A_{2n-1}^{(2)}$ ($n\geq 3$) : see section \ref{krexp}.

\subsection{Type $D_{n+1}^{(2)}$ ($n\geq 2$)}

The $q$-characters of fundamental representations for untwisted quantum affine algebras of type $D$ are known (see the formulas in \cite{ks2}).

First suppose that $i_0\leq n-1$.

Let $\mathbf B = \{ 1, \dots, n + 1, \overline{n+1}, \dots, \overline{1}\}$. 
We give the ordering $\prec$ on the set $\mathbf B$ by
\begin{equation*}
  1 \prec 2 \prec \cdots \prec n \prec \begin{matrix}n+1\\\overline{n+1}\end{matrix}\prec \overline{n} \prec\cdots \prec \overline{2}\prec \overline{1}.
\end{equation*}

For $a\in\CC^*$ and $b\in\CC^*$ such that $b^2 = a$. Let
\begin{equation*}
\begin{split}
\ffbox{i}_a = \begin{cases} Z_{1,a} &\text{ if $i=1$,} 
    \\
    Z_{i-1,aq^{2i}}^{-1} Z_{i,aq^{2(i-1)}} &\text{ if $2 \leq i \leq n - 1$,}
    \\
    Z_{n - 1,aq^{2n}}^{-1} Z_{n,bq^{n-1}}Z_{n,-bq^{n-1}}&\text{ if $i=n$,}
    \\
    Z_{n,bq^{n+1}}^{-1} Z_{n,-bq^{n - 1}}&\text{ if $i=n+1$,}
    \\
    Z_{n,-bq^{n+1}}^{-1} Z_{n,bq^{n - 1}}&\text{ if $i=\overline{n+1}$,}
    \\
    Z_{n-1,aq^{2n}} Z_{n,bq^{n+1}}^{-1}Z_{n,-bq^{n+1}}^{-1}&\text{ if $i=\overline{n}$,}
    \\
    Z_{j-1,aq^{4n  - 2j}} Z_{j,aq^{4n + 2 - 2j}}^{-1}&\text{ if $i=\overline{j}$ and $2 \le j \le n - 1$,}
    \\ 
    Z_{1,aq^{4n}}^{-1}&\text{ if $i=\overline{1}$.}
\end{cases}
\end{split}
\end{equation*}

For $1\leq i_0\leq n-1$, let $\text{Tab}_D(i_0)$ be the set of tableaux $T = (T_{i})_{1\leq i\leq i_0}$ with coefficients in $\mathbf B$ such that for $1\leq i\leq i_0 - 1$, $T_{i}\nsucceq T_{i+1}$.

For such a tableaux $T\in\text{Tab}_D(i_0)$ and $a\in\CC^*$ we set 
$$m_{T,a} = \prod_{1\leq i\leq i_0} \ffbox{T_{i}}_{aq^{2-2i}}.$$

\begin{prop} For $1\leq i_0\leq n-1$ and $a\in\CC^*$ we have :
$$\chi_q^\sigma(V_{i_0,a}) = \sum_{T\in \text{Tab}'(i_0,k)} m_{T,aq^{i_0-1}}.$$
\end{prop}

Following \cite{Naexa, kks} we define for $a\in\CC^*$ the half size box as
\newcommand{\fhbox}[1]{
\setbox9=\hbox{$\scriptstyle\overline{n\!\!-\!\!1}$}
\framebox[\yh][c]{\rule{0mm}{\ht9}${\scriptstyle #1}$}
}
\begin{equation*}
\begin{split}
   \fhbox{i}_a & = 
   \begin{cases}
      Z_{1,a^2q^{-2}} & \text{if $i=1$},
   \\
     Z_{i-1,a^2q^{2i-2}}^{-1} Z_{i,a^2q^{2i-4}} 
       & \text{if $2\le i\le n-1$},
   \\
     Z_{n-1,a^2q^{2n-4}}^{-1} & \text{if $i=n$},
   \\
     Z_{n,aq^{n}} & \text{if $i=n+1$},
   \end{cases}
\\
   \fhbox{\overline{i}}_a &= 
   \begin{cases}
     1 & \text{if $1\le i\le n-1$},
   \\
   Z_{n,-aq^{n+2}}^{-1}Z_{n,aq^{n+2}}^{-1} 
     & \text{if $i=n$},
   \\
   Z_{n,-aq^{n}} & \text{if $i=n+1$}.
   \end{cases}
\end{split}
\end{equation*}

Let $\mathcal B_{\operatorname{sp}}$ be the set of tableaux $T = (i_1,\dots,i_{n+1})$
satisfying
\begin{itemize}
    \item $i_a\in \mathbf B, \; i_1 \prec i_2 \prec \dots \prec i_{n+1}$,
    \item $i$ and $\overline{i}$ do not appear simultaneously,
    \item if $i_p = n$, then $n+1-p$ is even,
    \item if $i_p = \overline{n+1}$, then $n-p$ is odd.
\end{itemize}

For such a tableaux $T\in\mathcal B_{\operatorname{sp}}$, we define the monomial $m_T$ by
\begin{equation*}
  m_T = \prod_{p=1}^{n+1} \fhbox{i_p}_{aq^{n+2-2p}}.
\end{equation*}
\begin{prop} For $a\in\CC^*$ we have
$$\chi_q^\sigma(V_{n,a}) = \sum_{T\in \mathcal B_{\operatorname{sp}}} m_{T,a}.$$
\end{prop}

\subsection{Type $D_4^{(3)}$}\label{fundhe} In the following we will use the notation $i_a^p = Z_{i,a}^p$ for $i\in I_{\sigma}$, $a\in\CC^*$, $p\in\ZZ$.

We give $\chi_q^\sigma (V_{1,a})$ and $\chi_q^\sigma (V_{2,a})$ by drawing a graph as in \cite{Fre} : there is an arrow $m_1\rightarrow_i m_2$ is $m_1$ and $m_2$ satisfy $m_1 = A_{i,b}m_2$ for one $b\in\CC^*$, and $\tau_i(m_1)$, $\tau_i(m_2)$ are in the $q$-character of a certain simple $\hat{\U}_i$-module. Note for untwisted type $D_4^{(1)}$ a similar graph is given \cite{Nex}.

The result corresponds to the conjectural formulas in \cite{r}.

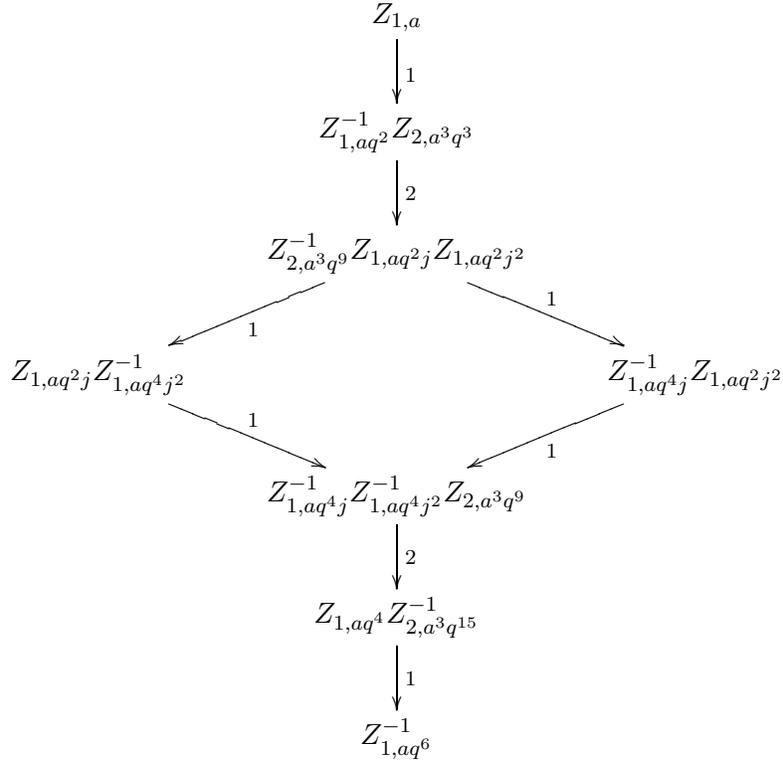
\begin{figure}[htbp]\label{fid}
\begin{center}
\begin{equation*}
\xymatrix{&Z_{1,a} \ar[d]^1&
\\&Z_{1,aq^2}^{-1}Z_{2,a^3q^3}\ar[d]^2&
\\&Z_{2,a^3q^9}^{-1}Z_{1,aq^2j}Z_{1,aq^2j^2}\ar[ld]^1\ar[rd]^1&
\\Z_{1,aq^2j}Z_{1,aq^4j^2}^{-1}\ar[rd]^1&&Z_{1,aq^4j}^{-1}Z_{1,aq^2j^2}\ar[ld]^1
\\&Z_{1,aq^4j}^{-1}Z_{1,aq^4j^2}^{-1}Z_{2,a^3q^9}\ar[d]^2&
\\&Z_{1,aq^4}Z_{2,a^3q^{15}}^{-1}\ar[d]^1&
\\&Z_{1,aq^6}^{-1}&}
\end{equation*}
\caption{Type $D_4^{(3)}$ : the graph of $\chi_q^\sigma(V_{1,a})$}
\end{center}
\end{figure}

\begin{figure}[htbp]
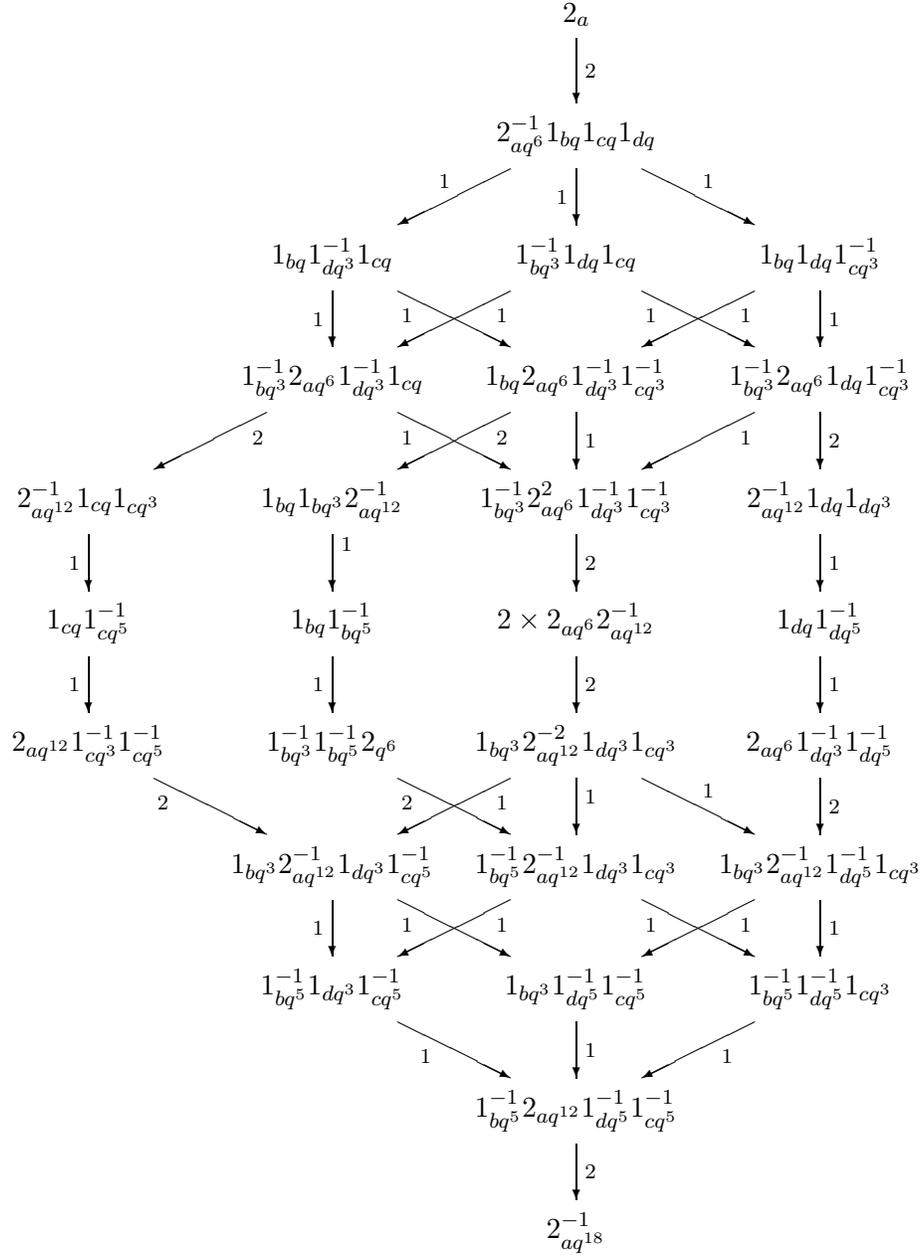
\label{fit}
\begin{center}
\leavevmode
\begin{equation*}
\divide\dgARROWLENGTH by 3
\dgARROWPARTS=6
\begin{diagram}
  \node[3]{2_a} \arrow{s,r}2
\\
  \node[3]{2_{aq^6}^{-1}1_{bq}1_{cq}1_{dq}}
  \arrow{sw,t,3}1 \arrow{s,l,3}1
  \arrow{se,t,3}1
\\
  \node[2]{1_{bq} 1_{dq^3}^{-1} 1_{cq}}
  \arrow{s,l}1 \arrow{se,b,1}1
  \node{1_{bq^3}^{-1} 1_{dq} 1_{cq}}
  \arrow{sw,b,1}1\arrow{se,b,1}1
  \node{1_{bq} 1_{dq} 1_{cq^3}^{-1}}
  \arrow{sw,b,1}1 \arrow{s,r}1
\\
  \node[2]{1_{bq^3}^{-1} 2_{aq^6} 1_{dq^3}^{-1} 1_{cq}}
  \arrow{sw,b,1}2 \arrow{se,b,1}1
  \node{1_{bq} 2_{aq^6} 1_{dq^3}^{-1} 1_{cq^3}^{-1}}
  \arrow{sw,b,1}2
  \arrow{s,r}1
  \node{1_{bq^3}^{-1} 2_{aq^6} 1_{dq} 1_{cq^3}^{-1}}
  \arrow{sw,b,1}1
  \arrow{s,r}2
\\
  \node{2_{aq^{12}}^{-1}1_{cq}1_{cq^3}}
  \arrow{s,l}1
  \node{1_{bq} 1_{bq^3} 2_{aq^{12}}^{-1}}
  \arrow{s,r,1}1
  \node{1_{bq^3}^{-1} 2_{aq^6}^2 1_{dq^3}^{-1}
  1_{cq^3}^{-1}}
  \arrow{s,r}2
  \node{2_{aq^{12}}^{-1} 1_{dq} 1_{dq^3}}
  \arrow{s,r}1
\\
  \node{1_{cq}1_{cq^5}^{-1}}
  \arrow{s,l}1
  \node{1_{bq}1_{bq^5}^{-1}}
  \arrow{s,l}1
  \node{2\times 2_{aq^6}2_{aq^{12}}^{-1}}
  \arrow{s,r}2
  \node{1_{dq}1_{dq^5}^{-1}}
  \arrow{s,r}1
\\
  \node{2_{aq^{12}}1_{cq^3}^{-1}1_{cq^5}^{-1}}
  \arrow{se,b,1}2
  \node{1_{bq^3}^{-1}1_{bq^5}^{-1} 2_{q^6}}
  \arrow{se,b,1}2
  \node{1_{bq^3}2_{aq^{12}}^{-2}1_{dq^3}1_{cq^3}}
  \arrow{sw,b,1}1\arrow{s,r,2}1
  \arrow{se,t,3}1
  \node{2_{aq^6}1_{dq^3}^{-1}1_{dq^5}^{-1}}
  \arrow{s,r}2
\\
  \node[2]{1_{bq^3}  2_{aq^{12}}^{-1} 1_{dq^3} 1_{cq^5}^{-1}}
  \arrow{s,l}1
  \arrow{se,b,1}1
  \node{1_{bq^5}^{-1} 2_{aq^{12}}^{-1} 1_{dq^3} 1_{cq^3}}
  \arrow{sw,b,1}1\arrow{se,b,1}1
  \node{1_{bq^3}  2_{aq^{12}}^{-1}1_{dq^5}^{-1}1_{cq^3}}
  \arrow{sw,b,1}1 \arrow{s,r}1
\\
  \node[2]{1_{bq^5}^{-1} 1_{dq^3} 1_{cq^5}^{-1}}
  \arrow{se,b,2}1
  \node{1_{bq^3} 1_{dq^5}^{-1} 1_{cq^5}^{-1}}
  \arrow{s,r}1
  \node{1_{bq^5}^{-1} 1_{dq^5}^{-1} 1_{cq^3}}
  \arrow{sw,b,2}1
\\
  \node[3]{1_{bq^5}^{-1}2_{aq^{12}}1_{dq^5}^{-1}1_{cq^5}^{-1}}
  \arrow{s,r}2
\\
  \node[3]{2_{aq^{18}}^{-1}}
\end{diagram}
\end{equation*}
\caption{Type $D_4^{(3)}$ : $\chi_q^\sigma(V_{2,a})$}
\end{center}
\end{figure}

Let $b\in\CC^*$ such that $b^3 = a$. Let $c = bj$ and $d = bj^2$. The graph of $\chi_q^\sigma(V_{1,a})$ is given in Figure \ref{fid} (see also the construction in \cite{kmoy}). The graph of $\chi_q^\sigma(V_{2,a})$ is given in Figure \ref{fit}.

\pagebreak

\subsection{Type $E_6^{(2)}$} Explicit formulas for fundamental representations of untwisted quantum affine algebras of type $E_6^{(1)}$ were obtained with a computer program by Nakajima, and by Hernandez-Schedler. For recent results in this direction see \cite{naco}. Nakajima's program is available online :

http://www.math.kyoto-u.ac.jp/\textasciitilde nakajima/Qchar/Qchar.html

\noindent The result for type $E_6^{(1)}$ from Hernandez-Schedler's program is available online :

http://www.math.uvsq.fr/\textasciitilde hernandez/e6.pdf 

\noindent With the results of this file and by using Theorem \ref{transi}, we get the twisted $q$-characters of fundamental representations in type $E_6^{(2)}$. The result corresponds to the conjectural formulas in \cite{r}.

\noindent In the present paper we list the result for the twisted $q$-character of $V_{1,a}$, $V_{2,a}$ and $V_{4,a}$. 
$V_{3,a}$ is of dimension $3732$ and $\chi_q^\sigma(V_{3,a})$ has $2925$ distinct monomials. So by lack of space, $\chi_q^\sigma(V_{3,a})$ is given online on the webpage :

http://www.math.uvsq.fr/\textasciitilde hernandez/e6twisted.pdf. 

\mk

\noindent Let us give $\chi_q^\sigma(V_{1,a})$, $\chi_q^\sigma(V_{2,a})$ and $\chi_q^\sigma(V_{4,a})$ :

\mk

$\chi_q^\sigma(V_{1,a})$ : there are 27 distinct monomials and the dimension of $V_{1,a}$ is 27 (all monomials have multiplicity $1$):

\mk

{\Small $1_{a}
+2_{aq}1^{-1}_{aq^2}
+3_{a^2q^4}2^{-1}_{aq^3}
+3^{-1}_{a^2q^8}4_{a^2q^6}2_{-aq^3}
+4^{-1}_{a^2q^{10}}2_{-aq^3}
+4_{a^2q^6}2^{-1}_{-aq^5}1_{-aq^4}
+3_{a^2q^8}4^{-1}_{a^2q^{10}}2^{-1}_{-aq^5}1_{-aq^4}
\\+4_{a^2q^6}1^{-1}_{-aq^6}
+3^{-1}_{a^2q^{12}}2_{aq^5}1_{-aq^4}
+3_{a^2q^9}4^{-1}_{a^2q^{10}}1^{-1}_{-aq^6}
+2^{-1}_{aq^7}1_{aq^6}1_{-aq^4}
+3^{-1}_{a^2q^{12}}2_{aq^5}2_{-aq^5}1^{-1}_{-aq^6}
+1^{-1}_{aq^8}1_{-aq^4}
\\+2^{-1}_{aq^7}2_{-aq^5}1_{aq^6}1^{-1}_{-aq^6}
+2_{aq^5}2^{-1}_{-aq^7}
+2_{-aq^5}1^{-1}_{aq^8}1^{-1}_{-aq^6}
+3_{a^2q^{12}}2^{-1}_{aq^7}2^{-1}_{-aq^7}1_{aq^6}
+3_{a^2q^{12}}2^{-1}_{-aq^7}1^{-1}_{aq^8}
\\+3^{-1}_{a^2q^{16}}4_{a^2q^{14}}1_{aq^6}
+3^{-1}_{a^2q^{16}}4_{a^2q^{14}}2_{aq^7}1^{-1}_{aq^8}
+4^{-1}_{a^2q^{18}}1_{aq^6}
+4^{-1}_{a^2q^{18}}2_{aq^7}1^{-1}_{aq^8}
+4_{a^2q^{14}}2^{-1}_{aq^9}
+3_{a^2q^{16}}4^{-1}_{a^2q^{18}}2^{-1}_{aq^9}
\\+3^{-1}_{a^2q^{20}}2_{-aq^9}
+2^{-1}_{-aq^{11}}1_{-aq^{10}}
+1^{-1}_{-aq^{12}}$.}

\mk

$\chi_q^\sigma(V_{4,a})$ : there are 78 distinct monomials and the dimension of $V_{4,a}$ is 79 (only one monomial has multiplicity $2$ : $3_{aq^{10}}3^{-1}_{aq^{14}}$) :

\mk

\noindent {\Small$
4_{a}
+3_{aq^2}4^{-1}_{aq^4}
+3^{-1}_{aq^6}2_{bq^2}2_{-bq^2}
+2^{-1}_{bq^4}2_{-bq^2}1_{bq^3}
+2_{bq^2}2^{-1}_{-bq^4}1_{-bq^3}
+3_{aq^6}2^{-1}_{bq^4}2^{-1}_{-bq^4}1_{bq^3}1_{-bq^3}
+2_{-bq^2}1^{-1}_{bq^5}
\\+2_{bq^2}1^{-1}_{-bq^5}
+3^{-1}_{aq^{10}}4_{aq^8}1_{bq^3}1_{-bq^3}
+3_{aq^6}2^{-1}_{-bq^4}1^{-1}_{bq^5}1_{-bq^3}
+3_{aq^6}2^{-1}_{bq^4}1_{bq^3}1^{-1}_{-bq^5}
+4^{-1}_{aq^{12}}1_{bq^3}1_{-bq^3}
\\+3^{-1}_{aq^{10}}4_{aq^8}2_{bq^4}1^{-1}_{bq^5}1_{-bq^3}
+3^{-1}_{aq^{10}}4_{aq^8}2_{-bq^4}1_{bq^3}1^{-1}_{-bq^5}
+3_{aq^6}1^{-1}_{bq^5}1^{-1}_{-bq^5}
+4^{-1}_{aq^{12}}2_{bq^4}1^{-1}_{bq^5}1_{-bq^3}
\\+4^{-1}_{aq^{12}}2_{-bq^4}1_{bq^3}1^{-1}_{-bq^5}
+4_{aq^8}2^{-1}_{bq^6}1_{-bq^3}
+3^{-1}_{aq^{10}}4_{aq^8}2_{bq^4}2_{-bq^4}1^{-1}_{bq^5}1^{-1}_{-bq^5}
+4_{aq^8}2^{-1}_{-bq^6}1_{bq^3}
\\+3_{aq^{10}}4^{-1}_{aq^{12}}2^{-1}_{bq^6}1_{-bq^3}
+4^{-1}_{aq^{12}}2_{bq^4}2_{-bq^4}1^{-1}_{bq^5}1^{-1}_{-bq^5}
+3_{aq^{10}}4^{-1}_{aq^{12}}2^{-1}_{-bq^6}1_{bq^3}
+4_{aq^8}2^{-1}_{bq^6}2_{-bq^4}1^{-1}_{-bq^5}
\\+4_{aq^8}2_{bq^4}2^{-1}_{-bq^6}1^{-1}_{bq^5}
+3^{-1}_{aq^{14}}2_{-bq^6}1_{-bq^3}
+3_{aq^{10}}4^{-1}_{aq^{12}}2^{-1}_{bq^6}2_{-bq^4}1^{-1}_{-bq^5}
+3_{aq^{10}}4^{-1}_{aq^{12}}2_{bq^4}2^{-1}_{-bq^6}1^{-1}_{bq^5}
\\+3^{-1}_{aq^{14}}2_{bq^6}1_{bq^3}
+3_{aq^{10}}4_{aq^8}2^{-1}_{bq^6}2^{-1}_{-bq^6}
+2^{-1}_{-bq^8}1_{-bq^3}1_{-bq^7}
+3^{-1}_{aq^{14}}2_{-bq^4}2_{-bq^6}1^{-1}_{-bq^5}
+3^2{aq^{10}}4^{-1}_{aq^{12}}2^{-1}_{bq^6}2^{-1}_{-bq^6}
\\+3^{-1}_{aq^{14}}2_{bq^4}2_{bq^6}1^{-1}_{bq^5}
+2^{-1}_{bq^8}1_{bq^3}1_{bq^7}
+3^{-1}_{aq^{14}}4_{aq^8}4_{aq^{12}}
+1_{-bq^3}1^{-1}_{-bq^9}
+2_{-bq^4}2^{-1}_{-bq^8}1^{-1}_{-bq^5}1_{-bq^7}
+23_{aq^{10}}3^{-1}_{aq^{14}}
\\+2_{bq^4}2^{-1}_{bq^8}1^{-1}_{bq^5}1_{bq^7}
+1_{bq^3}1^{-1}_{bq^9}
+4_{aq^8}4^{-1}_{aq^{16}}
+2_{-bq^4}1^{-1}_{-bq^5}1^{-1}_{-bq^9}
+3_{aq^{10}}2^{-1}_{-bq^6}2^{-1}_{-bq^8}1_{-bq^7}
+3^{-2}_{aq^{14}}4_{aq^{12}}2_{bq^6}2_{-bq^6}
\\+3_{aq^{10}}2^{-1}_{bq^6}2^{-1}_{bq^8}1_{bq^7}
+2_{bq^4}1^{-1}_{bq^5}1^{-1}_{bq^9}
+3_{aq^{10}}4^{-1}_{aq^{12}}4^{-1}_{aq^{16}}
+3_{aq^{10}}2^{-1}_{-bq^6}1^{-1}_{-bq^9}
+3^{-1}_{aq^{14}}4_{aq^{12}}2_{bq^6}2^{-1}_{-bq^8}1_{-bq^7}
\\+3^{-1}_{aq^{14}}4^{-1}_{aq^{16}}2_{bq^6}2_{-bq^6}
+3^{-1}_{aq^{14}}4_{aq^{12}}2^{-1}_{bq^8}2_{-bq^6}1_{bq^7}
+3_{aq^{10}}2^{-1}_{bq^6}1^{-1}_{bq^9}
+3^{-1}_{aq^{14}}4_{aq^{12}}2_{bq^6}1^{-1}_{-bq^9}
\\+4^{-1}_{aq^{16}}2_{bq^6}2^{-1}_{-bq^8}1_{-bq^7}
+4_{aq^{12}}2^{-1}_{bq^8}2^{-1}_{-bq^8}1_{bq^7}1_{-bq^7}
+4^{-1}_{aq^{16}}2^{-1}_{bq^8}2_{-bq^6}1_{bq^7}
+3^{-1}_{aq^{14}}4_{aq^{12}}2_{-bq^6}1^{-1}_{bq^9}
\\+4^{-1}_{aq^{16}}2_{bq^6}1^{-1}_{-bq^9}
+4_{aq^{12}}2^{-1}_{bq^8}1_{bq^7}1^{-1}_{-bq^9}
+3_{aq^{14}}4^{-1}_{aq^{16}}2^{-1}_{bq^8}2^{-1}_{-bq^8}1_{bq^7}1_{-bq^7}
+4_{aq^{12}}2^{-1}_{-bq^8}1^{-1}_{bq^9}1_{-bq^7}
\\+4^{-1}_{aq^{16}}2_{-bq^6}1^{-1}_{bq^9}
+3_{aq^{14}}4^{-1}_{aq^{16}}2^{-1}_{bq^8}1_{bq^7}1^{-1}_{-bq^9}
+4_{aq^{12}}1^{-1}_{bq^9}1^{-1}_{-bq^9}
+3^{-1}_{aq^{18}}1_{bq^7}1_{-bq^7}
+3_{aq^{14}}4^{-1}_{aq^{16}}2^{-1}_{-bq^8}1^{-1}_{bq^9}1_{-bq^7}
\\+3^{-1}_{aq^{18}}2_{-bq^8}1_{bq^7}1^{-1}_{-bq^9}
+3_{aq^{14}}4^{-1}_{aq^{16}}1^{-1}_{bq^9}1^{-1}_{-bq^9}
+3^{-1}_{aq^{18}}2_{bq^8}1^{-1}_{bq^9}1_{-bq^7}
+2^{-1}_{-bq^{10}}1_{bq^7}
+3^{-1}_{aq^{18}}2_{bq^8}2_{-bq^8}1^{-1}_{bq^9}1^{-1}_{-bq^9}
\\+2^{-1}_{bq^{10}}1_{-bq^7}
+2_{bq^8}2^{-1}_{-bq^{10}}1^{-1}_{bq^9}
+2^{-1}_{bq^{10}}2_{-bq^8}1^{-1}_{-bq^9}
+3_{aq^{18}}2^{-1}_{bq^{10}}2^{-1}_{-bq^{10}}
+3^{-1}_{aq^{22}}4_{aq^{20}}
+4^{-1}_{aq^{24}}$}

\mk

$\chi_q^\sigma(V_{2,a})$ : there are 351 distinct monomials and the dimension of $V_{2,a}$ is 378 :

\mk

{\Small $2_{a}
+
3_{a^2q^2}2^{-1}_{aq^2}1_{aq}
+
3^{-1}_{a^2q^6}4_{a^2q^4}2_{-aq^2}1_{aq}
+
3_{a^2q^2}1^{-1}_{aq^3}
+
4^{-1}_{a^2q^8}2_{-aq^2}1_{aq}
+
4_{a^2q^4}2^{-1}_{-aq^4}1_{aq}1_{-aq^3}
\\+
3^{-1}_{a^2q^6}4_{a^2q^4}2_{aq^2}2_{-aq^2}1^{-1}_{aq^3}
+
3_{a^2q^6}4^{-1}_{a^2q^8}2^{-1}_{-aq^4}1_{aq}1_{-aq^3}
+
4^{-1}_{a^2q^8}2_{aq^2}2_{-aq^2}1^{-1}_{aq^3}
+
4_{a^2q^4}2_{aq^2}2^{-1}_{-aq^4}1^{-1}_{aq^3}1_{-aq^3}
\\+
4_{a^2q^4}1_{aq}1^{-1}_{-aq^5}
+
4_{a^2q^4}2^{-1}_{aq^4}2_{-aq^2}
+
3^{-1}_{a^2q^{10}}2_{aq^4}1_{aq}1_{-aq^3}
+
3_{a^2q^6}4^{-1}_{a^2q^8}2_{aq^2}2^{-1}_{-aq^4}1^{-1}_{aq^3}1_{-aq^3}
\\+
3_{a^2q^6}4^{-1}_{a^2q^8}1_{aq}1^{-1}_{-aq^5}
+
3_{a^2q^6}4^{-1}_{a^2q^8}2^{-1}_{aq^4}2_{-aq^2}
+
3_{a^2q^6}4_{a^2q^4}2^{-1}_{aq^4}2^{-1}_{-aq^4}1_{-aq^3}
+
4_{a^2q^4}2_{aq^2}1^{-1}_{aq^3}1^{-1}_{-aq^5}
\\+
2^{-1}_{aq^6}1_{aq}1_{aq^5}1_{-aq^3}
+
3^{-1}_{a^2q^{10}}2_{aq^2}2_{aq^4}1^{-1}_{aq^3}1_{-aq^3}
+
3^{-1}_{a^2q^{10}}2_{aq^4}2_{-aq^4}1_{aq}1^{-1}_{-aq^5}
+
3^2_{a^2q^6}4^{-1}_{a^2q^8}2^{-1}_{aq^4}2^{-1}_{-aq^4}1_{-aq^3}
\\+
3_{a^2q^6}4^{-1}_{a^2q^8}2_{aq^2}1^{-1}_{aq^3}1^{-1}_{-aq^5}
+
3^{-1}_{a^2q^{10}}2_{-aq^2}2_{-aq^4}
+
3^{-1}_{a^2q^{10}}4_{a^2q^4}4_{a^2q^8}1_{-aq^3}
+
3_{a^2q^6}4_{a^2q^4}2^{-1}_{aq^4}1^{-1}_{-aq^5}
\\+
1_{aq}1^{-1}_{aq^7}1_{-aq^3}
+
2_{aq^2}2^{-1}_{aq^6}1^{-1}_{aq^3}1_{aq^5}1_{-aq^3}
+
2^{-1}_{aq^6}2_{-aq^4}1_{aq}1_{aq^5}1^{-1}_{-aq^5}
+2
3_{a^2q^6}3^{-1}_{a^2q^{10}}1_{-aq^3}
\\+
3^{-1}_{a^2q^{10}}2_{aq^2}2_{aq^4}2_{-aq^4}1^{-1}_{aq^3}1^{-1}_{-aq^5}
+
2_{aq^4}2^{-1}_{-aq^6}1_{aq}
+
3^2_{a^2q^6}4^{-1}_{a^2q^8}2^{-1}_{aq^4}1^{-1}_{-aq^5}
+
2_{-aq^2}2^{-1}_{-aq^6}1_{-aq^5}
\\+
4_{a^2q^4}4^{-1}_{a^2q^{12}}1_{-aq^3}
+
3^{-1}_{a^2q^{10}}4_{a^2q^4}4_{a^2q^8}2_{-aq^4}1^{-1}_{-aq^5}
+
2_{aq^2}1^{-1}_{aq^3}1^{-1}_{aq^7}1_{-aq^3}
+
3_{a^2q^6}2^{-1}_{aq^4}2^{-1}_{aq^6}1_{aq^5}1_{-aq^3}
\\+
3^{-2}_{a^2q^{10}}4_{a^2q^8}2_{aq^4}2_{-aq^4}1_{-aq^3}
+
3_{a^2q^6}2^{-1}_{-aq^4}2^{-1}_{-aq^6}1_{-aq^3}1_{-aq^5}
+
3_{a^2q^6}4^{-1}_{a^2q^8}4^{-1}_{a^2q^{12}}1_{-aq^3}
+
2_{-aq^4}1_{aq}1^{-1}_{aq^7}1^{-1}_{-aq^5}
\\+
2_{aq^2}2^{-1}_{aq^6}2_{-aq^4}1^{-1}_{aq^3}1_{aq^5}1^{-1}_{-aq^5}
+
3_{a^2q^{10}}2^{-1}_{aq^6}2^{-1}_{-aq^6}1_{aq}1_{aq^5}
+2
3_{a^2q^6}3^{-1}_{a^2q^{10}}2_{-aq^4}1^{-1}_{-aq^5}
+
2_{aq^2}2_{aq^4}2^{-1}_{-aq^6}1^{-1}_{aq^3}
\\+
2_{-aq^2}1^{-1}_{-aq^7}
+
4_{a^2q^4}4^{-1}_{a^2q^{12}}2_{-aq^4}1^{-1}_{-aq^5}
+
4_{a^2q^4}4_{a^2q^8}2^{-1}_{-aq^6}
+
2_{aq^2}2_{-aq^4}1^{-1}_{aq^3}1^{-1}_{aq^7}1^{-1}_{-aq^5}
\\+
3_{a^2q^6}2^{-1}_{aq^4}2^{-1}_{aq^6}2_{-aq^4}1_{aq^5}1^{-1}_{-aq^5}
+
3^{-2}_{a^2q^{10}}4_{a^2q^8}2_{aq^4}2^2_{-aq^4}1^{-1}_{-aq^5}
+
3_{a^2q^6}4^{-1}_{a^2q^8}4^{-1}_{a^2q^{12}}2_{-aq^4}1^{-1}_{-aq^5}
\\+
3_{a^2q^6}2^{-1}_{aq^4}1^{-1}_{aq^7}1_{-aq^3}
+
3^{-1}_{a^2q^{10}}4_{a^2q^8}2^{-1}_{aq^6}2_{-aq^4}1_{aq^5}1_{-aq^3}
+
3^{-1}_{a^2q^{10}}4^{-1}_{a^2q^{12}}2_{aq^4}2_{-aq^4}1_{-aq^3}
+2
3_{a^2q^6}2^{-1}_{-aq^6}
\\+
3^{-1}_{a^2q^{10}}4_{a^2q^8}2_{aq^4}2^{-1}_{-aq^6}1_{-aq^3}1_{-aq^5}
+
3_{a^2q^6}2^{-1}_{-aq^4}1_{-aq^3}1^{-1}_{-aq^7}
+
3_{a^2q^{10}}2^{-1}_{-aq^6}1_{aq}1^{-1}_{aq^7}
\\+
3_{a^2q^{10}}2_{aq^2}2^{-1}_{aq^6}2^{-1}_{-aq^6}1^{-1}_{aq^3}1_{aq^5}
+
3^{-1}_{a^2q^{14}}4_{a^2q^{12}}1_{aq}1_{aq^5}
+
3_{a^2q^{10}}4_{a^2q^4}4^{-1}_{a^2q^{12}}2^{-1}_{-aq^6}
+
3_{a^2q^6}1^{-1}_{-aq^5}1^{-1}_{-aq^7}
\\+
3_{a^2q^{10}}2_{aq^2}2^{-1}_{-aq^6}1^{-1}_{aq^3}1^{-1}_{aq^7}
+
3_{a^2q^6}3_{a^2q^{10}}2^{-1}_{aq^4}2^{-1}_{aq^6}2^{-1}_{-aq^6}1_{aq^5}
+
3_{a^2q^6}3_{a^2q^{10}}4^{-1}_{a^2q^8}4^{-1}_{a^2q^{12}}2^{-1}_{-aq^6}
\\+
3_{a^2q^6}2^{-1}_{aq^4}2_{-aq^4}1^{-1}_{aq^7}1^{-1}_{-aq^5}
+
3^{-1}_{a^2q^{10}}4_{a^2q^8}2^{-1}_{aq^6}2^2_{-aq^4}1_{aq^5}1^{-1}_{-aq^5}
+
3^{-1}_{a^2q^{10}}4^{-1}_{a^2q^{12}}2_{aq^4}2^2_{-aq^4}1^{-1}_{-aq^5}
\\+2
3^{-1}_{a^2q^{10}}4_{a^2q^8}2_{aq^4}2_{-aq^4}2^{-1}_{-aq^6}
+
3^{-1}_{a^2q^{10}}4_{a^2q^8}2_{-aq^4}1^{-1}_{aq^7}1_{-aq^3}
+
4^{-1}_{a^2q^{12}}2^{-1}_{aq^6}2_{-aq^4}1_{aq^5}1_{-aq^3}
\\+
4_{a^2q^8}2^{-1}_{aq^6}2^{-1}_{-aq^6}1_{aq^5}1_{-aq^3}1_{-aq^5}
+
4^{-1}_{a^2q^{12}}2_{aq^4}2^{-1}_{-aq^6}1_{-aq^3}1_{-aq^5}
+
3^{-1}_{a^2q^{10}}4_{a^2q^8}2_{aq^4}1_{-aq^3}1^{-1}_{-aq^7}
\\+
3^{-1}_{a^2q^{14}}4_{a^2q^{12}}2_{aq^6}1_{aq}1^{-1}_{aq^7}
+
3^{-1}_{a^2q^{14}}4_{a^2q^{12}}2_{aq^2}1^{-1}_{aq^3}1_{aq^5}
+
4^{-1}_{a^2q^{16}}1_{aq}1_{aq^5}
+
3^{-1}_{a^2q^{14}}4_{a^2q^4}2_{aq^6}
\\+
4_{a^2q^8}2_{aq^4}2^{-2}_{-aq^6}1_{-aq^5}
+
3^{-1}_{a^2q^{10}}4_{a^2q^8}2_{aq^4}2_{-aq^4}1^{-1}_{-aq^5}1^{-1}_{-aq^7}
+
3^{-1}_{a^2q^{14}}4_{a^2q^{12}}2_{aq^2}2_{aq^6}1^{-1}_{aq^3}1^{-1}_{aq^7}
\\+
3_{a^2q^6}3_{a^2q^{10}}2^{-1}_{aq^4}2^{-1}_{-aq^6}1^{-1}_{aq^7}
+
3_{a^2q^6}3^{-1}_{a^2q^{14}}4_{a^2q^{12}}2^{-1}_{aq^4}1_{aq^5}
+2
4_{a^2q^8}2^{-1}_{aq^6}2_{-aq^4}2^{-1}_{-aq^6}1_{aq^5}
\\+
3_{a^2q^6}3^{-1}_{a^2q^{14}}4^{-1}_{a^2q^8}2_{aq^6}
+2
4^{-1}_{a^2q^{12}}2_{aq^4}2_{-aq^4}2^{-1}_{-aq^6}
+
3^{-1}_{a^2q^{10}}4_{a^2q^8}2^2_{-aq^4}1^{-1}_{aq^7}1^{-1}_{-aq^5}
+
4^{-1}_{a^2q^{12}}2^{-1}_{aq^6}2^2_{-aq^4}1_{aq^5}1^{-1}_{-aq^5}
\\+
4^{-1}_{a^2q^{12}}2_{-aq^4}1^{-1}_{aq^7}1_{-aq^3}
+
4_{a^2q^8}2^{-1}_{-aq^6}1^{-1}_{aq^7}1_{-aq^3}1_{-aq^5}
+
3_{a^2q^{10}}4^{-1}_{a^2q^{12}}2^{-1}_{aq^6}2^{-1}_{-aq^6}1_{aq^5}1_{-aq^3}1_{-aq^5}
\\+
4_{a^2q^8}2^{-1}_{aq^6}1_{aq^5}1_{-aq^3}1^{-1}_{-aq^7}
+
4^{-1}_{a^2q^{12}}2_{aq^4}1_{-aq^3}1^{-1}_{-aq^7}
+
4^{-1}_{a^2q^{16}}2_{aq^6}1_{aq}1^{-1}_{aq^7}
+
4_{a^2q^{12}}2^{-1}_{aq^8}1_{aq}
\\+
4^{-1}_{a^2q^{16}}2_{aq^2}1^{-1}_{aq^3}1_{aq^5}
+
4_{a^2q^4}2^{-1}_{aq^8}1_{aq^7}
+
3^{-1}_{a^2q^{10}}3^{-1}_{a^2q^{14}}4_{a^2q^8}4_{a^2q^{12}}2_{-aq^4}1_{aq^5}
+
3^{-1}_{a^2q^{10}}3^{-1}_{a^2q^{14}}2_{aq^4}2_{aq^6}2_{-aq^4}
\\+
3_{a^2q^{10}}4_{a^2q^8}2^{-1}_{aq^6}2^{-2}_{-aq^6}1_{aq^5}1_{-aq^5}
+
3_{a^2q^{10}}4^{-1}_{a^2q^{12}}2_{aq^4}2^{-2}_{-aq^6}1_{-aq^5}
+
4_{a^2q^8}2^{-1}_{aq^6}2_{-aq^4}1_{aq^5}1^{-1}_{-aq^5}1^{-1}_{-aq^7}
\\+
4^{-1}_{a^2q^{12}}2_{aq^4}2_{-aq^4}1^{-1}_{-aq^5}1^{-1}_{-aq^7}
+
4^{-1}_{a^2q^{16}}2_{aq^2}2_{aq^6}1^{-1}_{aq^3}1^{-1}_{aq^7}
+
4_{a^2q^8}2_{aq^4}2^{-1}_{-aq^6}1^{-1}_{-aq^7}
+
4_{a^2q^{12}}2_{aq^2}2^{-1}_{aq^8}1^{-1}_{aq^3}
\\+
3_{a^2q^6}3^{-1}_{a^2q^{14}}4_{a^2q^{12}}2^{-1}_{aq^4}2_{aq^6}1^{-1}_{aq^7}
+2
4_{a^2q^8}2_{-aq^4}2^{-1}_{-aq^6}1^{-1}_{aq^7}
+
3_{a^2q^6}4^{-1}_{a^2q^{16}}2^{-1}_{aq^4}1_{aq^5}
\\+2
3_{a^2q^{10}}4^{-1}_{a^2q^{12}}2^{-1}_{aq^6}2_{-aq^4}2^{-1}_{-aq^6}1_{aq^5}
+
3_{a^2q^6}4^{-1}_{a^2q^8}2^{-1}_{aq^8}1_{aq^7}
+
4^{-1}_{a^2q^{12}}2^2_{-aq^4}1^{-1}_{aq^7}1^{-1}_{-aq^5}
\\+
3_{a^2q^{10}}4^{-1}_{a^2q^{12}}2^{-1}_{-aq^6}1^{-1}_{aq^7}1_{-aq^3}1_{-aq^5}
+
4_{a^2q^8}1^{-1}_{aq^7}1_{-aq^3}1^{-1}_{-aq^7}
+
3^{-1}_{a^2q^{14}}1_{aq^5}1_{-aq^3}1_{-aq^5}
\\+
3_{a^2q^{10}}4^{-1}_{a^2q^{12}}2^{-1}_{aq^6}1_{aq^5}1_{-aq^3}1^{-1}_{-aq^7}
+
3_{a^2q^{14}}4^{-1}_{a^2q^{16}}2^{-1}_{aq^8}1_{aq}
+
4_{a^2q^4}1^{-1}_{aq^9}
+
3_{a^2q^6}4_{a^2q^{12}}2^{-1}_{aq^4}2^{-1}_{aq^8}
\\+
3^{-1}_{a^2q^{10}}3^{-1}_{a^2q^{14}}4_{a^2q^8}4_{a^2q^{12}}2_{aq^6}2_{-aq^4}1^{-1}_{aq^7}
+
3_{a^2q^{10}}4_{a^2q^8}2^{-2}_{-aq^6}1^{-1}_{aq^7}1_{-aq^5}
+
3^2_{a^2q^{10}}4^{-1}_{a^2q^{12}}2^{-1}_{aq^6}2^{-2}_{-aq^6}1_{aq^5}1_{-aq^5}
\\+
4_{a^2q^8}2_{-aq^4}1^{-1}_{aq^7}1^{-1}_{-aq^5}1^{-1}_{-aq^7}
+
3_{a^2q^{10}}4^{-1}_{a^2q^{12}}2^{-1}_{aq^6}2_{-aq^4}1_{aq^5}1^{-1}_{-aq^5}1^{-1}_{-aq^7}
+
3^{-1}_{a^2q^{10}}4_{a^2q^8}4^{-1}_{a^2q^{16}}2_{-aq^4}1_{aq^5}
\\+2
3^{-1}_{a^2q^{14}}2_{-aq^4}1_{aq^5}
+
3^{-1}_{a^2q^{14}}4_{a^2q^8}4_{a^2q^{12}}2^{-1}_{-aq^6}1_{aq^5}1_{-aq^5}
+
3^{-1}_{a^2q^{10}}2_{aq^4}2^{-1}_{aq^8}2_{-aq^4}1_{aq^7}
\\+
3^{-1}_{a^2q^{14}}2_{aq^4}2_{aq^6}2^{-1}_{-aq^6}1_{-aq^5}
+
3_{a^2q^{10}}4_{a^2q^8}2^{-1}_{aq^6}2^{-1}_{-aq^6}1_{aq^5}1^{-1}_{-aq^7}
+
3_{a^2q^{10}}4^{-1}_{a^2q^{12}}2_{aq^4}2^{-1}_{-aq^6}1^{-1}_{-aq^7}
\\+
3_{a^2q^{14}}4^{-1}_{a^2q^{16}}2_{aq^2}2^{-1}_{aq^8}1^{-1}_{aq^3}
+
3_{a^2q^6}4^{-1}_{a^2q^{16}}2^{-1}_{aq^4}2_{aq^6}1^{-1}_{aq^7}
+2
3_{a^2q^{10}}4^{-1}_{a^2q^{12}}2_{-aq^4}2^{-1}_{-aq^6}1^{-1}_{aq^7}
+
3_{a^2q^6}4^{-1}_{a^2q^8}1^{-1}_{aq^9}
\\+
3^{-1}_{a^2q^{14}}2_{aq^6}1^{-1}_{aq^7}1_{-aq^3}1_{-aq^5}
+
3_{a^2q^{10}}4^{-1}_{a^2q^{12}}1^{-1}_{aq^7}1_{-aq^3}1^{-1}_{-aq^7}
+
3^{-1}_{a^2q^{14}}2_{-aq^6}1_{aq^5}1_{-aq^3}1^{-1}_{-aq^7}
+
3^{-1}_{a^2q^{18}}2_{-aq^8}1_{aq}
\\+
4^{-1}_{a^2q^{12}}4^{-1}_{a^2q^{16}}2_{-aq^4}1_{aq^5}
+
2^{-1}_{aq^6}2^{-1}_{aq^8}2_{-aq^4}1_{aq^5}1_{aq^7}
+
3_{a^2q^6}3_{a^2q^{14}}4^{-1}_{a^2q^{16}}2^{-1}_{aq^4}2^{-1}_{aq^8}
\\+
3^2_{a^2q^{10}}4^{-1}_{a^2q^{12}}2^{-2}_{-aq^6}1^{-1}_{aq^7}1_{-aq^5}
+
3_{a^2q^{10}}4^{-1}_{a^2q^{12}}2_{-aq^4}1^{-1}_{aq^7}1^{-1}_{-aq^5}1^{-1}_{-aq^7}
+
3^{-1}_{a^2q^{14}}2_{-aq^4}2_{-aq^6}1_{aq^5}1^{-1}_{-aq^5}1^{-1}_{-aq^7}
\\+
3^{-1}_{a^2q^{10}}4_{a^2q^8}4_{a^2q^{12}}2^{-1}_{aq^8}2_{-aq^4}
+
3^{-1}_{a^2q^{10}}4_{a^2q^8}4^{-1}_{a^2q^{16}}2_{aq^6}2_{-aq^4}1^{-1}_{aq^7}
+2
3^{-1}_{a^2q^{14}}2_{aq^6}2_{-aq^4}1^{-1}_{aq^7}
\\+
3^{-1}_{a^2q^{14}}4_{a^2q^8}4_{a^2q^{12}}2_{aq^6}2^{-1}_{-aq^6}1^{-1}_{aq^7}1_{-aq^5}
+
3_{a^2q^{10}}4_{a^2q^8}2^{-1}_{-aq^6}1^{-1}_{aq^7}1^{-1}_{-aq^7}
+2
3_{a^2q^{10}}3^{-1}_{a^2q^{14}}2^{-1}_{-aq^6}1_{aq^5}1_{-aq^5}
\\+
3^2_{a^2q^{10}}4^{-1}_{a^2q^{12}}2^{-1}_{aq^6}2^{-1}_{-aq^6}1_{aq^5}1^{-1}_{-aq^7}
+
4_{a^2q^8}4^{-1}_{a^2q^{16}}2^{-1}_{-aq^6}1_{aq^5}1_{-aq^5}
+
3^{-1}_{a^2q^{14}}4_{a^2q^8}4_{a^2q^{12}}1_{aq^5}1^{-1}_{-aq^7}
\\+
2_{aq^4}2^{-1}_{aq^8}2^{-1}_{-aq^6}1_{aq^7}1_{-aq^5}
+
3^{-1}_{a^2q^{10}}2_{aq^4}2_{-aq^4}1^{-1}_{aq^9}
+
3^{-1}_{a^2q^{14}}2_{aq^4}2_{aq^6}1^{-1}_{-aq^7}
+
3^{-1}_{a^2q^{18}}2_{aq^2}2_{-aq^8}1^{-1}_{aq^3}
\\+
2^{-1}_{aq^8}1_{-aq^3}1_{-aq^5}
+
3^{-1}_{a^2q^{14}}2_{aq^6}2_{-aq^6}1^{-1}_{aq^7}1_{-aq^3}1^{-1}_{-aq^7}
+
2^{-1}_{-aq^8}1_{aq^5}1_{-aq^3}
+
2^{-1}_{-aq^{10}}1_{aq}1_{-aq^9}
\\+
4^{-1}_{a^2q^{12}}4^{-1}_{a^2q^{16}}2_{aq^6}2_{-aq^4}1^{-1}_{aq^7}
+
3^{-2}_{a^2q^{14}}4_{a^2q^{12}}2_{aq^6}1_{aq^5}1_{-aq^5}
+
3_{a^2q^{10}}4^{-1}_{a^2q^{12}}4^{-1}_{a^2q^{16}}2^{-1}_{-aq^6}1_{aq^5}1_{-aq^5}
\\+
3_{a^2q^{10}}2^{-1}_{aq^6}2^{-1}_{aq^8}2^{-1}_{-aq^6}1_{aq^5}1_{aq^7}1_{-aq^5}
+
3^{-1}_{a^2q^{14}}2_{aq^6}2_{-aq^4}2_{-aq^6}1^{-1}_{aq^7}1^{-1}_{-aq^5}1^{-1}_{-aq^7}
+
2^{-1}_{aq^6}2_{-aq^4}1_{aq^5}1^{-1}_{aq^9}
\\+2
2^{-1}_{aq^8}2_{-aq^4}
+
3_{a^2q^6}3^{-1}_{a^2q^{18}}2^{-1}_{aq^4}2_{-aq^8}
+
3^{-1}_{a^2q^{10}}3_{a^2q^{14}}4_{a^2q^8}4^{-1}_{a^2q^{16}}2^{-1}_{aq^8}2_{-aq^4}
\\+2
3_{a^2q^{10}}3^{-1}_{a^2q^{14}}2_{aq^6}2^{-1}_{-aq^6}1^{-1}_{aq^7}1_{-aq^5}
+
3^2_{a^2q^{10}}4^{-1}_{a^2q^{12}}2^{-1}_{-aq^6}1^{-1}_{aq^7}1^{-1}_{-aq^7}
+
2_{-aq^4}2^{-1}_{-aq^8}1_{aq^5}1^{-1}_{-aq^5}
\\+2
3_{a^2q^{10}}3^{-1}_{a^2q^{14}}1_{aq^5}1^{-1}_{-aq^7}
+
4_{a^2q^8}4_{a^2q^{12}}2^{-1}_{aq^8}2^{-1}_{-aq^6}1_{-aq^5}
+
4_{a^2q^8}4^{-1}_{a^2q^{16}}2_{aq^6}2^{-1}_{-aq^6}1^{-1}_{aq^7}1_{-aq^5}
\\+
3^{-1}_{a^2q^{14}}4_{a^2q^8}4_{a^2q^{12}}2_{aq^6}1^{-1}_{aq^7}1^{-1}_{-aq^7}
+
4_{a^2q^8}4^{-1}_{a^2q^{16}}1_{aq^5}1^{-1}_{-aq^7}
+
2_{aq^4}2^{-1}_{-aq^6}1^{-1}_{aq^9}1_{-aq^5}
+
2_{aq^4}2^{-1}_{aq^8}1_{aq^7}1^{-1}_{-aq^7}
\\+
2_{aq^2}2^{-1}_{-aq^{10}}1^{-1}_{aq^3}1_{-aq^9}
+
2^{-1}_{aq^8}2_{-aq^6}1_{-aq^3}1^{-1}_{-aq^7}
+
2_{aq^6}2^{-1}_{-aq^8}1^{-1}_{aq^7}1_{-aq^3}
+
1_{aq}1^{-1}_{-aq^{11}}
+
2_{-aq^4}1^{-1}_{aq^7}1^{-1}_{aq^9}
\\+
3^{-1}_{a^2q^{10}}3^{-1}_{a^2q^{18}}4_{a^2q^8}2_{-aq^4}2_{-aq^8}
+
3^{-2}_{a^2q^{14}}4_{a^2q^{12}}2^2_{aq^6}1^{-1}_{aq^7}1_{-aq^5}
+
3_{a^2q^{10}}2^{-1}_{-aq^6}2^{-1}_{-aq^8}1_{aq^5}
\\+
3_{a^2q^{14}}4^{-1}_{a^2q^{12}}4^{-1}_{a^2q^{16}}2^{-1}_{aq^8}2_{-aq^4}
+
3_{a^2q^{10}}4^{-1}_{a^2q^{12}}4^{-1}_{a^2q^{16}}2_{aq^6}2^{-1}_{-aq^6}1^{-1}_{aq^7}1_{-aq^5}
+
3^{-2}_{a^2q^{14}}4_{a^2q^{12}}2_{aq^6}2_{-aq^6}1_{aq^5}1^{-1}_{-aq^7}
\\+
3_{a^2q^{10}}4^{-1}_{a^2q^{12}}4^{-1}_{a^2q^{16}}1_{aq^5}1^{-1}_{-aq^7}
+
3_{a^2q^{10}}2^{-1}_{aq^6}2^{-1}_{aq^8}1_{aq^5}1_{aq^7}1^{-1}_{-aq^7}
+
2^{-1}_{aq^8}2_{-aq^4}2_{-aq^6}1^{-1}_{-aq^5}1^{-1}_{-aq^7}
\\+
3^{-1}_{a^2q^{14}}4^{-1}_{a^2q^{16}}2_{aq^6}1_{aq^5}1_{-aq^5}
+
3^{-1}_{a^2q^{14}}4_{a^2q^{12}}2^{-1}_{aq^8}1_{aq^5}1_{aq^7}1_{-aq^5}
+
3_{a^2q^{10}}2^{-1}_{aq^6}2^{-1}_{-aq^6}1_{aq^5}1^{-1}_{aq^9}1_{-aq^5}
\\+2
3_{a^2q^{10}}2^{-1}_{aq^8}2^{-1}_{-aq^6}1_{-aq^5}
+
2_{aq^6}2_{-aq^4}2^{-1}_{-aq^8}1^{-1}_{aq^7}1^{-1}_{-aq^5}
+2
3_{a^2q^{10}}3^{-1}_{a^2q^{14}}2_{aq^6}1^{-1}_{aq^7}1^{-1}_{-aq^7}
\\+
3_{a^2q^6}2^{-1}_{aq^4}2^{-1}_{-aq^{10}}1_{-aq^9}
+
3_{a^2q^{14}}4_{a^2q^8}4^{-1}_{a^2q^{16}}2^{-1}_{aq^8}2^{-1}_{-aq^6}1_{-aq^5}
+
4_{a^2q^8}4_{a^2q^{12}}2^{-1}_{aq^8}1^{-1}_{-aq^7}
\\+
4_{a^2q^8}4^{-1}_{a^2q^{16}}2_{aq^6}1^{-1}_{aq^7}1^{-1}_{-aq^7}
+
2_{aq^4}1^{-1}_{aq^9}1^{-1}_{-aq^7}
+
2_{aq^2}1^{-1}_{aq^3}1^{-1}_{-aq^{11}}
+
3_{a^2q^{14}}2^{-1}_{aq^8}2^{-1}_{-aq^8}1_{-aq^3}
\\+
3_{a^2q^{10}}2^{-1}_{-aq^6}1^{-1}_{aq^7}1^{-1}_{aq^9}1_{-aq^5}
+
3_{a^2q^{10}}2_{aq^6}2^{-1}_{-aq^6}2^{-1}_{-aq^8}1^{-1}_{aq^7}
+
3^{-2}_{a^2q^{14}}4_{a^2q^{12}}2^2_{aq^6}2_{-aq^6}1^{-1}_{aq^7}1^{-1}_{-aq^7}
\\+
3_{a^2q^{10}}3_{a^2q^{14}}4^{-1}_{a^2q^{12}}4^{-1}_{a^2q^{16}}2^{-1}_{aq^8}2^{-1}_{-aq^6}1_{-aq^5}
+
3_{a^2q^{10}}4^{-1}_{a^2q^{12}}4^{-1}_{a^2q^{16}}2_{aq^6}1^{-1}_{aq^7}1^{-1}_{-aq^7}
+
3^{-1}_{a^2q^{18}}4^{-1}_{a^2q^{12}}2_{-aq^4}2_{-aq^8}
\\+
3^{-1}_{a^2q^{10}}4_{a^2q^8}2_{-aq^4}2^{-1}_{-aq^{10}}1_{-aq^9}
+
3^{-1}_{a^2q^{18}}4_{a^2q^8}2^{-1}_{-aq^6}2_{-aq^8}1_{-aq^5}
+
3^{-1}_{a^2q^{14}}4^{-1}_{a^2q^{16}}2^2_{aq^6}1^{-1}_{aq^7}1_{-aq^5}
\\+2
3^{-1}_{a^2q^{14}}4_{a^2q^{12}}2_{aq^6}2^{-1}_{aq^8}1_{-aq^5}
+
3^{-1}_{a^2q^{14}}4_{a^2q^{12}}2_{aq^6}2^{-1}_{-aq^8}1_{aq^5}
+
3^{-1}_{a^2q^{14}}4^{-1}_{a^2q^{16}}2_{aq^6}2_{-aq^6}1_{aq^5}1^{-1}_{-aq^7}
\\+
3^{-1}_{a^2q^{14}}4_{a^2q^{12}}2^{-1}_{aq^8}2_{-aq^6}1_{aq^5}1_{aq^7}1^{-1}_{-aq^7}
+
3_{a^2q^{10}}2^{-1}_{aq^6}1_{aq^5}1^{-1}_{aq^9}1^{-1}_{-aq^7}
+2
3_{a^2q^{10}}2^{-1}_{aq^8}1^{-1}_{-aq^7}
\\+
3_{a^2q^{14}}2^{-1}_{aq^8}2_{-aq^4}2^{-1}_{-aq^8}1^{-1}_{-aq^5}
+
4^{-1}_{a^2q^{16}}2^{-1}_{aq^8}1_{aq^5}1_{aq^7}1_{-aq^5}
+
3^{-1}_{a^2q^{14}}4_{a^2q^{12}}1_{aq^5}1^{-1}_{aq^9}1_{-aq^5}
+
3_{a^2q^6}2^{-1}_{aq^4}1^{-1}_{-aq^{11}}
\\+
3_{a^2q^{14}}4_{a^2q^8}4^{-1}_{a^2q^{16}}2^{-1}_{aq^8}1^{-1}_{-aq^7}
+
3^{-1}_{a^2q^{18}}4_{a^2q^{16}}1_{-aq^3}
+
4_{a^2q^8}2^{-1}_{-aq^6}2^{-1}_{-aq^{10}}1_{-aq^5}1_{-aq^9}
+
4_{a^2q^{12}}2^{-2}_{aq^8}1_{aq^7}1_{-aq^5}
\\+
3_{a^2q^{10}}1^{-1}_{aq^7}1^{-1}_{aq^9}1^{-1}_{-aq^7}
+
3_{a^2q^{10}}3_{a^2q^{14}}2^{-1}_{aq^8}2^{-1}_{-aq^6}2^{-1}_{-aq^8}
+
3^{-1}_{a^2q^{14}}4_{a^2q^{12}}2_{aq^6}1^{-1}_{aq^7}1^{-1}_{aq^9}1_{-aq^5}
\\+
3_{a^2q^{10}}3_{a^2q^{14}}4^{-1}_{a^2q^{12}}4^{-1}_{a^2q^{16}}2^{-1}_{aq^8}1^{-1}_{-aq^7}
+
3^{-1}_{a^2q^{14}}4_{a^2q^{12}}2^2_{aq^6}2^{-1}_{-aq^8}1^{-1}_{aq^7}
+
3^{-1}_{a^2q^{14}}4^{-1}_{a^2q^{16}}2^2_{aq^6}2_{-aq^6}1^{-1}_{aq^7}1^{-1}_{-aq^7}
\\+2
3^{-1}_{a^2q^{14}}4_{a^2q^{12}}2_{aq^6}2^{-1}_{aq^8}2_{-aq^6}1^{-1}_{-aq^7}
+
3_{a^2q^{10}}3^{-1}_{a^2q^{18}}4^{-1}_{a^2q^{12}}2^{-1}_{-aq^6}2_{-aq^8}1_{-aq^5}
+2
4^{-1}_{a^2q^{16}}2_{aq^6}2^{-1}_{aq^8}1_{-aq^5}
\\+
4^{-1}_{a^2q^{12}}2_{-aq^4}2^{-1}_{-aq^{10}}1_{-aq^9}
+
3^{-1}_{a^2q^{10}}4_{a^2q^8}2_{-aq^4}1^{-1}_{-aq^{11}}
+
3^{-1}_{a^2q^{18}}4_{a^2q^8}2_{-aq^8}1^{-1}_{-aq^7}
+
4^{-1}_{a^2q^{16}}2_{aq^6}2^{-1}_{-aq^8}1_{aq^5}
\\+
4_{a^2q^{12}}2^{-1}_{aq^8}2^{-1}_{-aq^8}1_{aq^5}1_{aq^7}
+
4^{-1}_{a^2q^{16}}2^{-1}_{aq^8}2_{-aq^6}1_{aq^5}1_{aq^7}1^{-1}_{-aq^7}
+
3^{-1}_{a^2q^{14}}4_{a^2q^{12}}2_{-aq^6}1_{aq^5}1^{-1}_{aq^9}1^{-1}_{-aq^7}
\\+
3^{-1}_{a^2q^{18}}4_{a^2q^{16}}2_{-aq^4}1^{-1}_{-aq^5}
+
4^{-1}_{a^2q^{16}}1_{aq^5}1^{-1}_{aq^9}1_{-aq^5}
+
4^{-1}_{a^2q^{20}}1_{-aq^3}
+
4_{a^2q^{12}}2^{-2}_{aq^8}2_{-aq^6}1_{aq^7}1^{-1}_{-aq^7}
\\+
3^{-1}_{a^2q^{14}}3^{-1}_{a^2q^{18}}2_{aq^6}2_{-aq^8}1_{-aq^5}
+
3_{a^2q^{10}}4^{-1}_{a^2q^{12}}2^{-1}_{-aq^6}2^{-1}_{-aq^{10}}1_{-aq^5}1_{-aq^9}
+
3_{a^2q^{14}}4^{-1}_{a^2q^{16}}2^{-2}_{aq^8}1_{aq^7}1_{-aq^5}
\\+
3^{-1}_{a^2q^{14}}4_{a^2q^{12}}2_{aq^6}2_{-aq^6}1^{-1}_{aq^7}1^{-1}_{aq^9}1^{-1}_{-aq^7}
+
4^{-1}_{a^2q^{16}}2_{aq^6}1^{-1}_{aq^7}1^{-1}_{aq^9}1_{-aq^5}
+
4_{a^2q^8}2^{-1}_{-aq^6}1_{-aq^5}1^{-1}_{-aq^{11}}
\\+
4_{a^2q^8}2^{-1}_{-aq^{10}}1^{-1}_{-aq^7}1_{-aq^9}
+
4_{a^2q^{12}}2^{-1}_{aq^8}1^{-1}_{aq^9}1_{-aq^5}
+
3_{a^2q^{10}}3^{-1}_{a^2q^{18}}4_{a^2q^{16}}2^{-1}_{-aq^6}
+2
4_{a^2q^{12}}2_{aq^6}2^{-1}_{aq^8}2^{-1}_{-aq^8}
\\+
3_{a^2q^{10}}3^{-1}_{a^2q^{18}}4^{-1}_{a^2q^{12}}2_{-aq^8}1^{-1}_{-aq^7}
+2
4^{-1}_{a^2q^{16}}2_{aq^6}2^{-1}_{aq^8}2_{-aq^6}1^{-1}_{-aq^7}
+
4^{-1}_{a^2q^{16}}2^2_{aq^6}2^{-1}_{-aq^8}1^{-1}_{aq^7}
+
4^{-1}_{a^2q^{12}}2_{-aq^4}1^{-1}_{-aq^{11}}
\\+
3_{a^2q^{14}}4^{-1}_{a^2q^{16}}2^{-1}_{aq^8}2^{-1}_{-aq^8}1_{aq^5}1_{aq^7}
+
4_{a^2q^{12}}2^{-1}_{-aq^8}1_{aq^5}1^{-1}_{aq^9}
+
4^{-1}_{a^2q^{16}}2_{-aq^6}1_{aq^5}1^{-1}_{aq^9}1^{-1}_{-aq^7}
+
4^{-1}_{a^2q^{20}}2_{-aq^4}1^{-1}_{-aq^5}
\\+
4_{a^2q^8}1^{-1}_{-aq^7}1^{-1}_{-aq^{11}}
+
3^{-1}_{a^2q^{14}}3^{-1}_{a^2q^{18}}4_{a^2q^{12}}4_{a^2q^{16}}2_{aq^6}
+
3^{-1}_{a^2q^{14}}3^{-1}_{a^2q^{18}}2_{aq^6}2_{-aq^6}2_{-aq^8}1^{-1}_{-aq^7}
\\+
3_{a^2q^{14}}4_{a^2q^{12}}2^{-2}_{aq^8}2^{-1}_{-aq^8}1_{aq^7}
+
3_{a^2q^{14}}4^{-1}_{a^2q^{16}}2^{-2}_{aq^8}2_{-aq^6}1_{aq^7}1^{-1}_{-aq^7}
+
4_{a^2q^{12}}2_{aq^6}2^{-1}_{-aq^8}1^{-1}_{aq^7}1^{-1}_{aq^9}
\\+
4^{-1}_{a^2q^{16}}2_{aq^6}2_{-aq^6}1^{-1}_{aq^7}1^{-1}_{aq^9}1^{-1}_{-aq^7}
+
4_{a^2q^{12}}2^{-1}_{aq^8}2_{-aq^6}1^{-1}_{aq^9}1^{-1}_{-aq^7}
+
3^{-1}_{a^2q^{18}}2^{-1}_{aq^8}2_{-aq^8}1_{aq^7}1_{-aq^5}
\\+
3^{-1}_{a^2q^{14}}2_{aq^6}2^{-1}_{-aq^{10}}1_{-aq^5}1_{-aq^9}
+
3_{a^2q^{10}}4^{-1}_{a^2q^{12}}2^{-1}_{-aq^6}1_{-aq^5}1^{-1}_{-aq^{11}}
+
3_{a^2q^{10}}4^{-1}_{a^2q^{12}}2^{-1}_{-aq^{10}}1^{-1}_{-aq^7}1_{-aq^9}
\\+
3_{a^2q^{14}}4^{-1}_{a^2q^{16}}2^{-1}_{aq^8}1^{-1}_{aq^9}1_{-aq^5}
+
3_{a^2q^{10}}4^{-1}_{a^2q^{20}}2^{-1}_{-aq^6}
+2
3_{a^2q^{14}}4^{-1}_{a^2q^{16}}2_{aq^6}2^{-1}_{aq^8}2^{-1}_{-aq^8}
+
3^{-1}_{a^2q^{18}}1_{aq^5}1_{aq^7}
\\+
3_{a^2q^{14}}4^{-1}_{a^2q^{16}}2^{-1}_{-aq^8}1_{aq^5}1^{-1}_{aq^9}
+
3_{a^2q^{10}}4^{-1}_{a^2q^{12}}1^{-1}_{-aq^7}1^{-1}_{-aq^{11}}
+
3^2_{a^2q^{14}}4^{-1}_{a^2q^{16}}2^{-2}_{aq^8}2^{-1}_{-aq^8}1_{aq^7}
\\+
3_{a^2q^{14}}4^{-1}_{a^2q^{16}}2_{aq^6}2^{-1}_{-aq^8}1^{-1}_{aq^7}1^{-1}_{aq^9}
+
3^{-1}_{a^2q^{14}}4_{a^2q^{12}}4^{-1}_{a^2q^{20}}2_{aq^6}
+2
3^{-1}_{a^2q^{18}}2_{aq^6}
+
3^{-1}_{a^2q^{18}}4_{a^2q^{12}}4_{a^2q^{16}}2^{-1}_{aq^8}1_{aq^7}
\\+
3^{-1}_{a^2q^{18}}2^{-1}_{aq^8}2_{-aq^6}2_{-aq^8}1_{aq^7}1^{-1}_{-aq^7}
+
3^{-1}_{a^2q^{14}}2_{aq^6}2_{-aq^6}2^{-1}_{-aq^{10}}1^{-1}_{-aq^7}1_{-aq^9}
+
3_{a^2q^{14}}4_{a^2q^{12}}2^{-1}_{aq^8}2^{-1}_{-aq^8}1^{-1}_{aq^9}
\\+
3_{a^2q^{14}}4^{-1}_{a^2q^{16}}2^{-1}_{aq^8}2_{-aq^6}1^{-1}_{aq^9}1^{-1}_{-aq^7}
+
2^{-1}_{aq^8}2^{-1}_{-aq^{10}}1_{aq^7}1_{-aq^5}1_{-aq^9}
+
3^{-1}_{a^2q^{18}}2_{-aq^8}1^{-1}_{aq^9}1_{-aq^5}
\\+
3^{-1}_{a^2q^{14}}2_{aq^6}1_{-aq^5}1^{-1}_{-aq^{11}}
+
3^{-1}_{a^2q^{18}}2_{aq^8}1_{aq^5}1^{-1}_{aq^9}
+
4^{-1}_{a^2q^{16}}4^{-1}_{a^2q^{20}}2_{aq^6}
+
2_{aq^6}2^{-1}_{-aq^8}2^{-1}_{-aq^{10}}1_{-aq^9}
\\+
3^{-1}_{a^2q^{14}}2_{aq^6}2_{-aq^6}1^{-1}_{-aq^7}1^{-1}_{-aq^{11}}
+
3^{-1}_{a^2q^{18}}2_{aq^6}2_{aq^8}1^{-1}_{aq^7}1^{-1}_{aq^9}
+2
3_{a^2q^{14}}3^{-1}_{a^2q^{18}}2^{-1}_{aq^8}1_{aq^7}
\\+
3^2_{a^2q^{14}}4^{-1}_{a^2q^{16}}2^{-1}_{aq^8}2^{-1}_{-aq^8}1^{-1}_{aq^9}
+
4_{a^2q^{12}}4^{-1}_{a^2q^{20}}2^{-1}_{aq^8}1_{aq^7}
+
3^{-1}_{a^2q^{18}}4_{a^2q^{12}}4_{a^2q^{16}}1^{-1}_{aq^9}
\\+
2^{-1}_{aq^8}2_{-aq^6}2^{-1}_{-aq^{10}}1_{aq^7}1^{-1}_{-aq^7}1_{-aq^9}
+
3^{-1}_{a^2q^{18}}2_{-aq^6}2_{-aq^8}1^{-1}_{aq^9}1^{-1}_{-aq^7}
+
2^{-1}_{-aq^{10}}1^{-1}_{aq^9}1_{-aq^5}1_{-aq^9}
\\+
2^{-1}_{aq^8}1_{aq^7}1_{-aq^5}1^{-1}_{-aq^{11}}
+
2^{-1}_{aq^{10}}1_{aq^5}
+
3^{-2}_{a^2q^{18}}4_{a^2q^{16}}2_{-aq^8}1_{aq^7}
+
3_{a^2q^{14}}4^{-1}_{a^2q^{16}}4^{-1}_{a^2q^{20}}2^{-1}_{aq^8}1_{aq^7}
\\+
3_{a^2q^{14}}2^{-1}_{aq^8}2^{-1}_{-aq^8}2^{-1}_{-aq^{10}}1_{aq^7}1_{-aq^9}
+
2^{-1}_{aq^8}2_{-aq^6}1_{aq^7}1^{-1}_{-aq^7}1^{-1}_{-aq^{11}}
+
2_{aq^6}2^{-1}_{-aq^8}1^{-1}_{-aq^{11}}
+
2_{aq^6}2^{-1}_{aq^{10}}1^{-1}_{aq^7}
\\+2
3_{a^2q^{14}}3^{-1}_{a^2q^{18}}1^{-1}_{aq^9}
+
4_{a^2q^{12}}4^{-1}_{a^2q^{20}}1^{-1}_{aq^9}
+
2_{-aq^6}2^{-1}_{-aq^{10}}1^{-1}_{aq^9}1^{-1}_{-aq^7}1_{-aq^9}
+
1^{-1}_{aq^9}1_{-aq^5}1^{-1}_{-aq^{11}}
+
3_{a^2q^{14}}2^{-1}_{aq^8}2^{-1}_{aq^{10}}
\\+
3^{-2}_{a^2q^{18}}4_{a^2q^{16}}2_{aq^8}2_{-aq^8}1^{-1}_{aq^9}
+
3_{a^2q^{14}}4^{-1}_{a^2q^{16}}4^{-1}_{a^2q^{20}}1^{-1}_{aq^9}
+
3_{a^2q^{14}}2^{-1}_{-aq^8}2^{-1}_{-aq^{10}}1^{-1}_{aq^9}1_{-aq^9}
\\+
2_{-aq^6}1^{-1}_{aq^9}1^{-1}_{-aq^7}1^{-1}_{-aq^{11}}
+
3^{-1}_{a^2q^{18}}4^{-1}_{a^2q^{20}}2_{-aq^8}1_{aq^7}
+
3^{-1}_{a^2q^{18}}4_{a^2q^{16}}2^{-1}_{-aq^{10}}1_{aq^7}1_{-aq^9}
\\+
3_{a^2q^{14}}2^{-1}_{aq^8}2^{-1}_{-aq^8}1_{aq^7}1^{-1}_{-aq^{11}}
+
3^{-1}_{a^2q^{18}}4_{a^2q^{16}}2^{-1}_{aq^{10}}2_{-aq^8}
+
3^{-1}_{a^2q^{18}}4^{-1}_{a^2q^{20}}2_{aq^8}2_{-aq^8}1^{-1}_{aq^9}
\\+
3^{-1}_{a^2q^{18}}4_{a^2q^{16}}2_{aq^8}2^{-1}_{-aq^{10}}1^{-1}_{aq^9}1_{-aq^9}
+
3_{a^2q^{14}}2^{-1}_{-aq^8}1^{-1}_{aq^9}1^{-1}_{-aq^{11}}
+
4^{-1}_{a^2q^{20}}2^{-1}_{-aq^{10}}1_{aq^7}1_{-aq^9}
\\+
3^{-1}_{a^2q^{18}}4_{a^2q^{16}}1_{aq^7}1^{-1}_{-aq^{11}}
+
4^{-1}_{a^2q^{20}}2^{-1}_{aq^{10}}2_{-aq^8}
+
4_{a^2q^{16}}2^{-1}_{aq^{10}}2^{-1}_{-aq^{10}}1_{-aq^9}
+
4^{-1}_{a^2q^{20}}2_{aq^8}2^{-1}_{-aq^{10}}1^{-1}_{aq^9}1_{-aq^9}
\\+
3^{-1}_{a^2q^{18}}4_{a^2q^{16}}2_{aq^8}1^{-1}_{aq^9}1^{-1}_{-aq^{11}}
+
4^{-1}_{a^2q^{20}}1_{aq^7}1^{-1}_{-aq^{11}}
+
3_{a^2q^{18}}4^{-1}_{a^2q^{20}}2^{-1}_{aq^{10}}2^{-1}_{-aq^{10}}1_{-aq^9}
+
4_{a^2q^{16}}2^{-1}_{aq^{10}}1^{-1}_{-aq^{11}}
\\+
4^{-1}_{a^2q^{20}}2_{aq^8}1^{-1}_{aq^9}1^{-1}_{-aq^{11}}
+
3^{-1}_{a^2q^{22}}1_{-aq^9}
+
3_{a^2q^{18}}4^{-1}_{a^2q^{20}}2^{-1}_{aq^{10}}1^{-1}_{-aq^{11}}
+
3^{-1}_{a^2q^{22}}2_{-aq^{10}}1^{-1}_{-aq^{11}}
+
2^{-1}_{-aq^{12}}$.}

\end{document}